\def\no{\if01}
\def\iftwelvept{\no}

\def\ifusepdf{\no}
\def\ifpsfont{\no}

\iftwelvept
\documentclass[leqno,12pt]{amsart}
\else
\documentclass[leqno,10pt]{amsart}
\fi
\usepackage{amssymb}
\usepackage{amscd}
\usepackage{latexsym}
\usepackage{verbatim}
\usepackage[all]{xy}

\ifusepdf
\usepackage{hyperref}
\else\fi
\ifpsfont
\usepackage[T1]{fontenc}
\usepackage{times}
\else\fi

\iftwelvept
\else\fi

\theoremstyle{plain}
\newtheorem{Theorem}{Theorem}[section]

\newtheorem{Proposition}[Theorem]{Proposition}
\newtheorem{Lemma}[Theorem]{Lemma}
\newtheorem{Corollary}[Theorem]{Corollary}

\theoremstyle{definition}

\newtheorem{Definition}[Theorem]{Definition}
\newtheorem{Remark}[Theorem]{Remark}
\newtheorem{Construction}[Theorem]{Construction}
\newtheorem{Example}[Theorem]{Example}

\newtheorem{Notation}[Theorem]{Notation}

\newcommand{\ZZ}{\mathbf{Z}}

\newcommand{\RR}{\mathbb{R}}
\newcommand{\RRR}{\mathbf{R}}

\newcommand{\DD}{\mathbb{D}}

\newcommand{\NNNN}{\operatorname{N}}
\newcommand{\GG}{\mathcal{G}}
\newcommand{\HH}{\operatorname{\mathcal{HH}}}

\newcommand{\HHH}{\operatorname{\mathcal{H}}_\bullet}
\newcommand{\HHHH}{\operatorname{\mathcal{H}}^\bullet}

\newcommand{\FF}{\mathcal{F}}

\newcommand{\AAA}{\mathcal{A}}

\newcommand{\MM}{\mathcal{M}}

\newcommand{\DDD}{\mathcal{D}}
\newcommand{\uni}{\mathbf{1}}
\newcommand{\CCC}{\mathcal{C}}

\newcommand{\PR}{\operatorname{Pr}^{\textup{L}}}

\newcommand{\MMM}{\mathcal{M}}
\newcommand{\NNN}{\mathcal{N}}

\newcommand{\GGG}{\mathcal{G}}
\newcommand{\PPP}{\mathcal{P}}

\newcommand{\Ker}{\operatorname{Ker}}

\newcommand{\Spec}{\operatorname{Spec}}

\newcommand{\Perf}{\operatorname{Perf}}

\newcommand{\SP}{\operatorname{Sp}}
\newcommand{\SPS}{\operatorname{Sp}^{\Sigma}}

\newcommand{\Mod}{\operatorname{Mod}}

\newcommand{\SSS}{\mathcal{S}}
\newcommand{\wSSS}{\widehat{\mathcal{S}}}

\newcommand{\colim}{\operatorname{colim}}

\newcommand{\Cat}{\textup{Cat}_{\infty}}

\newcommand{\Map}{\operatorname{Map}}

\newcommand{\Fun}{\operatorname{Fun}}
\newcommand{\Alg}{\operatorname{Alg}}

\newcommand{\End}{\operatorname{End}}

\newcommand{\FIN}{\operatorname{Fin}_\ast}

\newcommand{\wCat}{\widehat{\textup{Cat}}_{\infty}}

\newcommand{\CAlg}{\operatorname{CAlg}}
\newcommand{\CCAlg}{\operatorname{CAlg}^{\le0}}

\newcommand{\conn}{\textup{con}}

\newcommand{\BAR}{\operatorname{Bar}}

\newcommand{\EXT}{\operatorname{Art}^{\textup{tsz}}}
\newcommand{\TSZ}{\mathsf{TSZ}}

\newcommand{\DEFOR}{\mathcal{D}ef}
\newcommand{\DEFORM}{\mathcal{D}ef}

\newcommand{\FST}{\widehat{\mathsf{St}}^{\ast}}

\newcommand{\LM}{\operatorname{\mathcal{LM}}}
\newcommand{\RM}{\operatorname{\mathcal{RM}}}
\newcommand{\BM}{\operatorname{\mathcal{BM}}}

\newcommand{\Ind}{\operatorname{Ind}}

\newcommand{\PRST}{\PR_{\textup{St}}}
\newcommand{\Free}{\operatorname{Free}}

\newcommand{\assoc}{\operatorname{As}}

\newcommand{\eone}{\mathbf{E}_1}
\newcommand{\etwo}{\mathbf{E}_2}
\newcommand{\eenu}{\mathbf{E}_n}
\newcommand{\einf}{\mathbf{E}_\infty}

\newcommand{\KS}{\mathbf{KS}}
\newcommand{\Lie}{\mathbf{Lie}}

\newcommand{\LMod}{\operatorname{LMod}}
\newcommand{\RMod}{\operatorname{RMod}}
\newcommand{\BMod}{\operatorname{BMod}}
\newcommand{\ST}{\operatorname{\mathcal{S}t}}

\newcommand{\hhhh}{\mathfrak{h}}

\newcommand{\Def}{\operatorname{Def}}

\newcommand{\Proof}{{\sl Proof.}\quad}
\newcommand{\QED}{{\unskip\nobreak\hfil\penalty50\quad\null\nobreak\hfil
{$\Box$}\parfillskip0pt\finalhyphendemerits0\par\medskip}}

\begin{document}

\title{Moduli theory associated to Hochschild pairs}

\author{Isamu Iwanari}






\address{Mathematical Institute, Tohoku University, Sendai, Miyagi, 980-8578,
Japan}

\email{isamu.iwanari.a2@tohoku.ac.jp}

\begin{abstract}
Let $\CCC$ be an $A$-linear stable $\infty$-category and 
let $(\HH^\bullet(\CCC/A),\HH_\bullet(\CCC/A))$ be
the pair of the Hochschild cohomology spectrum (Hochschild cochain complex)
and the Hochschild homology spectrum (Hochschild chain complex).
The purpose of this paper is to establish a moduli-theoretic interpretation of the algebraic 
structure on the Hochschild pair $(\HH^\bullet(\CCC/A),\HH_\bullet(\CCC/A))$
of $\CCC$.
The notions of cyclic deformations and equivariant deformations (of the Hochschild chain complex)
associated to deformations of $\CCC$ play a central role.
\end{abstract}

\maketitle

\section{Introduction}

Since \cite{Ger} Gerstenhaber developed the deformation theory of associative algebras
by means of Hochschild cohomology \cite{Hoch}.
The Hochschild cohomology of an associative algebra governs the deformation
theory of algebra.
The relationship between deformations of associative algebras
and Hochschild cohomology has been fundamental and successful in many branches of mathematics: e.g.,
deformation quantizations of Poisson manifolds \cite{KP}.
The algebraic structure on the Hochschild cochain complex
(which computes Hochschild cohomology) plays an important role
in recent developments in  deformations of associative algebras
based on Hochschild cohomology. The Hochschild cochain complex
admits a structure of an algebra over the little disk operad
(see e.g. Introduction of \cite{KaSc} and references therein for many proofs
of the existence of an algebraic structure).
The Hochschild cohomology of associative algebras is invariant under Morita equivalences.
Based on the invariance, Hochschild cohomology is generalized to
those of abelian categories, suitable enriched categories
and stable $\infty$-categories, etc.
Moreover, deformation theories of abelian categories and stable $\infty$-categories
have been developed by using Hochschild cohomology theories of them \cite{LVB}, \cite{KL}, \cite[X]{DAG}.
In particular, in \cite{DAG}
the deformation theory of stable $\infty$-categories is formulated in the local version of derived geometry over
$\etwo$-algebras, i.e., algebras over the little disk operad $\etwo$, and it is shown that 
the deformation theory is controlled by the $\etwo$-algebra of Hochschild
cochain complex.

Let $A$ be a commutative ring spectrum (or simply we may suppose that $A$ is an ordinary commutative ring).
Let $\CCC$ be an $A$-linear small stable idempotent-complete  $(\infty,1)$-category/$\infty$-category.
Let $\HH^\bullet(\CCC/A)$ and $\HH_\bullet(\CCC/A)$
denote the Hochschild cohomology spectrum and the Hochschild homology spectrum of $\CCC$ over $A$,
respectively (in differential graded setting, they are
the Hochschild cochain complex and the Hochschild chain complex).
We refer to them simply as Hochschild cohomology and Hochschild homology, respectively
though these are not graded abelian groups.
The Hochschild homology $\HH_\bullet(\CCC/A)$ admits an $S^1$-action
given by the Connes operator.
Moreover, there exists a certain module action of the $\etwo$-algebra $\HH^\bullet(\CCC/A)$ on $\HH_\bullet(\CCC/A)$.
These algebraic structures on
the Hochschild pair $(\HH^\bullet(\CCC/A),\HH_\bullet(\CCC/A))$
is defined as an algebra over a topological operad $\KS$
having two colors, that is called Kontsevich-Soibelman
operad \cite{KS}.
For the construction of this algebra structure we refer the reader to \cite{I} and references therein
(we will use the construction in \cite{I}, cf. Section~\ref{KSalgebrasection} for a quick review).
The algebra $(\HH^\bullet(\CCC/A),\HH_\bullet(\CCC/A))$ over $\KS$
can be thought of as an analogue of the pair $(\oplus_{p\ge0} \wedge^p T_M, \oplus_{q\ge0}\Omega_{M}^q)$ of multivector fields and differential forms of a differential (or algebraic) manifold $M$.
Needless to say,
in algebraic geometry and differential geometry, the structures on $(\oplus_{p\ge0} \wedge^p T_M, \oplus_{q\ge0}\Omega_{M}^q)$ given by Cartan calculus are important. For example,
the algebraic description of local period maps for a family of algebraic manifolds due to Griffiths \cite{Gri}
has been fundamental in Torelli problems.

\vspace{1mm}

The purpose of this paper is to provide a moduli-theoretic interpretation of the algebraic 
structure on the Hochschild pair $(\HH^\bullet(\CCC/A),\HH_\bullet(\CCC/A))$.
To get a feeling of our results and moduli problems we consider,
we will describe
the formulation in an informal way.
For precise presentations, we refer the reader to the text.
Let $A$ be a connective noetherian commutative differential graded 
(dg) algebra over a field $k$ of characteristic zero.
Let $\CCC=\CCC_A$ be an $A$-linear small stable idempotent-complete
$\infty$-category.
We may regard $\CCC_A$ as a family of stable $\infty$-categories over $\Spec A$. 
We denote by $\ST_R$ the $\infty$-category of $R$-linear stable idempotent-complete
$\infty$-categories. 
For an augmented commutative dg $A$-algebra $R\to A$,
we let $\Def_{\CCC}(R)$ denote the $\infty$-groupoid of deformations of $\CCC$ to $R$ that consists of
$(\CCC_R \in \ST_R, \phi: \CCC_R\otimes_RA\simeq \CCC)$
where $\CCC_R\otimes_RA$ is the base change to $\ST_A$.
Let $\EXT_A$ be the $\infty$-category of trivial square-zero extensions of $A$.
The assignment $[R\to A]\mapsto \Def_{\CCC}(R)$
defines the deformation functor
\[
\Def_{\CCC}:\EXT_A\to \SSS
\]
associated to $\CCC$, where $\SSS$ is the $\infty$-category of $\infty$-groupoids/spaces.
Let $\HH_\bullet(\CCC/A)$ be the Hochschild homology of $\CCC$ that is a dg $A$-module. It admits an action of $S^1$, which corresponds to Connes operator.
Let $\Def^{S^1}(\HH_\bullet(\CCC/A)):\EXT_A\to \SSS$ be the deformation functor
which assigns $R\to A$ to the $\infty$-groupoid $\Def^{S^1}(\HH_\bullet(\CCC/A))(R)$ of $S^1$-equivariant deformations $(N_R, N\otimes_RA\simeq \HH_\bullet(\CCC/A))$ of $\HH_\bullet(\CCC/A)$ such that $N_R$ is a dg $R$-module endowed with $S^1$-action. 

Suppose that we are given a deformation $\CCC_R$ of $\CCC$ to $\ST_R$.
Then the relative Hochschild homology $\HH_\bullet(\CCC_R/R)$
is a deformation of $\HH_\bullet(\CCC/A)$ to $R$.
The moduli-theoretic assignment $\CCC_R\mapsto \HH_\bullet(\CCC_R/R)$
determines a natural transformation
\[
M_{\CCC}:\Def_{\CCC}\longrightarrow \Def^{S^1}(\HH_\bullet(\CCC/A)).
\]

There is a categorical correspondence between dg Lie algebras and pointed formal stacks
\cite[X]{DAG}, \cite{Gai2}, \cite{H}
 (cf. Section~\ref{formalstack}).
In particular, to a dg Lie algebra $L$ over $A$ one can associate
a pointed formal stack $\mathcal{F}_L$ that is defined as
a functor $\mathcal{F}_L:\EXT_A\to \SSS$.
This beautiful correspondence allows us to use dg Lie algebras in the study of deformation
theories, i.e., local moduli theory.
Let $\GG_{\CCC}$ be the dg Lie algebra over $A$ associated to the $\etwo$-algebra $\HH^\bullet(\CCC/A)$ whose underlying complex of $\GG_{\CCC}$ is $\HH^\bullet(\CCC/A)[1]$.
Let $\End^L(\HH_\bullet(\CCC/A))$ be the endomorphism dg Lie algebra
which has the $S^1$-action arising from that of $\HH_\bullet(\CCC/A)$.
Let $\End^L(\HH_\bullet(\CCC/A))^{S^1}$ denote the homotopy fixed points.
We associate formal stacks $\FF_{\GG_{\CCC}}$ and $\FF_{\End^L(\HH_\bullet(\CCC/A))^{S^1}}$.
There are canonical morphisms
\[
\Def_{\CCC}\to \FF_{\GG_{\CCC}}\ \ \  \textup{and}\ \ \  \Def^{S^1}(\HH_\bullet(\CCC/A))\to \FF_{\End^L(\HH_\bullet(\CCC/A))^{S^1}}
\]
in the functor category $\Fun(\EXT_A,\SSS)$.
Intuitively, they exhibit $\FF_{\GG_{\CCC}}$ and $\FF_{\End^L(\HH_\bullet(\CCC/A))^{S^1}}$
as formal stacks closest to $\Def_{\CCC}$ and $\Def^{S^1}(\HH_\bullet(\CCC/A))$, respectively \cite{DAG}. 
For instance,
each $\Def^{S^1}(\HH_\bullet(\CCC/A))(R)\to \FF_{\End^L(\HH_\bullet(\CCC/A))^{S^1}}(R)$ is fully faithful
when we regard spaces as $\infty$-groupoids.
Moreover,
$\Def^{S^1}(\HH_\bullet(\CCC/A))\to \FF_{\End^L(\HH_\bullet(\CCC/A))^{S^1}}$
is an equivalence under a mild condition. In this way, we may think that
$\GG_{\CCC}$ and $\End^L(\HH_\bullet(\CCC/A))^{S^1}$
controll deformations of $\CCC$ and $S^1$-equivariant deformations of $\HH_\bullet(\CCC/A)$.

The algebraic input is an action of the dg Lie algebra
$\GG_{\CCC}$ on $\HH_\bullet(\CCC/A)$.
This action is determined by the algebra $(\HH^\bullet(\CCC/A),\HH_\bullet(\CCC/A))$
over $\KS$ in a purely algebraic way.
The corresponding dg Lie algebra map
 $\GG_{\CCC}\to \End^L(\HH_\bullet(\CCC/A))^{S^1}$
gives rise to a morphism of formal stacks
\[
\mathfrak{A}_{\CCC}:\mathcal{F}_{\GG_{\CCC}}\longrightarrow \FF_{\End^L(\HH_\bullet(\CCC/A))^{S^1}}.
\]
In this situation, the simplest version of the main result can be stated as follows:

\begin{Theorem}
\label{intromain}
The diagram
\[
\xymatrix{
 \Def_{\CCC} \ar[r]^(0.4){M_{\CCC}} \ar[d] & \Def^{S^1}(\HH_\bullet(\CCC/A)) \ar[d] \\
 \FF_{\GG_{\CCC}} \ar[r]^(0.4){\mathfrak{A}_{\CCC}} & \FF_{\End^L(\HH_\bullet(\CCC/A))^{S^1}} \\
}
\]
commutes up to canonical homotopy.
\end{Theorem}
The morphism $M_{\CCC}:\Def_{\CCC} \to \Def^{S^1}(\HH_\bullet(\CCC/A))$
has a clear moduli-theoretic meaning.
On the other hand, 
$\FF_{\GG_{\CCC}} \to \FF_{\End^L(\HH_\bullet(\CCC/A))^{S^1}}$ is obtained by 
the $S^1$-equivariant Lie algebra action of $\GG_{\CCC}$ on $\HH_\bullet(\CCC/A)$
which comes from the algebraic structure on the Hochschild pair
$(\HH^\bullet(\CCC/A),\HH_\bullet(\CCC/A))$.
Therefore, we may think that the above result reveals a moduli-theoretic aspects of 
a certain portion of $(\HH^\bullet(\CCC/A),\HH_\bullet(\CCC/A))$.

\vspace{2mm}

{\it Cyclic deformations.}
An important observation is that
$M_{\CCC}$ factors as the sequence
\[
\xymatrix{
 \Def_{\CCC} \ar[r]^(0.4){M^\circlearrowleft_{\CCC}} &  \Def^{\circlearrowleft}(\HH_\bullet(\CCC/A))\ar[r]^{N_{\CCC}} &   \Def^{S^1}(\HH_\bullet(\CCC/A)) 
}
\]
such that
 $\Def^{\circlearrowleft}(\HH_\bullet(\CCC/A))$ are refinements of $\Def^{S^1}(\HH_\bullet(\CCC/A))$.
This factorization and $\Def^{\circlearrowleft}(\CCC)$ are crucial in this work though they do not appear in the statement of Theorem~\ref{intromain}.
Moreover, the refined diagram not only provides an appropriate framework which allows us to prove main
results but also
amplifies Theorem~\ref{intromain} with respect
to practical uses.
Recall that $\HH_\bullet(\CCC/A)$ has an $S^1$-action corresponding to Connes operator.
The functor $\Def^{S^1}(\HH_\bullet(\CCC/A))$
describes $S^1$-equivariant deformations of $\HH_\bullet(\CCC/A)$.
The functor $\Def^\circlearrowleft(\HH_\bullet(\CCC/A))$
describes what we call cyclic deformations of $\HH_\bullet(\CCC/A)$.
The notion of cyclic deformation
is the key to revealing the moduli-theoretic meaning, and 
it also has practical singnificance.
Since this notion is relatively new,
we briefly introduce the definition of cyclic deformations.
For an augmented $A$-algebra $R\to A$,
a cyclic deformation of $\HH_\bullet(\CCC/A)$ to $R$
is a pair $(N, N\otimes_{R\otimes S^1}A\simeq \HH_\bullet(\CCC/A))$
such that $N$ is a $(R\otimes_A S^1)$-module endowed
with an $S^1$-action which commutes with the canonical
action on $R\otimes_A S^1$.
 Here $R\otimes_A S^1$ is the tensor of $R$ by $S^1$ in the $\infty$-category of commutative
dg algebras over $A$ so that
 $R\otimes_A S^1\simeq R\otimes_{R\otimes_A R}R$.
Let $\Def^{\circlearrowleft}(\HH_\bullet(\CCC/A)):\EXT_A \to \SSS$ denote the deformation functor which assigns
to $R\in \EXT_A$ to the $\infty$-groupoid of
cyclic deformations of $\HH_\bullet(\CCC/A)$.
Each deformation of $\CCC_R$ of $\CCC$ in $\Def_{\CCC}(R)$ maps to
a cyclic deformation $\HH_\bullet(\CCC_R/A)$ of $\HH_\bullet(\CCC/A)$ 
that belongs to $\Def^{\circlearrowleft}(\CCC)(R)$ (it is important to notice that
the associated cyclic deformation is not $\HH_\bullet(\CCC_R/R)$ but $\HH_\bullet(\CCC_R/A)$).
The base change $\otimes_{R\otimes_AS^1}R$ sends
$\HH_\bullet(\CCC_R/A)$ to $\HH_\bullet(\CCC_R/R)$, that is
an $S^1$-equivariant deformation of $\HH_\bullet(\CCC/A)$ to $R$.
The transition
\[
\CCC_R \rightsquigarrow  \{\HH_\bullet(\CCC_R/A) \curvearrowleft S^1\} \rightsquigarrow \{\HH_\bullet(\CCC_R/R) \curvearrowleft S^1\}
\]
induces the above factorization (the final procedure forgets the $S^1$-action).

The functor $\Def^{\circlearrowleft}(\HH_\bullet(\CCC/A))$
has a natural transfomation to a functor $\FF_{A\oplus\End(\HH_\bullet(\CCC/A))}^{\circlearrowleft}$
that is defined in an algebraic way
(cf. Notation~\ref{cyclicliefunctor}, Construction~\ref{fromcube}).
We will give several different presentations of
$\mathcal{F}^{\circlearrowleft}_{A\oplus \End(\HH_\bullet(\CCC/A))}$
which are related to one another via Koszul-type dualities.
These different descriptions enable us to prove Theorem~\ref{intromain}
and the following generalization (cf. Theorem~\ref{supermain}, Remark~\ref{refinedremark}):

\begin{Theorem}
\label{intromain2}
There exists a diagram in the functor category $\Fun(\EXT_A,\SSS)$
\[
\xymatrix{
\Def_{\CCC}  \ar[r]^(0.35){M^{\circlearrowleft}_{\CCC}} \ar[d]_{J_{\CCC}} & \Def^{\circlearrowleft}(\HH_\bullet(\CCC/A)) \ar[d]^{J^{\circlearrowleft}_{\HH_\bullet(\CCC/A)}}  \ar[r]^{N_{\CCC}}  &   \Def(\HH_\bullet(\CCC/A))^{S^1} \ar[d]^{J^{S^1}_{\HH_\bullet(\CCC/A)}}   \\
\mathcal{F}_{\GG_{\CCC}} \ar[r]^(0.3){\mathfrak{T}_{\CCC}}  &   \mathcal{F}_{A\oplus \End(\HH_\bullet(\CCC/A))}^{\circlearrowleft} \ar[r]   &  \mathcal{F}_{\End^L(\HH_\bullet(\CCC/A))^{S^1}} 
}
\]
which commutes and is an extension of the diagram in Theorem~\ref{intromain}.
\end{Theorem}
It should be worth emphasizing that the upper row has a moduli-theoretic interpretation
while the lower row is defined in a purely algebraic way.
In fact, the lower row admits a description in terms of dg Lie algebras.
See Propotion~\ref{algebraint2}, Proposition~\ref{algebraint3}, Remark~\ref{modulimeaning1},
and Remark~\ref{modulimeaning2}.
Let $\GG_{\CCC}^{S^1}$ denote the cotensor by $S^1$, i.e., $\GG_{\CCC}\simeq \GG_{\CCC}\times_{\GG_{\CCC}\times\GG_{\CCC}}\GG_{\CCC}$.
The morphism 
$\mathfrak{T}_{\CCC}:\mathcal{F}_{\GG_{\CCC}} \to  \mathcal{F}_{A\oplus \End(\HH_\bullet(\CCC/A))}^{\circlearrowleft}$
is induced by the Lie algebra action $\GG_{\CCC}^{S^1}\to \End^L(\HH_\bullet(\CCC/A))$
determined by the Hochschild pair $(\HH^\bullet(\CCC/A), \HH_\bullet(\CCC/A))$, which extends the action of $\GG_{\CCC}$ on $\HH_\bullet(\CCC/A)$.
Informally, Theorem~\ref{intromain2} especially means that 
a moduli-theoretic interpretation of
the Lie algebra action of $\GG_{\CCC}^{S^1}$ on $\HH_\bullet(\CCC/A)$
is given by
$M^{\circlearrowleft}_{\CCC}:\Def_{\CCC}\to  \Def^{\circlearrowleft}(\HH_\bullet(\CCC/A))$
determined by the assignment $\CCC_R\rightsquigarrow \HH_\bullet(\CCC_R/A)$.
Likewise, an algebraic incarnation of $\Def_{\CCC}\to  \Def^{S^1}(\HH_\bullet(\CCC/A))$
is the composite $\mathcal{F}_{\GG_{\CCC}} \to  \mathcal{F}_{\End^L(\HH_\bullet(\CCC/A))^{S^1}}$
induced by the $S^1$-equivariant Lie algebra map $\GG_{\CCC}\stackrel{\textup{diagonal}}{\longrightarrow} \GG_{\CCC}^{S^1} \to \End^L(\HH_\bullet(\CCC/A))$.
The situation may be depicted as the following table:
\begin{center}
\begin{tabular}{c|c}
algebra & moduli \\ \hline
$\GG_{\CCC}\curvearrowright \HH_\bullet(\CCC/A)$ & $\CCC_R\rightsquigarrow  \HH_\bullet(\CCC_R/R)$: equivariant deformation \\
$\GG^{S^1}_{\CCC}\curvearrowright \HH_\bullet(\CCC/A)$ & $\CCC_R\rightsquigarrow  \HH_\bullet(\CCC_R/A)$: cyclic deformation \\
\end{tabular}
\end{center}

\vspace{3mm}

{\it Applications.}
Let us mention what are uses of our results.
One use is to applications to the study of Hochschild homology and Hochschild cohomology
of stable $\infty$-categories.
The relationship between moduli-theoretic methods in derived geometry
and algebraic structures of Hochschild pair is useful.
In fact, our main result has already been applied to the study of Hochschild homology:
in \cite{IKS} we use the pair $(\HH^\bullet(\CCC/A),\HH_\bullet(\CCC/A))$
and main results of this paper to prove that $\HH_\bullet(\CCC/A)$ admits an equivariant 
deformation to the derived loop space $LX=\Map(S^1,\Spec A)$ in a natural way.
Moreover, we show that the associated periodic cyclic homology/complex
$\mathcal{HP}_\bullet(\CCC/A)$ has a $D$-module structure.
It is worth mentioning that the data of cyclic deformations plays an essential role
in these applications.
The results of this paper will be used in our future works.

Another use is to the study of the period map in noncommutative algebraic geometry.
The negative cyclic homology gives a Hodge-type filtration on the periodic cyclic homology 
of a stable $\infty$-category (cf. \cite{IP}).
Theorem~\ref{intromain2} can be used to construct and study 
the period map of a family of stable $\infty$-categorie 
which is a map into the classifying (moduli) space of Hodge-type filtrations 
on the periodic cyclic homology.
Combining with the recent progress on the study of $\HH_\bullet(\CCC/A)$, one may expect
applications of algebraic structures on $\HH_\bullet(\CCC/A)$ to the period map in noncommutative
algebraic
geometry. 

\vspace{2mm}

This paper is organized as follows: Section 2 collects conventions and some of the notation that we will
use. In Section 3, we will review
some of background material
and will formulate adequate basis for our main goal.
In Section 4 we will give a brief guide to subsequent sections for the reader's convenience.
In Section 5, we give an axiomatic formulation of deformation problems and
their Koszul duals.
In Section 6 we apply the axiomatic formulation defined in Section 5
to deformations of categories and cyclic deformations of modules.
The main result is Proposition~\ref{firstkoszuldual}.
In Section 7 we give detail analysis of the interaction between Koszul
duality with changes of operads and with Hochschild homology.
In Section 8 we prove main results of this paper.

\section{Notation and Convention}

{\it $(\infty,1)$-categories.}
Throughout this paper we use 
the language of $(\infty,1)$-categories.
We use the theory of {\it quasi-categories} as a model of $(\infty,1)$-categories.
We assume that the reader is familiar with this theory.
We will use the notation similar to that used in \cite{I}, \cite{IKS}.
A quasi-category is a simplicial set which
satisfies the weak Kan condition of Boardman-Vogt.
Following \cite{HTT}, we shall refer to quasi-categories
as {\it $\infty$-categories}.
Our main references are \cite{HTT}
 and \cite{HA}.
 To an ordinary category, we can assign an $\infty$-category by taking
its nerve, and therefore
when we treat ordinary categories we often omit the nerve $\NNNN(-)$
and directly regard them as $\infty$-categories.
 
 Here is a list of (some) of the conventions and notation that we will use:

\begin{itemize}

\item $\SSS$: the $\infty$-category of (small) spaces.

\item $\ZZ$: the ring of integers, $\RRR$ denotes the set of real numbers which we regard as either a topological space or a ring.

\item $\Delta$: the category of linearly ordered non-empty
finite sets (consisting of $[0], [1], \ldots, [n]=\{0,\ldots,n\}, \ldots$)

\item $\Delta^n$: the standard $n$-simplex

\item $\textup{N}$: the simplicial nerve functor (cf. \cite[1.1.5]{HTT})

\item $\SSS$: $\infty$-category of small spaces/$\infty$-groupoids. We denote by $\widehat{\SSS}$
the $\infty$-category of large spaces (cf. \cite[1.2.16]{HTT}).

\item $\CCC^\simeq$: the largest Kan subcomplex of an $\infty$-category $\CCC$

\item $\CCC^{op}$: the opposite $\infty$-category of an $\infty$-category. We also use the superscript ``op" to indicate the opposite category for ordinary categories and enriched categories.


\item $\SP$: the stable $\infty$-category of spectra.

\item $\Fun(A,B)$: the function complex for simplicial sets $A$ and $B$. If $A$ and $B$ are $\infty$-categories, we regard $\Fun(A,B)$ as the functor category.

\item $\Map_{\mathcal{C}}(C,C')$: the mapping space from an object $C\in\mathcal{C}$ to $C'\in \mathcal{C}$ where $\mathcal{C}$ is an $\infty$-category.
We usually view it as an object in $\mathcal{S}$ (cf. \cite[1.2.2]{HTT}).

\item $\FIN$: the category of pointed finite
sets $\langle 0 \rangle, \langle 1 \rangle,\ldots \langle n \rangle,...$
where $\langle n \rangle=\{*,1,\ldots n \}$ with the base point $*$.
We write $\Gamma$ for $\NNNN(\FIN)$. $\langle n\rangle^\circ=\langle n\rangle\backslash*$.
Notice that the (nerve of) Segal's gamma category is the opposite category
of our $\Gamma$.

\end{itemize}

{\it Operads and Algebras.}
We will use operads.
We employ
the theory of {\it $\infty$-operads} which is
thoroughly developed in \cite{HA}.
The notion of $\infty$-operads gives
one of the models of colored operads.
Here is a list of (some) of the notation about $\infty$-operads and algebras over them that we will use:

\begin{itemize}

\item Let $\mathcal{M}^\otimes\to \mathcal{O}^\otimes$ be a fibration of $\infty$-operads. We denote by
$\Alg_{/\mathcal{O}^\otimes}(\mathcal{M}^\otimes)$ the $\infty$-category of algebra objects (cf. \cite[2.1.3.1]{HA}).  We often write
$\Alg_{\mathcal{O}^\otimes}(\MMM)$ for $\Alg_{/\mathcal{O}^\otimes}(\MMM^\otimes)$.

\item $\CAlg(\mathcal{M}^\otimes)$: $\infty$-category of commutative
algebra objects in a symmetric
monoidal $\infty$-category $\mathcal{M}^\otimes\to \NNNN(\FIN)=\Gamma$.
When the symmetric monoidal structure is clear,
we usually write $\CAlg(\mathcal{M})$ for $\CAlg(\mathcal{M}^\otimes)$.

\item $\Mod_R^\otimes(\mathcal{M}^\otimes)\to \Gamma$: symmetric monoidal
$\infty$-category of
$R$-module objects,
where $\mathcal{M}^\otimes$
is a symmetric monoidal $\infty$-category such that (1)
the underlying $\infty$-category admits a colimit for any simplicial diagram, and (2)
its tensor product functor $\mathcal{M}\times\mathcal{M}\to \mathcal{M}$
preserves
colimits of simplicial diagrams separately in each variable.
Here $R$ belongs to $\CAlg(\mathcal{M}^\otimes)$
cf. \cite[3.3.3, 4.5.2]{HA}.
We write $\Mod_R(\mathcal{M}^\otimes)$ for the underlying $\infty$-category.
When $\MMM^\otimes$ is the symmetric monoidal $\infty$-category $\SP^\otimes$
of spectra, 
we write $\Mod_R^\otimes$
for
$\Mod_R^\otimes(\SP^\otimes)$.

\item $\CAlg_R$: $\infty$-category of commutative
algebra objects in the symmetric monoidal $\infty$-category $\Mod_R^\otimes$
where $R$ is a commutative ring spectrum, that is, an object of $\CAlg(\SP)$.
We write $\CAlg^+_R$ for $(\CAlg_R)_{/R}\simeq \CAlg(\SP)_{R//R}$.
When $R$ is the Eilenberg-MacLane spectrum $HC$ with a commutative ring $C$,
then we write $\CAlg_C$ for $\CAlg_{HC}$.
If $F$ is a field of characteristic zero,
the $\infty$-category $\CAlg_F$ is
equivalent to the $\infty$-category
obtained from the model category of commutative
differential graded $F$-algebras by inverting quasi-isomorphisms
(cf. \cite[7.1.4.11]{HA}). Therefore,
we often think of a commutative differential graded (dg) algebra
as an object of $\CAlg_F$ and refer to an object of $\CAlg_F$ as a commutative dg algebra.
We denote by $\CCAlg_{F}$ is the full subcategory of $\CAlg_F$ 
which consists of connective commutative dg algebras.
For $A\in \CCAlg_{F}$, we say that $A$ is noetherian if $H^0(A)$
is a usual
noetherian ring and if $H^{-n}(A)$ is trivial when $n>>0$ and is of finite type over $H^0(A)$ for any $n$.

\item $\mathbf{E}^\otimes_n$: the $\infty$-operad of little $n$-cubes (cf. \cite[5.1]{HA}).
For a symmetric monoidal $\infty$-category $\CCC^\otimes$,
we write $\Alg_{n}(\CCC)$ or $\Alg_{\eenu}(\CCC)$ for
the $\infty$-category of algebra objects over $\mathbf{E}^\otimes_n$ in $\CCC^\otimes$.
We refer to an object of $\Alg_{n}(\CCC)$ as an $\eenu$-algebra in $\CCC$.
If we denote by $\assoc^\otimes$ the associative operad (\cite[4.1.1]{HA}),
there is the standard equivalence $\assoc^\otimes\simeq \eone^\otimes$
of $\infty$-operads. We usually identify $\Alg_{1}(\CCC)$
with the $\infty$-category $\Alg_{\assoc^\otimes}(\CCC)$, that is, the 
$\infty$-category of associative algebras in $\CCC$.
We write $\Alg_{n}^+(\CCC)$ for the $\infty$-category $\Alg_n(\CCC)_{/\uni_{\CCC}}$
of augmeneted objects where $\uni_{\CCC}$ denotes the unit algebra.

\item $\LM^\otimes$: the $\infty$-operad defined in \cite[4.2.1.7]{HA}.
Roughly, an algebra over $\LM^\otimes$ is a pair $(A,M)$
such that $A$ is an unital associative algebra
and $M$ is a left $A$-module.
For a symmetric monoidal $\infty$-category $\CCC^\otimes \to \Gamma$,
we write $\LMod(\CCC^\otimes)$ or $\LMod(\CCC)$
for $\Alg_{\LM^\otimes}(\CCC^\otimes)$.
There is the natural inclusion of $\infty$-operads $\assoc^\otimes\to \LM^\otimes$.
This inclusion determines $\LMod(\CCC)\to \Alg_{\assoc^\otimes}(\CCC)\simeq \Alg_1(\CCC)$
which sends $(A,M)$ to $A$.
For $A\in \Alg_1(\CCC)$, we define $\LMod_A(\CCC)$
to be the fiber of $\LMod(\CCC)\to \Alg_1(\CCC)$ over $A$ in $\Cat$.

\item $\RM^\otimes$, $\BM^\otimes$: these $\infty$-operads are variants of $\LM^\otimes$
which are used to define structures of right modules over associative algebras
and bimodules over pairs of algebras \cite[4.2.1.36, 4.3.1]{HA}.
Informally, an algebra over $\RM^\otimes$ is a pair $(A,M)$
such that $A$ is an unital associative algebra and
$M$ is a right $A$-module.
In a similar vein, an algebra over $\BM^\otimes$ is a triple $(A,M,B)$
such that $A$ and $B$ are unital associative algebras,
and $M$ is an $A$-$B$-bimodule.
For a symmetric monoidal $\infty$-category $\CCC^\otimes \to \Gamma$,
we write $\RMod(\CCC^\otimes)$ (or simply $\RMod(\CCC)$) and $\BMod(\CCC^\otimes)$
(or $\BMod(\CCC)$)
for $\Alg_{\RM^\otimes}(\CCC^\otimes)$ and $\Alg_{\BM^\otimes}(\CCC^\otimes)$,
respectively.
There is a canonical functor $\BMod(\CCC^\otimes)\to \Alg_1(\CCC^\otimes)\times \Alg_1(\CCC^\otimes)$ which sends $(A,M,B)$ to $(A,B)\in \Alg_1(\CCC^\otimes)\times \Alg_1(\CCC^\otimes)$.

\item Unless otherwise stated, $k$ is a base field of characteristic zero.

\end{itemize}

{\it Group actions.}
Let $G$ be a group object in $\SSS$ (see e.g. \cite[7.2.2.1]{HTT} for the notion of group objects).
The main example in this paper is the circle $S^1$.
Let $\mathcal{C}$ be an $\infty$-category.
For an object $C\in \CCC$,
a $G$-action on $C$ means a lift of $C\in \CCC$
to $\Fun(BG, \CCC)$, where $BG$ is the classifying space of $G$.
A $G$-equivariant morphism means a morphism in $\Fun(BG, \CCC)$.
We often identify $\Fun(BG, \CCC)$ as the limit (``$G$-invariants'') of the trivial $G$-action on $\CCC$
and write $\CCC^G$ for $\Fun(BG, \CCC)$ (e.g., $\Mod_A^{S^1}$).
(We remark that when we regard $G$ as a group object, $\CCC^G$ is not the cotensor with the space
$G$.)

\section{Background}

\label{sectionbackground}

In this section, we review some theories which we will use.

\subsection{$\infty$-categories of stable $\infty$-categories}
\label{sectionmodule}
Let $\ST$ be the $\infty$-category of small
stable idempotent-complete $\infty$-categories
whose morphisms are exact functors.
This $\infty$-category is compactly generated.
Let $\CCC$ be a small
stable idempotent-complete $\infty$-category
and let 
$\Ind(\CCC)$ denote the $\infty$-category of Ind-objects \cite[5.3.5]{HTT}.
Then $\Ind(\CCC)$ is a compactly generated stable $\infty$-category.
The inclusion $\CCC\to \Ind(\CCC)$ identifies
the essential image with the full subcategory $\Ind(\CCC)^\omega$ spanned by
compact obejcts in $\Ind(\CCC)$.
Given $\CCC,\CCC'\in \ST$, if
we write $\Fun^{\textup{ex}}(\CCC,\CCC')$ for the full subcategory
spanned by exact functors,
the left Kan extension \cite[5.3.5.10]{HTT} gives rise to
a fully faithful functor $\Fun^{\textup{ex}}(\CCC,\CCC')\to \Fun^{\textup{L}}(\Ind(\CCC),\Ind(\CCC'))$ whose essential image consists of
those functors $F$ such that the essential image of $F$ are contained in $\CCC'$.
Here $\Fun^{\textup{L}}(\Ind(\CCC),\Ind(\CCC'))\subset \Fun(\Ind(\CCC),\Ind(\CCC'))$ denotes the full subcategory
consisting of those functors which preserve small colimits.
We let $\PR$ denote the $\infty$-category of presentable $\infty$-categories 
such that mapping spaces are spaces of functors which preserve small colimits
(i.e., left adjoint functors) \cite[5.5.3]{HTT}. It has a closed symmetric monoidal structure
whose internal Hom/mapping objects are given by $\Fun^{\textup{L}}(-,-)$,
see \cite[4.8.1.15]{HA}.
We set $\PRST=\Mod_{\SP^\otimes}(\PR)$,
which can be regarded as the full subcategory of $\PR$
that consists of stable presentable $\infty$-categories.
We denote by $\textup{Cgt}_{\textup{St}}^{\textup{L,cpt}}\subset \PR_{\textup{St}}$ is the subcategory spanned by compactly generated
stable $\infty$-categories, whose mapping spaces consist of those functors which
preserve small colimits and preserve compact objects.
There exists a sequence
\[
\ST\to \textup{Cgt}_{\textup{St}}^{\textup{L,cpt}} \subset \PR_{\textup{St}}
\]
where
the left arrow is an equivalence given by $\Ind$-construction
$\CCC\mapsto \Ind(\CCC)$.
The subcategory $\textup{Cgt}_{\textup{St}}^{\textup{L,cpt}} \subset \PR_{\textup{St}}$ is closed under
the tensor product so that 
$\textup{Cgt}_{\textup{St}}^{\textup{L,cpt}}$ and $\ST$
inherit symmetric monoidal structures from that of $\PRST$.
The stable $\infty$-category of compact spectra is a unit object in $\ST$.
Given two objects $\CCC$ and $\CCC'$ of $\ST$, the tensor product
$\CCC\otimes \CCC'$ is naturally equivalent to
the full subcategory
$(\Ind(\CCC)\otimes \Ind(\CCC'))^{\omega}\subset \Ind(\CCC)\otimes \Ind(\CCC')$ spanned by compact objects. The tensor product functor $\otimes :\ST\times \ST\to \ST$
preserves small colimits separately in each variable since $\Fun^{\textup{ex}}(-,-)$ gives us internal
Hom/mapping objects.
For more details, we refer the readers to
\cite[Section 3]{BGT1}, \cite[4.8]{HA}.

Let $A\in \CAlg(\SP)$.
Let $\Mod^\otimes_A\in \CAlg(\PRST)$ be the symmetric monoidal $\infty$-category of $A$-modules in $\SP$.
Let $\Perf^\otimes_A\in \CAlg(\ST)$ be the symmetric monoidal $\infty$-category of compact $A$-modules in $\SP$.
Let $\Mod^\otimes_{\Perf_A^\otimes}(\ST)$ be
the symmetric monoidal $\infty$-category of $\Perf_A^\otimes$-modules in $\ST$.
This symmetric monoidal $\infty$-category is presentable, and
the tensor product functor preserves small colimits separately in each variable. 
We refer to an object of the underlying $\infty$-category $\Mod_{\Perf_A^\otimes}(\ST)$
as an $A$-linear small stable $\infty$-category.
For ease of notation, put $\ST_A^\otimes=\Mod^\otimes_{\Perf_A^\otimes}(\ST)$
and $\ST_A:=\Mod_{\Perf_A^\otimes}(\ST)$.
We refer the reader to \cite{I} for the description of $\ST_A^\otimes$ by means of spectral
categories.
We write $(\PR_A)^{\otimes}$ for $\Mod_{\Mod_A^\otimes}^\otimes(\PRST)$.
We refer to an object of the underlying $\infty$-category $\PR_A$
as an $A$-linear stable presentable $\infty$-categories.

For $B\in \Alg_{1}(\Mod_A)$, we denote by $\LMod_B(\Mod_A)$ (resp. $\RMod_A(\Mod_A)$)
(or simply $\LMod_B$ and (resp. $\RMod_B$)) the $\infty$-category of left $B$-modules
(resp. right $B$-module spectra) (thare is a canonical equivalence $\LMod_B(\Mod_A)\simeq \LMod_B(\SP)$ induced by the forgetful functor $\Mod_A\to \SP$).
There is a symmetric monoidal functor $\Theta_A:\Alg_{1}(\Mod_A)\to (\PR_A)_{\Mod_A/}$
which carries $B$ to $\alpha:\Mod_A\to \LMod_{B}(\Mod_A)$
such that $\alpha$ is determined by $\SP\stackrel{\otimes A}{\to} \Mod_A\to \LMod_{B}(\Mod_A)$ which carries the sphere spectrum
to $B$. The homomorphism $f:B\to C$ maps to the base change functor $f^*:\LMod_B\to \LMod_C$.
It is a fully faithful left adjoint functor, see \cite[4.8.5.11]{HA}.
Since $\Theta_A$ is a symmetric monoidal functor, 
by Dunn additivity theorem \cite[5.1.2]{HA}, it determines
$\Alg_{n}(\Mod_A)\simeq \Alg_{n-1}(\Alg_{1}(\Mod_A))\to  \Alg_{n-1}((\PR_A)_{\Mod_A/})\simeq \Alg_{n-1}(\PR_A)$.
In particular, for $B\in \Alg_{2}(\Mod_A)$, $\LMod_B$ belongs to $\Alg_{1}(\PR_A)$,
that is, an associative monoidal $\infty$-category.
For $B\in \Alg_{1}(\Mod_A)$, we denote by $\Perf_{B}$
the smallest stable idempotent-complete subcategory of $\LMod_B$ which contains $B$
(it is more appropriate to write $\textup{LPerf}_B$ instead of $\Perf_B$).
The full subcategory $\Perf_B$ coincides with the full subcategory of compact objects
so that $\Perf_B\to \LMod_B$ is extended to an equivalence $\Ind(\Perf_B)\stackrel{\sim}{\to} \LMod_B$.
Moreover, if $B$ is a commutative algebra (which comes from $\CAlg(\Mod_A)$), $\Perf_B$
can be identified with the full subcategory of dualizable objects in the symmetric monoidal $\infty$-category
$\Mod_B\simeq \LMod_B$. 
When $B\in \Alg_{1}(\Mod_A)$ is promoted to $\Alg_{2}(\Mod_A)$,
the tensor product functor $\otimes:\LMod_B\times \LMod_B\to \LMod_B$
induces $\Perf_B\times \Perf_B\to \Perf_B$ so that $\Perf_B$ inherits a monoidal
structure from that of $\LMod_B$.
The functor $\Theta_{A}:\Alg_{1}(\Mod_A)\to (\PR_A)_{\Mod_A/}$
factors as $\Alg_{1}(\Mod_A) \to \LMod_{\Mod_A^\otimes}(\textup{Cgt}_{\textup{St}}^{\textup{L,cpt}})_{\Mod_A/} \to (\PR_A)_{\Mod_A/}$.
It gives rise to a symmetric monoidal functor
$\Alg_{1}(\Mod_A)\to (\ST_A)_{\Perf_{A}/}\simeq   \LMod_{\Mod_A^\otimes}(\textup{Cgt}_{\textup{St}}^{\textup{L,cpt}})_{\Mod_A/} $.
This functor sends $B$ to $\Perf_A\to \Perf_B$.
We denote the symmetric monoidal
functor $\Alg_{1}(\Mod_A)\to \ST_A$ (induced by forgetting the functor from $\Perf_{A}$)
by $\Perf_{(-)}$.


\subsection{}
\label{stablederived}
Let $\AAA$ be a small $\infty$-category.
Let $\PPP(\AAA)$ denote the functor category $\Fun(\AAA^{op},\SSS)$ where the $\SSS$ is the $\infty$-category of spaces/$\infty$-groupoids. There
is the Yoneda embedding $\mathfrak{h}_{\AAA}:\AAA\to \PPP(\AAA)$.
Let $\PPP_{\Sigma}(\AAA)\subset \PPP(\AAA)$ be the full subcategory spanned by
those functors $\AAA^{op}\to \SSS$ which preserve finite products \cite[5.5.8.8]{HTT}.
The $\infty$-category  $\PPP_{\Sigma}(\AAA)$ is compactly generated,
and $\PPP_{\Sigma}(\AAA)\subset \PPP(\AAA)$ is characterized as
the smallest full subcategory which contains
the essential image of the Yoneda embedding and is closed under sifted colimits.

Suppose that $\AAA$ admits finite coproducts and a zero object $0$.
Consider the set of morphisms $S=\{\mathfrak{h}_{\AAA}(0)\sqcup_{\mathfrak{h}_{\AAA}(M)}\mathfrak{h}
_{\AAA}(0)\to \mathfrak{h}_{\AAA}(0\sqcup_{M}0)\}$
where $M\in \AAA$
such that the pushout $0\sqcup_{M}0$ exists in $\AAA$.
Let $\mathcal{P}^{st}_{\Sigma}(\AAA)$ be
the presentable $\infty$-category obtained from $\mathcal{P}_{\Sigma}(\AAA)$
by inverting morphisms in $S$ (see e.g. \cite[5.5.4]{HTT} for the localization).
The $\infty$-category
$\mathcal{P}^{st}_{\Sigma}(\AAA)$ can be regarded as the full subcategory of $\PPP_\Sigma(\AAA)$
spanned by $S$-local objects. In other words, $\mathcal{P}^{st}_{\Sigma}(\AAA)$ is the full subcategory
which consists of functors $F$ such that the canonical morphism
$F(0\sqcup_M0)\to \ast\times_{F(M)}\ast$ is an equivalence for any $M\in \AAA$ such that
$0\sqcup_M0$ exists in $\AAA$.
Note that any object of the essential image of the Yoneda embedding
$\AAA\to \PPP_\Sigma(\AAA)$ is $S$-local.
Let $\DDD$ be a presentable $\infty$-category.
Let $\Fun^{\textup{L}}(\mathcal{P}^{st}_{\Sigma}(\AAA),\DDD)$ be the full subcategory
of $\Fun(\mathcal{P}^{st}_{\Sigma}(\AAA),\DDD)$ which consists of colimit-preserving functors (i.e.,
left adjoint functors).
Let $\Fun^{st}(\AAA,\DDD)$ be the full subcategory of $\Fun(\AAA,\DDD)$
spanned by those functors $f$ which preserve finite coproducts 
and carry pushouts of the form $0\sqcup_M0$ ($M\in \AAA$)
to $f(0)\sqcup_{f(M)}f(0)$.
Taking into account the universal properties of $\PPP_\Sigma$
and the localization  \cite[5.5.8.15, 5.5.4.20]{HA},
we see that the composition with the fully faithful functor $\AAA\hookrightarrow \PPP^{st}_\Sigma(\AAA)$ induced by the
Yoneda embedding determines an equivalence
\[
\Fun^{\textup{L}}(\PPP^{st}_\Sigma(\AAA),\DDD)\stackrel{\sim}{\to} \Fun^{st}(\AAA,\DDD).
\]

\subsection{Formal stacks}
\label{formalstack}
Let $A$ be a connective commutative dg algebra $A$ over a field $k$ of characteristic zero.
We use a correspondence between formal stacks and dg Lie algebras over $A$, that is proved by Hennion \cite{H}. In Gaitsgory and Rozenblyum \cite{Gai2}, a similar correspondence is established in the Ind-coherent setting. 
These are generalizations of the correspondence between dg Lie algebras and formal moduli problems
over a field of characteristic zero in Lurie \cite[X]{DAG}. 

Let $Lie_A$ be the $\infty$-category of dg Lie algebras: there are several approaches to define it.
One approach is to obtain it from the model category of dg Lie algebras by inverting quasi-isomorphisms.
Another one is to consider the $\infty$-category of algebras in $\Mod_A$ over the Lie operad $\Lie$.
Let
$\EXT_A$ be the full subcategory of $\CCAlg_{A//A}:=(\CCAlg_A)_{/A}$,
which is spanned by  trivial square zero extensions $A=A\oplus 0\hookrightarrow A\oplus M\stackrel{p_1}{\to} A$ such that
$M$ is a connective $A$-module of the form $\oplus_{1\le i\le n}A^{\oplus r_i}[d_i]$.
By abuse of notation, we often write $R$ for an object $A\to R\to A$ of $\CCAlg_{A//A}$.
Similarly, we often omit the augmentations from the notation.
Let $\TSZ_A$ denote the opposite category of $\EXT_A$.
We define the $\infty$-category $\FST_A$ of formal stacks over $A$ to be
$\PPP_\Sigma^{st}(\TSZ)$. We often regard $\TSZ_A$ as a full subcategory of $\FST_A$.
By definition, $\FST_A$ is the full subcategory of $\Fun(\EXT_A,\SSS)$
so that we think of a formal stack as a functor $\EXT_A\to \SSS$.

Let $\Free_{Lie}:\Mod_A\to Lie_A$ be the free Lie algebra functor
which is a left adjoint to the forgetful functor $Lie_A\to \Mod_A$. 
Let $\Mod_{A}^{f}\subset \Mod_A$ be the full subcategory that consists of
objects of the form $\oplus_{1\le i\le n}A^{\oplus r_i}[d_i]$ ($d_i\le -1$).
Let $Lie_A^f$ be the full subcategoy of $Lie_A$, which is the essential image of the restriction of
the free Lie algebra functor  $\Mod_A^{f}\to Lie_A$.
According to \cite[1.2.2]{H}, the inclusions $\Mod_A^f\hookrightarrow \Mod_A$ and 
$Lie_A^f\hookrightarrow Lie_A$ are extended to equivalences
$\PPP_\Sigma^{st}(\Mod_A^f) \stackrel{\sim}{\to} \Mod_A$ and $\PPP_\Sigma^{st}(Lie_A^f) \stackrel{\sim}{\to}Lie_A$ in an essentially unique way (cf. Section~\ref{stablederived}).
Let
\[
\xymatrix{
Ch^\bullet:Lie_A  \ar@<0.5ex>[r] &   (\CAlg_{A//A})^{op}:\DD_\infty   \ar@<0.5ex>[l]
}
\]
be the adjoint pair where the left adjoint $Ch^\bullet$
is the ``Chevalley-Eilenberg cochain functor''
which carries $L\in Lie_A$ to the Chevalley-Eilenberg cochain complex $Ch^\bullet(L)$
(i.e., the $A$-linear dual of Chevalley-Eilenberg chain complex), see e.g. \cite[1.4]{H}, \cite[X, 2.2]{DAG}.
Thanks to \cite[1.5.6]{H}, this adjoint pair induces an adjoint pair
\[
\xymatrix{
\mathcal{F}:Lie_A  \ar@<0.5ex>[r] &   \FST_A :\mathcal{L}.\ar@<0.5ex>[l]  
}
\]
Moreover, if $A$ is noetherian (cf. ection~\ref{sectionbackground}) both $\mathcal{F}$ and $\mathcal{L}$ are categorical equivalences. 
The left adjoint $\mathcal{F}$ is defined as follows. The restriction of the functor $Ch^\bullet$ to $Lie_A^f$
induces $Lie_A^f\to \TSZ\subset (\CAlg_{A//A})^{op}$ such that $Lie_A^f\stackrel{Ch^\bullet}{\to} \TSZ_A\hookrightarrow \FST_A$ belongs to $\Fun^{st}(Lie_A^f,\FST_A)$.
There exists an essentially unique left adjoint functor $\mathcal{F}:Lie_A\simeq \PPP_\Sigma^{st}(Lie_A^f)\to \FST_A$
which extends $Lie_A^f\to \FST_A$ (cf. the final equivalence in Section~\ref{stablederived}).
The right adjoint $\mathcal{L}$ is induced by the composition with $Lie_A^f\to \TSZ_A$ where we think of $Lie_A$ as the full subcategory
of $\Fun((Lie_A^f)^{op},\SSS)$.
If $A$ is noetherian, $\mathcal{F}$ and $\mathcal{L}$ are
reduced to a pair of mutually inverse functors $Ch^\bullet:Lie_A^f\simeq \TSZ:\DD_\infty$.
For $L\in Lie_A$, we usually write $\mathcal{F}_L$ for the image $\mathcal{F}(L)$.

\subsection{Koszul duals of $\eenu$-algebras}
\label{koszuldual}
Let $n\ge1$ be a natural number and  let $\Alg_{n}(\Mod_A)$ be the $\infty$-category of $\eenu$-algebras in $\Mod_A$.
We write $\Alg^+_{n}(\Mod_A)$ for the $\infty$-category
$\Alg_{n}(\Mod_A)_{/A}$ of augmented $\eenu$-algebras.
We review the Koszul duals of augmented $\eenu$-algebras.
Let $p_{B}:B\to A$ and $p_C:C\to A$ be augmented $\eenu$-algebras.
Let $\Map_{\Alg_{n}(\Mod_A)}(B\otimes_AC,A)\to \Map_{\Alg_{n}(\Mod_A)}(B,A)\times \Map_{\Alg_{n}(\Mod_A)}(C,A)$ be the map induced by the compositions with $B=B\otimes_AA\to B\otimes_AC$
and $C=A\otimes_AC\to B\otimes_AC\to A$.
We denote by $\textup{Pair}(p_B,p_C)$ or simply $\textup{Pair}(B,C)$ the fiber product
\[
\Map_{\Alg_{n}(\Mod_A)}(B\otimes_AC,A)\times_{\Map_{\Alg_{n}(\Mod_A)}(B,A)\times \Map_{\Alg_{n}(\Mod_A)}(C,A)}\{(p_B,p_C)\}.
\]
We shall refer to $\textup{Pair}(B,C)$ as the space of pairing of $p_B:B\to A$ and $p_C:C\to A$.
The functor $\Alg^+_{n}(\Mod_A)^{op}\to \SSS$ given by $C\mapsto \textup{Pair}(B,C)$
is representable by an object $\DD_{\eenu}(B)\in \Alg_{\eenu}^+(\Mod_A)$.
We shall call $\DD_{\eenu}(B)$ the $\eenu$-Kosuzl dual of $B\to A$.
For ease of notation, we write $\DD_{n}$ for $\DD_{\eenu}$.
There is a universal/tautological pairing $B\otimes_A\DD_{n}(B)\to A$ which corresponds to
the identity map $\textup{id}\in \Map_{\Alg_{n}(\Mod_A)}(\DD_{n}(B),\DD_{n}(B))$.
The Koszul dual $\DD_{n}(B)$ can also be interpreted as a centralizer of $B\to A$, see \cite[5.3.1]{HA}.
Thanks to the construction in
\cite[X, 4.4.6]{DAG} or \cite[5.2.5.5]{HA} or Construction~\ref{smartconstruction}, the assignment $B\mapsto \DD_{n}(B)$ is promoted to a ($\eenu$-Koszul duality) functor
$\DD_{n}:\Alg^+_{n}(\Mod_A)^{op} \to \Alg^+_{n}(\Mod_A)$ whose right adjoint
is $\DD_{n}:\Alg^+_{n}(\Mod_A)\to \Alg^+_{n}(\Mod_A)^{op}$.
The Koszul duals can be described in terms of bar constructions.
Let $\textup{Bar}:\Alg_{1}^+(\Mod_A)\to \Mod_{A}$ be the functor
given by $[B\to A] \mapsto A\otimes_BA$, which we refer to as the bar construction of augmented algebras.
Then $A\otimes_BA$ admits a structure of a coaugmented coalgebra in a suitable way, and $\textup{Bar}$
is promoted to $\textup{Bar}:\Alg_{1}^+(\Mod_A)\to \Alg^+_{1}((\Mod_A)^{op})^{op}$
(here we abuse notation by using the same symbol).
Interating bar construction we have
induces $\textup{Bar}^n:\Alg_{n}^+(\Mod_A)\to \Alg_{n}^+((\Mod_A)^{op})^{op}$.
According to \cite[X, 4.4.20]{DAG}, for $B\in \Alg^+_{n}(\Mod_A)$, the $\eenu$-Koszul dual $\DD_{n}(B)$
is equivalent to
the $A$-linear dual $\mathcal{H}om_{A}(\textup{Bar}^n(B),A)$.

We remark that there exists a Koszul duality functor in more general setting.
Let $\mathcal{M}$ be a symmetric monoidal presentable $\infty$-category whose tensor product
$\mathcal{M}\times\mathcal{M}\to \mathcal{M}$ preserves small colimits separately in each variable.
If we replace $\Mod_A$ by $\mathcal{M}$, 
there exists a $\eenu$-Koszul duality functor $(\Alg_{n}(\mathcal{M})_{/\uni})^{op} \to (\Alg_{n}(\mathcal{M})_{/\uni})$
which carries  an augmented algebra $B\to \uni$ 
to $\DD_{n}(B)\in \Alg_{n}(\mathcal{M})_{/\uni}$ 
where $\uni$ is a unit algebra (see \cite[5.2.5.8]{HA}).

\subsection{Universal enveloping algebras}
\label{UEA}

We will define an adjoint pair
\[
\xymatrix{
U_n:Lie_A  \ar@<0.5ex>[r] & \Alg_{n}^+(\Mod_A)  : res_{\eenu/Lie}   \ar@<0.5ex>[l]
}
\]
where the left adjoint $U_n$ sends $L$ to a universal enveloping $\eenu$-algebra of $L$.
We assume that $A$ is noetherian.
We first define the right adjoint.
Let $Y_{\eenu}^\vee:\Alg_{n}^+(\Mod_A)\to \Fun(\Alg_{n}^+(\Mod_A),\SSS)$ be the functor defined as the composite
\[
\Alg_{n}^+(\Mod_A)\to \Fun( \Alg_{n}^+(\Mod_A)^{op},\SSS)\to \Fun( \Alg_{n}^+(\Mod_A),\SSS)
\]
where the left functor is the Yoneda embedding and the right functor is induced by the composition with $\DD_{n}$.
Namely, $Y_{\eenu}^\vee(R)$ is given by $B\mapsto \Map_{\Alg_{n}^+(\Mod_A)}(\DD_{n}(B),R)$.

Let $\CAlg^+_{A}=\CAlg(\Mod_A)_{/A}\to \Alg_{n}^+(\Mod_A)$ be the forgetful functor.
Composition with this functor induces $\operatorname{Res}_{n/\infty}:\Fun( \Alg_{n}^+(\Mod_A),\SSS)\to \Fun(\EXT_A,\SSS)$.
Thus we obtain
\[
res_{n/\infty}:\Alg_{n}^+(\Mod_A)\stackrel{Y_{\eenu}^\vee}{\longrightarrow} \Fun(\Alg_{n}^+(\Mod_A),\SSS)\to \Fun(\EXT_A,\SSS).
\]
The essential image of this composite is contained in $\FST_A\subset \Fun(\EXT_A,\SSS)$.
To see this, it will suffice to verify that
\begin{enumerate}
\item  For a finite product $B_1\times \cdots \times B_m$ in $\EXT_A$,
a finite coproduct $\DD_{n}(B_1)\sqcup \cdots \sqcup \DD_{n}(B_m)$ in $\Alg_{n}^+(\Mod_A)$
is its $\eenu$-Koszul dual $\DD_{n}(B_1\times \cdots \times B_m)$.

\item For $B$ in $\EXT_A$, $\DD_{n}(A\times_{B}A)\simeq \DD_{n}(A)\sqcup_{\DD_{n}(B)}\DD_{n}(A)\simeq A\sqcup_{\DD_{n}(B)}A$.

\end{enumerate}

These facts follows from \cite[X, 4.4.5, 4.5.6]{DAG} when $A$ is a field of characteristic zero.
The general case for $A\in \CCAlg_k$ follows from Remark~\ref{Koszuldeduction} below. 
Consequently, the composition $\Alg_{n}^+(\Mod_A)\to \FST_A$ and the equivalence $\FST_A\simeq Lie_A$ ($A$ is noetherian)
gives us
\[
res_{\eenu/Lie}:\Alg_{n}^+(\Mod_A)\longrightarrow Lie_A.
\]

\begin{Lemma}
The functor $res_{\eenu/Lie}$ admits a left adjoint functor $U_n:Lie_A\to :\Alg_{n}^+(\Mod_A)$.
\end{Lemma}

\Proof
Note first that both $\Alg_{n}^+(\Mod_A)$ and $Lie_A$
are presentable $\infty$-categories, so that
by the adjoint functor theorem, it is enough to prove that $\Alg_{n}^+(\Mod_A)\to \FST_A\simeq Lie_A$ preserves small limits and filtered colimits.
Note that $Y_{\eenu}^\vee$ preserves small limits, and the inclusion $\FST_A\subset \Fun(\EXT_A,\SSS)$
preserves small limits (since it admits a left adjoint). Thus $\Alg_{n}^+(\Mod_A)\to \FST_A$ preserves small limits.
Taking into account the definition $\PPP_{\Sigma}^{st}(\TSZ_A)=\FST_A$,
we see that $\FST_A\hookrightarrow \Fun(\EXT_A,\SSS)$ is stable under filtered colimits.
Thus it will suffice to show that $res_{n/\infty}$ preserves filtered colimits.
It is enough to show that $\DD_{n}(B)$ is a compact object in $\Alg_{n}^+(\Mod_A)$ for any $B(\to A) \in \EXT_A$.
But when $B\to A$ is $A\oplus M\to A$, $\DD_{n}(B)$ is a free augmented $\eenu$-algebra object $\Free_{\eenu}(M^\vee[-n])$ (see Remark~\ref{Koszuldeduction}).
Thus our claim follows from the compactness of $M$ in $\Mod_A$.
\QED
We obtain an adjoint pair
\[
\xymatrix{
U_n:Lie_A  \ar@<0.5ex>[r] &    \Alg_{n}^+(\Mod_A): res_{\eenu/Lie}.   \ar@<0.5ex>[l]  
}
\]
By composing $\textup{forget}:\Alg_{n}^+(\Mod_A)\rightleftarrows \Alg_{n}(\Mod_A)$,
we also have
\[
\xymatrix{
U_n:Lie_A  \ar@<0.5ex>[r] &    \Alg_{n}(\Mod_A): res_{\eenu/Lie}.   \ar@<0.5ex>[l]  
}
\]
where by abuse of notation we use the same symbols $U_n$ and $res_{\eenu/Lie}$. 
We shall refer to $U_n$ as the universal enveloping $\eenu$-algebra functor.

\begin{Remark}
\label{Koszuldeduction}
Let $B(\to A)$ be an augmented $\eenu$-algebra over $A$ and let $B\to \DD_{n}\DD_{n}(B)$ be the biduality morphism (induced by the universality).
Let $M$ be an $A$-module spectrum.
Suppose that $M$ is of the form $\oplus_{1\le i \le m}A^{\oplus r_i}[d_i]$ ($d_i\le -n$).
Let $\Free_{\eenu}:\Mod_A\to \Alg^+_n(\Mod_A)$ denote the free functor,
 that is a left adjoint of the forgetful functor.
We set $B=\Free_{\eenu}(M)$. We show that the biduality map $B\to \DD_{n}\DD_{n}(B)$
is an equivalence.
Then (the proof of) \cite[DAG X 4.5.6]{DAG} shows that
there is a canonical equivalence $\DD_{n}(\Free_{\eenu}(M))\simeq A\oplus M^\vee[-n]$
where $A\oplus M^\vee[-n]$ indicates the trivial square zero extension.
Consequently, we have the canonical morphism
$\Free_{\eenu}(M) \to \DD_{n}\DD_{n}(\Free_{\eenu}(M))\simeq \DD_{n}(A\oplus M^\vee[-n])$. 
It will suffice to prove that $\Free_{\eenu}(M) \to \DD_{n}(A\oplus M^\vee[-n])$ is an equivalence.
When $A$ is $k$, then the desired equivalence follows from \cite[4.4.5, 4.5.6]{DAG}.
We will describe how to deduce the general case from the case when $A=k$.
Let $\Free_{k}:\Mod_k\to \Alg_{n}^+(\Mod_k)$ be the free functor (i.e., a left adjoint of the forgetful functor), and we write $\Free_A$ for the above $\Free_{\eenu}$. Let $\DD_{n,k}:\Alg_{n}^+(\Mod_k)^{op} \to \Alg_{n}^+(\Mod_k)$
be the Koszul duality functor.
Note canonical equivalences
\[
\Free_A(\oplus_{1\le i \le m}A^{\oplus r_i}[d_i])\simeq \Free_k(\oplus_{1\le i \le m}k^{\oplus r_i}[d_i])\otimes_kA\simeq \DD_{n,k}(k\oplus (\oplus_{1\le i \le m}k^{\oplus r_i}[d_i])^\vee[-n])\otimes_kA.
\]
To prove that $\Free_A(\oplus_{1\le i \le m}A^{\oplus r_i}[d_i])\to \DD_{n}(A\oplus (\oplus_{1\le i \le m}A^{\oplus r_i}[d_i])^\vee[-n])$ is an equivalence,
it will suffice to show that for a connective perfect $k$-module $N$, $\DD_{n,k}(k \oplus N)\otimes_kA\to \DD_{n}(A\otimes_k(k \oplus N))$ is 
an equivalence in $\Mod_A$. It is obtained from
\begin{eqnarray*}
\DD_{n,k}(k \oplus N)\otimes_kA &\simeq& \mathcal{H}om_k(\BAR_k^n(k\oplus N),k)\otimes_kA \\
&\simeq& \mathcal{H}om_k(\BAR_k^n(k\oplus N),A) \\
&\simeq& \mathcal{H}om_A(\BAR_A^n(A\otimes_k(k \oplus N)),A) \\
&\simeq & \DD_{n}(A\otimes_k(k \oplus N))
\end{eqnarray*}
where $\BAR_k$ and $\BAR_A$ indicate the bar construction of augmented algebras in $\Mod_k$ and $\Mod_A$,
respectively.
Here we denote by $\mathcal{H}om_k$ and $\mathcal{H}om_A$ the internal Hom objects in $\Mod_k$ and $\Mod_A$, respectively.
The second equivalence (key point) follows from the fact that the homology of $\BAR_k^n(k\oplus N)$
is a finite dimensional $k$-vector space in each degree.
Namely, each degree of $\mathcal{H}om_k(\BAR_k^n(k\oplus N),k)\simeq \Free_k(N^\vee[-n])$ is finite dimensional: it can be deduced from the description
of free $\eenu$-algebras in terms of configuration spaces of Euclidean spaces
and the estimate of the amplitude and the finite dimensionality of homology groups of configuration spaces (see e.g. \cite[X 4.1.15]{DAG}).
The first and final equivalences are the presentation of Koszul duals in terms of iterated bar constructions
(see Section~\ref{koszuldual}).
\end{Remark}

\begin{Proposition}
\label{UEAProp}
The followings hold.

\begin{enumerate}
\item The left adjoint functor $U_n:Lie_A\to \Alg^+_{n}(\Mod_A)$ is an essentially unique colimit-preserving functor
 which extends the composite
\[
Lie_A^{f}\stackrel{Ch^\bullet}{\longrightarrow} (\CAlg_A^+)^{op}\stackrel{\textup{forget}}{\longrightarrow}(\Alg_{n}^+(\Mod_A) )^{op} \stackrel{\DD_{n}}{\longrightarrow} \Alg_{n}^+(\Mod_A).
\]

\item The composite $res_{\eenu/Lie}:\Alg^+_{n}(\Mod_A)\to Lie_A\stackrel{\textup{forget}}{\to} \Mod_A$ is equivalent to the functor defined by the formula $[B\to A]\mapsto \Ker(B\to A)[n-1]$.

\end{enumerate}
\end{Proposition}

\Proof
We prove (1).
By the equivalences $\Fun^{\textup{L}}(\PPP^{st}_\Sigma(Lie_A^f),\Alg_{n}^+(\Mod_A))\simeq \Fun^{st}(Lie_A^f,\Alg_{n}^+(\Mod_A))$ and $\PPP^{st}_\Sigma(Lie_A^f)\simeq Lie_A$,
it will suffice to show that $Lie_A^f\hookrightarrow Lie_A\stackrel{U_n}{\to}  \Alg_{n}^+(\Mod_A)$
is equivalent to $\DD_{n}\circ Ch^\bullet|_{Lie_A^f}$.
Unwinding the definition of $res_{\eenu/Lie}$,
we see that for $C\in \Alg^+_n(\Mod_A)$, $res_{\eenu/Lie}(C)\in \PPP_\Sigma^{st}(Lie_A^f)\subset \Fun(Lie_A^f,\SSS)$ (given by $L\mapsto \Map_{Lie_A}(L,res_{\eenu/Lie}(C))$) is
equivalent to the functor given by $L\mapsto \Map_{\Alg^+_{n}(\Mod_A)}(\DD_{n}(Ch^\bullet(L)),C)$. 
It follows that the restriction of $U_n$ to $Lie_A^f$ is equivalent to $\DD_{n}\circ Ch^\bullet|_{Lie_A^f}$.

Next, we prove (2).
To observe (2),
we first consider 
\[
(\Mod_{A})^{op}\stackrel{\Free_{Lie}}{\longrightarrow} (Lie_A)^{op} \stackrel{Ch^\bullet}{\longrightarrow} \CAlg_A^+.
\]
Taking into account the universal property of $\mathcal{P}^{st}_{\Sigma}$ (see Section~\ref{stablederived}), 
$(\Mod_A^{f})^{op} \to (Lie_A^{f})^{op} \to \EXT_A$ induces colimit-preserving functors of presentable $\infty$-categories
\[
\mathcal{P}^{st}_{\Sigma}(\Mod_A^f)\to \mathcal{P}^{st}_{\Sigma}(Lie_A^f) \to \mathcal{P}^{st}_{\Sigma}(\TSZ_A)=\FST_A
\]
where the left functor can be identified with $\Free_{Lie}:\Mod_A\to Lie_A$,
and the final equality is the definition of $\FST_A$.
Note that the right adjoint functors of this sequence are the restriction of the sequence 
\[
\Fun(\EXT_A,\SSS)\to \Fun((Lie_A^f)^{op},\SSS)\to \Fun((\Mod_A^f)^{op},\SSS)
\]
given by the composition with $(\Mod_A^{f})^{op}\to (Lie_A^{f})^{op} \to \EXT_A$.
In order to verify that $\Alg^+_{n}(\Mod_A)\to Lie_A\to \Mod_A$ sends $B\to A$ to $\Ker(B\to A)[n-1]$,
by (1) and left Kan extension to $\mathcal{P}^{st}_{\Sigma}(\Mod_A^f)\simeq \Mod_A$
it is enough to prove that the composite
\[
(\Mod_{A}^f) \stackrel{\Free_{Lie}}{\to} (Lie^f_A) \stackrel{Ch^\bullet}{\to} (\CAlg_A^+)^{op} \stackrel{\textup{forget}}{\to} \Alg_{n}^+(\Mod_A)^{op} \stackrel{\DD_{n}}{\to} \Alg_{n}^+(\Mod_A)
\]
is equivalent to the ``shifted free functor'' given by $M\mapsto \Free_{\eenu}(M[1-n])$.
The composite 
$(\Mod_{A}^{f})^{op} \to (Lie_A)^{op} \to \CAlg_A^+$
is equivalent to $M \mapsto[ A\oplus M^\vee[-1]\to A]$ (see \cite[1.4.11]{H}).
Thus, we are reduced to proving that
$\EXT_A \to \Alg_{n}^+(\Mod_A) \stackrel{\DD_{n}}{\to}(\Alg_{n}^+(\Mod_A))^{op}$
is equivalent to the functor given by $[A\oplus N\to A]\mapsto \Free_{\eenu}(N^{\vee}[-n])$.
This is a consequence of the canoncial equivalence $\DD_{n}(A\oplus N)\simeq \Free_{\eenu}(N^\vee[-n])$ where $N$ is a connective $A$-module of the form
$\oplus_{1\le i \le n}A^{\oplus r_i}[d_i]$ (see Remark~\ref{Koszuldeduction}). 
\QED

\begin{Remark}
The functor $U_1:Lie_A\to \Alg^+_1(\Mod_A)$ defined above is equivalent 
to
the standard definition of the universal enveloping algebras 
$\mathcal{U}:Lie_A\to \Alg^+_1(\Mod_A)$ (see e.g. \cite[X, 2.1.7]{DAG}, \cite{H}
for the universal enveloping algebras of dg Lie algebras).
We first consider the case when $A=k$.
The result from \cite[3.3.2]{DAG}
says that $Ch^\bullet:(Lie_k)^{op}\stackrel{Ch^\bullet}{\to} \CAlg_k^+\stackrel{\textup{forget}}{\to} \Alg_1^+(\Mod_k)$
is equivalent to $\DD_1\circ \mathcal{U}:(Lie_k)^{op}\to  \Alg_1^+(\Mod_k)$.
If $L\in Lie_k^f$, it follows from \cite[3.1.5]{DAG} that the the canonical
map $\mathcal{U}(L)\to \DD_1\circ\DD_1(\mathcal{U}(L))$ is an equivalence.
Thus,
there exists 
natural equivalences
$\DD_1\circ Ch^\bullet\simeq \DD_1\circ \DD_1\circ \mathcal{U}\simeq \mathcal{U}$
between functors $\Fun(Lie_k^f,\Alg_1^+(\Mod_k))$
(here we omit the forgetful functor $\CAlg_k^+\to \Alg_1^+(\Mod_k)$ from the notation).
Since $\mathcal{U}$ is a left adjoint, 
we deduce from Proposition~\ref{UEAProp} (1) 
that $\mathcal{U}\simeq U_1$ in $\Fun(\PPP^{st}_{\Sigma}(Lie_k^f),\Alg^+_1(\Mod_k))\simeq \Fun(Lie_k,\Alg^+_1(\Mod_k))$.
Next we consider the general case.
The proof of \cite[3.3.2]{DAG} reveals that it holds in the general case.
According to Remark~\ref{Koszuldeduction},
$\mathcal{U}(L)\to \DD_1\circ \DD_1(\mathcal{U}(L))$ is an equivalence for $L\in Lie_A^f$.
The general case follows from these observations and the argument in the case $A=k$.
\end{Remark}

\subsection{Hochschild homology and Hochschild cohomology}
\label{hochschildhomologysection}

Let $\ST_A$ denote $\Mod_{\Perf_A^\otimes}(\ST)$.
Let
\[
\HH_\bullet(-/A):\ST_A \longrightarrow \Mod_A^{S^1}=\Fun(BS^1,\Mod_A)
\]
be the symmetric monoidal functor which carries $\mathcal{C} \in \ST_A$ to the Hochschild homology $A$-module spectrum $\HH_\bullet(\mathcal{C}/A)$
over $A$.
We refer the reader to \cite[Section 6, 6.14]{I} for the construction
of the Hochschild chain/homology functor $\HH_\bullet(-/A)$
(in {\it loc. cit.} we use the symbol $\HH_\bullet(\mathcal{C})$ instead of  $\HH_\bullet(\mathcal{C}/A)$).
Let $\Perf_{(-)}:\Alg_{1}(\Mod_A)\to \ST_A$ denote the symmetric monoidal functor
given by $R \mapsto \Perf_R$
(see Section~\ref{sectionmodule}).
We will write  $\HH_\bullet(-/A)$
also for
the composite functor 
\[
\Alg_{1}(\Mod_A)\stackrel{\Perf_{(-)}}{\longrightarrow} \ST_A\stackrel{\HH_\bullet(-/A)}{\longrightarrow}\Mod_A^{S^1}.
\]
According to \cite[Remark 7.10, Proposition 7.11]{I} and its proof,
$\HH_\bullet(-/A):\Alg_{1}(\Mod_A)\to \Mod_A^{S^1}$ is canonically equivalent to
the functor $\int_{S^1}:\Alg_{1}(\Mod_A) \to  \Mod_A^{S^1}$ which assings to $R$ the factorization homology $\int_{S^1}R$.
Moreover, we will use:

\begin{Lemma}
\label{commutativemodulitransform2}
For ease of notation, we let $h:\CAlg(\Mod_A)\simeq \CAlg(\Alg_{1}(\Mod_A))\to \CAlg(\Mod^{S^1}_A)$
denote the functor induced by $\HH_\bullet(-/A):\Alg_{1}(\Mod_A)\to \Mod^{S^1}_A$.
(The equivalence $\CAlg(\Mod_A)\simeq \CAlg(\Alg_{1}(\Mod_A))$ follows from
Dunn additivity theorem.)
Let $i:\CAlg(\Mod_A)\to \CAlg(\Mod^{S^1}_A)$
denote the functor which carries $R$ to $R$ with the trivial $S^1$-action.
Then $h:\CAlg(\Mod_A)\to \CAlg(\Mod^{S^1}_A)$
is equivalent to the functor $\otimes_AS^1$ gieven by the tensor with 
$S^1\in \SSS$, and there is a natural transformation $\sigma:h\to i$
induced by the $S^1$-equivariant map $S^1\to \ast$ into the contractible space.
\end{Lemma}

\Proof
We will describe the construction of $\sigma:h\to i$.
To this end, we describe
$h:\CAlg(\Mod_A)\simeq \CAlg(\Alg_{1}(\Mod_A))\to \CAlg(\Mod^{S^1}_A)$
in terms of cyclic sets.
For this purpose,
we briefly review the Hochschild cyclic objects:
the construction of $\HH_\bullet(-/A):\Alg_{1}(\Mod_A)\to \Mod^{S^1}_A$
in \cite{I}.
We use the theory of symmetric spectra. The readers who do not know symmetric spectra
are invited to skip the construction on the first reading.
We use the notation and terminology in \cite[Section 6]{I}. We refer the reader to {\it loc. cit.} for details.
The symmetric monoidal functor $\HH_\bullet(-/A)$ is obtained from the composite of symmetric monoidal functors
\[
\Alg_{1}(\SPS(\mathbb{A})^c)\stackrel{\mathcal{HH}(-)_\bullet}{\longrightarrow} \Fun(\Lambda^{op},\SPS(\mathbb{A})^c)\to \Fun(\Lambda^{op},\SPS(\mathbb{A})^c[W^{-1}])\stackrel{L}{\to} \Fun(BS^1,\SPS(\mathbb{A})[W^{-1}])
\]
by inverting weak equivalences in $\Alg_{1}(\SPS(\mathbb{A})^c)$
(see the construction before \cite[Lemma 6.11]{I}).
 Here $\mathbb{A}$ is a cofibrant commutative symmetric spectrum
 which is a model of $A$, $\SPS(\mathbb{A})^c$ is the category of cofibrant $\mathbb{A}$-module symmetric
spectra, and $\Lambda$ is the cyclic category.
We denote by $\SPS(\mathbb{A})^c[W^{-1}]$ the (symmetric monoidal) $\infty$-category
obtained from $\SPS(\mathbb{A})^c$ by inverting weak equivalences (see \cite[1.3.4, 4.1.7, 4.1.8]{HA}
for localizations with respect to weak equivalences).
There is a canonical functor $\SPS(\mathbb{A})^c\to \SPS(\mathbb{A})^c[W^{-1}]$
that induces the middle functor. The third (right) functor is the symmetric monoidal functor
determined by left Kan extensions along the groupoid completion $\Lambda^{op}\to BS^1$. 
There are relationships $\SPS(\mathbb{A})[W^{-1}]\simeq \Mod_A$ and
$\Alg_{1}(\SPS(\mathbb{A})^c)[W^{-1}]\simeq \Alg_{1}(\Mod_A)$.
The first functor $\mathcal{HH}(-)_\bullet$ carries $\mathbb{R}\in \Alg_{1}(\SPS(\mathbb{A})^c)$
to the Hochschild cyclic object $\mathcal{HH}(\mathbb{R})_\bullet:\Lambda^{op}\to \SPS(\mathbb{A})^c$ (see e.g. \cite[Definition 6.18]{I} for the formula: in {\it loc. cit.}, spectral categories are used in the definition instead of $\Alg_{1}(\SPS(\mathbb{A})^c)$).
Consider the composite functor  
$\CAlg(\SPS(\mathbb{A})^c)\simeq \CAlg(\Alg_{1}(\SPS(\mathbb{A})^c))\to \Fun(\Lambda^{op},\CAlg(\SPS(\mathbb{A})^c))$ induced by $\mathcal{HH}(-)_\bullet$.
Let $\Delta^1/\partial\Delta^1$ be the standard simplicial model of the circle $S^1$.
The simpliclal set $\Delta^1/\partial\Delta^1:\Delta^{op}\to \textup{Sets}$
is extended to a cyclic set $C:\Lambda^{op}\to \textup{Sets}$ (see \cite[6.1.9]{L}).
Then
$\CAlg(\SPS(\mathbb{A})^c))\to \Fun(\Lambda^{op},\CAlg(\SPS(\mathbb{A})^c))$
is given by the formula
\[
\RR\mapsto \mathcal{HH}(\mathbb{R})_\bullet:=[\Lambda^{op}\ni [p] \mapsto C([p])\otimes \RR \in \CAlg(\SPS(\mathbb{A})^c)]
\]
where 
$C([p])\otimes \RR$ is the tensor product with the finite set $C([p])$ in $\CAlg(\SPS(\mathbb{A})^c)$,
that is, the $(p+1)$-fold (homotopy) coproduct 
($(p+1)$-fold smash product) $\wedge^{p+1}\mathbb{R}$.

To define $h\to i$, we use the ``contraction'' of $C$ to the terminal cyclic set.
Let $C_{\ast}$ denote the cyclic set which is the constant functor whose value is
the set consisting of one element $\ast$ (that is, the terminal/final cyclci set).
The (unique) map $C\to C_\ast$ determines
the natural transformation $\sigma_\Lambda:\mathcal{HH}(-)_\bullet\to \textup{const}$
between functors $\CAlg(\SPS(\mathbb{A})^c)\to \Fun(\Lambda^{op},\CAlg(\SPS(\mathbb{A})^c))$
where $\textup{const}(\mathbb{R})$ is the constant cyclic object with value $\mathbb{R}$.
For any $\RR$ and $p\ge0$,
$\mathcal{HH}(\mathbb{R})_p=\wedge^{(p+1)}\mathbb{R}\to \mathbb{R}$
induced by $\sigma_{\Lambda}$
is the multiplication map $\wedge^{(p+1)}\mathbb{R}\to \mathbb{R}$

By inverting weak equivalences, we have the diagram
\[
\xymatrix{
\CAlg(\SPS(\mathbb{A})^c)[W^{-1}]\ar[r] \ar[d]^{\simeq} &  \Fun(\Lambda^{op},\CAlg(\SPS(\mathbb{A})))[W^{-1}]\ar[d] \\
\CAlg(\Alg_{1}(\SPS(\mathbb{A})^c)[W^{-1}])\ar[r] & \CAlg(\Fun(\Lambda^{op},\SPS(\mathbb{A})^c[W^{-1}]))
}
\]
which commutes up to canonical homotopy, where the horizontal arrows are defined
by taking Hochschild cyclic objects $\mathcal{HH}(-)_\bullet$, and $[W^{-1}]$
indicates the localization of weak equivalences.
Let 
\[
h_\Lambda:\CAlg_A\simeq \CAlg(\Alg_{1}(\SPS(\mathbb{A})^c)[W^{-1}])\to \CAlg(\Fun(\Lambda^{op},\SPS(\mathbb{A})^c[W^{-1}]))\simeq \Fun(\Lambda^{op},\CAlg_A)
\]
denote the functor obtained by the composition with the lower horizontal functor.
(By the universal property of $[W^{-1}]$, this functor is given by
$R\mapsto [\Lambda^{op}\ni [p] \mapsto C([p])\otimes R \in \CAlg_A]$.)
Let $\textup{const}_\infty$ denote the functor $\CAlg_A\to  \Fun(\Lambda^{op},\CAlg_A)$
which carries $R$ to the constant functor with value $R$.
 Then $\sigma_{\Lambda}$ determines $\sigma'_{\Lambda}:h_{\Lambda}\to\textup{const}_\infty$.
 The composition of $\sigma_\Lambda$ and $\Fun(\Lambda^{op},\CAlg_A)\to \Fun(BS^1,\CAlg_A)$ induced by
$L:\Fun(\Lambda^{op},\Mod_A)\to \Fun(BS^1,\Mod_A)$
determines $\sigma:h\to i$ between functors $\CAlg_A\to \Fun(BS^1,\CAlg_A)\simeq \CAlg(\Mod^{S^1}_A)$.

Let $L_{\SSS}:\Fun(\Lambda^{op},\SSS) \to \Fun(BS^1,\SSS)$
be a left adjoint functor of the functor $\Fun(BS^1,\SSS)\to \Fun(\Lambda^{op},\SSS)$
induced by the composition with $\Lambda^{op}\to BS^1$.
Then $L_{\SSS}$
sends $C$ to $S^1$ with an $S^1$-action.
Indeed, we note that the cyclic set $C:\Lambda^{op}\to \textup{Sets}\subset \SSS$
is represented by by $[0]\in \Lambda$.
The left adjoint functor $L_{\SSS}$ is an essentially unique colimit-preserving functor
which extends $\Lambda \to (BS^1)^{op}\hookrightarrow \Fun(BS^1,\SSS)$
where the second functor is Yoneda embedding.
If we denote the unique object of $BS^1$ by $\star$,
it follows that $L_{\SSS}(C)$ is corepresented by $\star$.
That is, $L_{\SSS}(C)$ amounts to
$\Map_{BS^1}(\star,\star)=S^1$ with the $S^1$-action
determined by the multiplication $S^1\times S^1\to S^1$.

For a cyclic object $E:\Lambda^{op}\to \mathcal{M}$
in an $\infty$-category $\mathcal{M}$, $F:BS^1\to \MMM$ is a left Kan extension
of $E$ along $\Lambda^{op}\to BS^1$ if and only if
the composite $\ast=\Delta^0\to BS^1\to \MMM$
is a colimit of the restricition
$E|_{\Delta^{op}}:\Delta^{op}\to \MMM$ (see e.g.  \cite[Lemma 6.9 (ii)]{I}).
(Since the restriction $C|_{\Delta^{op}}$ is $\Delta^1/\partial\Delta^1$,
it follows that the underlying space of $L_{\SSS}(C)$ is $S^1=|\Delta^1/\partial\Delta^1|$.) 
Note that $\varinjlim_{[p]\in \Delta^{op}}(C([p])\otimes R)\simeq \varinjlim_{[p]\in \Delta^{op}}(C([p]))\otimes R\simeq S^1\otimes R$ in $\CAlg_A$.
We then see that there is a canonical equivalence $h(R)\stackrel{\sim}{\to} L_{\SSS}(C)\otimes R$
in $\CAlg(\Mod_A^{S^1})$
because the underlying map
$\varinjlim_{[p]\in \Delta^{op}}(C([p])\otimes R)\to S^1\otimes R$
is an equivalence.
In particular, the $S^1$-action on $h(R)$ is induced by that of $L_{\SSS}(C)\simeq S^1$.
We may regard $\sigma:h\to i$ as the natural transformation
informally given by $h(R)\simeq R\otimes S^1\to R\otimes \ast=i(R)$ induced by the $S^1$-equivariant map
$S^1\to \ast$.
\QED

We review the Hochschild cohomology of $\CCC\in \ST_A$.
Let $\Ind(\CCC)$ be the Ind-category that belongs to
$\PR_A$. Moreover, it is compactly generated.
We denote by $\theta_A:\Alg_2(\Mod_A)\to \Alg_1(\PR_A)$
the functor informally given by $B\mapsto \LMod^\otimes_B$ (Section~\ref{sectionmodule}).
By definition, the endomorphism algebra object $\mathcal{E}nd_A(\Ind(\CCC))\in \Alg_1(\PR_A)$
endowed with a tautological action on $\Ind(\CCC)$ is a final object of
$\LMod(\PR_A)\times_{\PR_A}\{\Ind(\CCC)\}$.
There exists a final object of $\Alg_2(\Mod_A)\times_{\Alg_1(\PR_A)}\LMod(\PR_A)\times_{\PR_A}\{\Ind(\CCC)\}$.
To see this, note first that
there exists a right adjoint $\Alg_{1}(\PR_A)\to \Alg_{2}(\Mod_A)$ of
$\theta_A$, see \cite[4.8.5.11, 4.8.5.16]{HA}. Namely, there exists a adjoint pair
\[
\xymatrix{
\theta_A:\Alg_{2}(\Mod_A)  \ar@<0.5ex>[r] &    \Alg_{1}(\PR_A):E_A  \ar@<0.5ex>[l]  
}
\]
such that $E_A$ sends $\MM^\otimes$ to the endomorphism algebra of the unit object $\uni_{\MM}$.
Then its final object is given by $E_A(\mathcal{E}nd_A(\Ind(\CCC)))\in \Alg_2(\Mod_A)$
with the left module action of $\LMod_{E_A(\mathcal{E}nd_A(\Ind(\CCC)))}$ on $\Ind(\CCC)$
determined by the counit map $\theta_A(E_A(\mathcal{E}nd_A(\Ind(\CCC))))\to \mathcal{E}nd_A(\Ind(\CCC))$.
The Hochschild cohomology $\HH^\bullet(\CCC/A)$ 
is defined to be $E_A(\mathcal{E}nd_A(\Ind(\CCC)))$.
We refer to  $\HH^\bullet(\CCC/A)$ as the Hochschild cohomology of $\CCC$ (over $A$).

Though we use the word ``homology'' and ``cohomology'',
$\HH_\bullet(\CCC/A)$ and $\HH^\bullet(\CCC/A)$
are not graded modules obtained by passing to (co)homology but spectra (or chain complexes) with algebraic structures.
If there is no confusion likely to arise, in Section 6--8, for ease of notation,
we often use the simple symbols $\HHH$ and $\HHHH$
instead of $\HH_{\bullet}(\CCC/A)$ and $\HH^\bullet(\CCC/A)$,
respectively (cf. Notation~\ref{simplenotation1}).

\section{Guide}

Before proceeding to Section 5--8, we will highlight several points in view of Theorem~\ref{intromain2} in an informal way.

(i) Deformation problems we will study have the following form or its variant.
Let $\mathcal{M}^\otimes$ be a monoidal $\infty$-category which
satisfies a suitable condition.
Let $M$ be an object of $\mathcal{M}$. We will consider the deformation problem
of $M$.
Let $B\to \uni$ be an augmented algebra object of $\mathcal{M}$, that is,
an object of $\Alg_1^+(\mathcal{M})=\Alg_1(\mathcal{M})_{/\uni}$
where $\uni$ is the unit algebra.
The $\infty$-category of deformations of $M$ to $B$ is
$\RMod_{B}(\mathcal{M})\times_{\textup{red}, \mathcal{M}}\{M\}$.
Here $\textup{red}:\RMod_{B}(\mathcal{M})\to \mathcal{M}$
is given by the reduction functor $N\mapsto N\otimes_B\uni$. 

The vertical arrows in Theorem~\ref{intromain2} are given
by descent type data (Koszul dual type data) associated to deformations. 
Rougly, data is defined as follows.
Suppose that $\uni$ is a $B$-$B'$-bimodule.
Let $N$ be a right $B$-module. Then
$N\otimes_B\uni$ is a right $B'$-module.
Since $\uni$ can be thought of as a $B$-$\DD_1(B)^{op}$-bimodule (see Section~\ref{koszuldual}),
$N\otimes_B\uni$ is a right $\DD_1(B)^{op}$-module, that is,
a left $\DD_1(B)$-module. This is the descent/Koszul type data associated to
the deformation $N$.
The assignment
\[
(B\to \uni,\ N,\ M\simeq N\otimes_B\uni) \mapsto (\DD_1(B)\to \uni,\ \DD_1(B) \curvearrowright M)
\]
defines the ``duality functor''
\begin{eqnarray*}
\DD^{\textup{mod}}_{M,\mathcal{M}^\otimes}:\mathcal{D}ef'_{M}(\mathcal{M})^{op}:=(\Alg_1^+(\mathcal{M})\times_{\Alg_1(\mathcal{M})}\RMod(\mathcal{M})\times_{\textup{red},\mathcal{M}}\{M\})^{op} \ \ \ \ \ \ \ \ \ \ \ \ \ \ \ \ \ \ \\
\ \ \ \ \ \ \ \ \ \ \ \ \ \to \Alg_1^+(\mathcal{M})\times_{\Alg_1(\mathcal{M})}\LMod(\mathcal{M})\times_{\textup{forget},\mathcal{M}}\{M\}=:\LMod^+(\mathcal{M})_M
\end{eqnarray*}
which extends the $\eone$-Koszul duality functor $\DD_{1,\mathcal{M}^\otimes}:\Alg_1^+(\mathcal{M})^{op}\to \Alg_1^+(\mathcal{M})$ (cf. Section~\ref{koszuldual}).
Let $\mathcal{N}^\otimes$ be another monoidal $\infty$-category
and let $F:\mathcal{M}^\otimes\to \mathcal{N}^\otimes$
be a monoidal functor such that it induces $\mathcal{D}ef'_M(F):\mathcal{D}ef'_{M}(\mathcal{M})\to \mathcal{D}ef'_{F(M)}(\mathcal{N})$.
The duality functors are functorial in the sense that
$F$ gives rise to natural transformations
\[
\gamma^{\textup{mod}}:\LMod^+(F)_M\circ \DD^{\textup{mod}}_{M, \mathcal{M}^\otimes} \to \DD^{\textup{mod}}_{F(M),\mathcal{N}^\otimes}\circ \mathcal{D}ef'_M(F), \ \ \ \gamma:\Alg_1^+(F)\circ \DD_{1,\mathcal{M}^\otimes} \to \DD_{1,\mathcal{N}^\otimes}\circ \Alg_1^+(F)
\]
where $\LMod^+(F)_M$ is the induced functor
$\LMod^+(\mathcal{M})_M\to \LMod^+(\mathcal{N})_{F(M)}$,
and $\Alg_1^+(F)$ is the induced functor
$\Alg_1^+(\mathcal{M})\to \Alg_1^+(\mathcal{N})$.
These natural transformations need not to be equivalences.
Consider the case
where $\mathcal{M}^\otimes=\ST_A^\otimes$, $\mathcal{N}^\otimes=\Mod_A^{S^1}$
and $F=\HH_\bullet(-/A):\ST_A\to \Mod_A^{S^1}$.
One of the key steps to the construction of the left commutative square in Theorem~\ref{intromain2} is to prove that
$\gamma^{\textup{mod}}$ is an equivalence after the restriction
along the functor $\EXT_A\to \Alg_1^+(\ST_A)$ given by $R\mapsto \Perf_R^\otimes$. It is reduced to the Koszul duality result
which says that $\gamma$ is an equivalence after the restriction
(cf. Propoistion~\ref{firstkoszuldual}).

(ii)
We need the interaction of Koszul duality with the change of operads
including
the Lie operad, the $\eone$-operad (the associative operad),
the $\etwo$-operad and the $\einf$-operad (commutative operad)
and with Hochschild chain functor.
Arguably, the key equivalences are given by the commutative diagram
of $S^1$-equivariant associative algebras
\[
\xymatrix{
U_1(\DD_{\infty}(R)^{S^1}) \ar[r]_{\simeq}  & \DD_{1}(R\otimes_AS^1)  \\
\HH_\bullet(U_2(\DD_{\infty}(R))/A) \ar[r]_{\simeq} \ar[u]_{\simeq} &  \HH_\bullet(\DD_{2}(R)/A) \ar[u]^{\simeq}  
}
\]
for $R\in \EXT_A$ (see Proposition~\ref{commutativediagram}).
These equivalences allow us to describe
the bottom row of the diagram in Theorem~\ref{intromain2}
in terms of dg Lie algebras (which coherently commutes with moduli-theoretic
construction).
We find that they are controlled by
the algebraic structure on $(\HH^\bullet(\CCC/A),\HH_\bullet(\CCC/A))$.
In particular, the Koszul dual of the cyclic deformation
of $\HH_\bullet(\CCC/A)$ associated to a deformation
of an $A$-linear stable $\infty$-category $\CCC$ is determined by the $S^1$-equivariant action of  $\GG_{\CCC}^{S^1}$
on $\HH_\bullet(\CCC/A)$ that arises from $(\HH^\bullet(\CCC/A),\HH_\bullet(\CCC/A))$ (cf. Proposition~\ref{algebraint2}).
The Lie algebra theoretic presentation of the transition from cyclic
deformations to $S^1$-equivariant deformations is the restriction
along the diagonal map $\GG_{\CCC}\to \GG_{\CCC}^{S^1}$
(cf. Remark~\ref{liediagonal}).

\section{Deformations in abstract contexts}

\label{DAC}

We will define the formalism of deformations of an object $M$
in a monoidal $\infty$-category $\MMM^\otimes$.
The formalism will be given in Section~\ref{abstractdeform}.
For this purpose,  we first recall the relative tensor product
of bimodules in Section~\ref{bimodulesection}.

Next, we turn to consider the Koszul duality functors.
In Section~\ref{sectionpairing}, we review the notion of pairing of $\infty$-categories. Under a good condition, the pairing of $\infty$-categories
determines a ``duality functor''.
In Section~\ref{pairing1} and Section~\ref{pairing2},
we will give examples of pairings and duality functors
used in subsequent sections.

\subsection{Tensor products of bimodules}
\label{bimodulesection}
Let $\MMM^\otimes\to \assoc^\otimes$ be a monoidal $\infty$-category.
(All examples of $\MMM$ in this paper come from symmetric monoidal $\infty$-categories.) 
Suppose that it has geometric realizations/colimits of simplicial objects, and 
the tensor product functor $\otimes:\MMM\times \MMM\to \MMM$ preserves geometric
realizations of simplicial objects.
Let $\uni$ be a unit object of $\mathcal{M}$.

Consider the canonical projection $\BMod(\MMM)\stackrel{(\pi_1,\pi_2)}{\to} \Alg_{1}(\MMM)\times \Alg_{1}(\MMM)$
which carries a $A$-$B$-bimodule ${}_AM_B$ to $(A,B)$.
We denote by ${}_A\BMod_B(\MMM)$ the fiber of $(\pi_1,\pi_2)$ over $(A,B)$.
Let
\[
\BMod(\MMM)\times_{\pi_2,\Alg_{1}(\MMM),\pi_1} \BMod(\MMM)\to \BMod(\MMM)
\]
be the relative tensor product functor.
We refer the readers to \cite[4.4.2]{HA} for the construction of relative tensor products.
If ${}_AM_B$ and ${}_BN_C$ are an $A$-$B$-bimodule and a $B$-$C$-bimodule, respectively,
then it sends $({}_AM_B,{}_BN_C)$ to an $A$-$C$-bimodule ${}_AM_B\otimes_B{}_BN_C$
whose underlying object in $\MMM$ is the tensor product $M\otimes_BN$ obtained by bar
construction. 
The forgetful functors induce canonical equivalence $\BMod_{\uni}(\MMM)\stackrel{\sim}{\to} \LMod(\MMM)$, ${}_{\uni}\BMod(\MMM)\stackrel{\sim}{\to}\RMod(\MMM)$,
and ${}_{\uni}\BMod_{\uni}(\MMM)\stackrel{\sim}{\to} \MMM$.
In particular, the relative tensor product functor of bimodules induces
\[
\RMod(\MMM)\times_{\pi_2,\Alg_{1}(\MMM),\pi_1} \LMod(\MMM)\to \MMM
\]
which carries $(M_B,{}_BN)$ to $M\otimes_BN$.
We write $\LMod(\MMM)_{\uni}$ for $\LMod(\MMM)\times_{\MMM}\{\uni\}$
where $\LMod(\MMM)\to \MMM$ is the forgetful functor.
By \cite[4.7.1.40]{HA}, there is an equivalence $\LMod(\MMM)_{\uni}\simeq \Alg_1(\MMM)_{/\uni}=\Alg_1^+(\MMM)$ which commutes with the projection to $\Alg_1(\MMM)$
up to canonical homotopy.
Thus,
the relative tensor product functor determines the ``reduction functor''
\[
\mathsf{r}_{\MMM}:\RMod^+(\MMM):=\RMod(\MMM)\times_{\pi_2,\Alg_{1}(\MMM)} \Alg_1^+(\MMM)\simeq \RMod(\MMM)\times_{\pi_2,\Alg_{1}(\MMM)} \LMod(\MMM)_{\uni}\to \MMM
\]
which carries $(N_B, B, B\to \uni)$ to $N_B\otimes_B\uni$.

\subsection{Abstract deformation functors}
\label{abstractdeform}

We continue to assume that $\MMM^\otimes\to \assoc^\otimes$ is a monoidal $\infty$-category such that
it has geometric realizations/colimits of simplicial objects, and 
the tensor product functor $\otimes:\MMM\times \MMM\to \MMM$ preserves geometric
realizations of simplicial objects.
Let $M$ be an object of $\MMM$.

\begin{Definition}
\label{ADF}
We will define the notion of deformatins of $M$ in a general situation. 
Set $\RMod_M^+(\MMM)=\{M\}\times_{\MMM}\RMod^+(\MMM)$
where $\RMod^+(\MMM)\to \MMM$ is the reduction functor $\mathsf{r}_\MMM$
(cf. Section~\ref{bimodulesection}).
Let $\RMod^+(\MMM)^\dagger$ be
the subcategory spanned by coCartesian morphisms over $\Alg_1^+(\MMM)$ 
whose objects are the same as those of $\RMod^+(\MMM)$.
Consider the left fibration
\[
\mathcal{D}ef_{M}(\MMM):=\{M\}\times_{\MMM}\RMod^+(\MMM)^\dagger \to \Alg_1^+(\MMM).
\]
This left fibration is classified by the functor
$\Def_{M}(\MMM):\Alg_1^+(\MMM)\to \SSS$ informally given by
\[
[B\to \uni]\mapsto \RMod_B(\MMM)^\simeq \times_{\MMM^\simeq}\{M\}
\]
where $\RMod_B(\MMM)\to \MMM$ is given by the reduction functor
that sends $N_B$ to $N_B\otimes_B\uni$.
We refer to $\Def_{M}(\MMM)$ as the deformation functor
of $M$ over $\Alg_1^+(\MMM)$.
We call $\Def_{M}(\MMM)(B\to\uni)$ the space of deformations of $M$ along  $B\to \uni$ (or simply to $B\in \Alg_1(\MMM)$).
For ease of notation, we often
write $\Def_{M}(\MMM)(B)$ for $\Def_{M}(\MMM)(B\to\uni)$.
\end{Definition}

An object of the space of deformations of $M$ to $B\in \Alg_1(\MMM)$
is described as $(M'\in \RMod_B(\MMM),\ B\to \uni,\ M\simeq M'\otimes_{B}\uni)$.

We will define a change of domains of defomation functors.
\begin{Definition}
\label{ADF2}
Consider the context in Definition~\ref{ADF}.
Let $\NNN^\otimes\to \assoc^\otimes$ be another monoidal $\infty$-category such that
it has geometric realizations/colimits of simplicial objects, and 
the tensor product functor $\otimes:\NNN\times \NNN\to \NNN$ preserves geometric realizations of simplicial objects. 
Let $F:\NNN^\otimes\to \MMM^\otimes$ be a monoidal functor.
Let
$\Alg_1^+(\NNN)\to \Alg_1^+(\NNN)$ be the induced functor.
Then the base change $\mathcal{D}ef_{M}(\MMM)\times_{\Alg_1^+(\MMM)}\Alg_1^+(\NNN)\to \Alg_1^+(\NNN)$ of the left fibration is classified by
$\Alg_1^+(\NNN)\to \Alg_1^+(\MMM)\stackrel{\Def_{M}(\MMM)}{\longrightarrow} \SSS$ of $M$ over $\Alg_1^+(\NNN)$.
We will refer to the composite $\Alg_1^+(\NNN) \to \SSS$
as the deformation functor of $M$ over $\Alg_1^+(\NNN)$. 
\end{Definition}

\subsection{}
\label{sectionpairing}
In order to define Koszul duality functors, following \cite[X]{DAG}, we review the notion of pairings.
Let $\CCC$ and $\DDD$ be $\infty$-categories.
Let $\lambda:\mathcal{M} \to \CCC\times \DDD$ be a right fibration.
Following \cite[X, 3.1.1]{DAG}, we refer to the right fibration to $\CCC\times \DDD$ as a pairing of $\infty$-categories $\CCC$ and $\DDD$.
Let $\CCC^{op}\times \DDD^{op} \to \SSS$ be a functor which corresponds to
the right fibration $\lambda$.
We say that $\lambda$ is left representable (resp. right representable)
if $\CCC^{op}\to \Fun(\DDD^{op},\SSS)$  (resp. $\DDD^{op}\to \Fun(\CCC^{op},\SSS)$)
factors as $\CCC^{op}\stackrel{\DD_{\lambda}}{\to} \DDD\stackrel{Y_{\mathcal{D}}}{\hookrightarrow} \Fun(\DDD^{op},\SSS)$
(resp. $\DDD^{op}\stackrel{\DD_{\lambda}'}{\to} \CCC \stackrel{Y_{\mathcal{F}}}{\hookrightarrow} \Fun(\CCC^{op},\SSS)$).
Here $Y_{\mathcal{D}}$ and $Y_{\mathcal{F}}$ are Yoneda embeddings.
We refer to $\DD_{\lambda}$ and $\DD_{\lambda}'$ as the duality functor associated to $\lambda$.
For example, Koszul duality functors can be described/defined by duality functors associated to pairings.

\begin{Lemma}
\label{exchangeduality}
Let $f:\CCC_1\to \mathcal{C}_2$ and $g:\DDD_1\to \mathcal{D}_2$ be functors between $\infty$-categories and
let
\[
\xymatrix{
\mathcal{M}  \ar[r] \ar[d]_{\lambda}   &  \mathcal{N} \ar[d]^{\mu} \\
\mathcal{C}_1\times \mathcal{D}_1 \ar[r]^{f \times g} & \mathcal{C}_2\times \mathcal{D}_2
}
\]
be a commutative diagram in $\wCat$ such that the vertical arrows are pairing of $\infty$-categories.
We refer to the commutative diagram as a morphism of pairings.
Suppose that both $\lambda$ and $\mu$ are left representable. 
The diagram determines a natural transformation
\[
g\circ \DD_{\lambda}\to \DD_{\mu} \circ f.
\]
\end{Lemma}

\Proof
Indeed, consider the functors $S_{\lambda}:\CCC_1^{op}\times \DDD_1^{op}\to \SSS$ and $S_{\mu}:\mathcal{C}_2^{op}\times \mathcal{D}_2^{op} \to \SSS$ corresponding to $\lambda $ and $\mu$.
The diagram induces a natural transformation $S_{\lambda}\to S_{\mu}\circ (f\times g)$.
Let $T_\lambda:\CCC_1^{op}\to \Fun(\DDD_1^{op},\SSS)$ and 
$T_\mu:\mathcal{C}_2^{op}\to \Fun(\mathcal{D}_2^{op},\SSS)$ be functors determined by $S_{\lambda}$
and $S_{\mu}$, respectively.
Let $\psi:\Fun(\mathcal{D}_2^{op},\SSS)\to \Fun(\mathcal{D}_1^{op},\SSS)$
be the functor induced by the composition with $g$.
Let $\phi$ be a left adjoint to $\psi$ (if necessary we replace $\SSS$ by the $\infty$-category
of spaces in a larger universe).
The natural transformation $S_{\lambda}\to S_{\mu}\circ (f\times g)$
induces $T_{\lambda}\to \psi \circ T_{\mu} \circ f$. It gives rise to
$\phi\circ T_{\lambda}\to T_{\mu} \circ f$.
Since the left adjoint $\phi$ is a left Kan extension of $\mathfrak{h}_{\mathcal{D}_2}\circ g:\mathcal{D}_1 \to \mathcal{D}_2\to \Fun(\mathcal{D}_2^{op},\SSS)$, we have $\phi\circ \mathfrak{h}_{\mathcal{D}_1}\simeq \mathfrak{h}_{\mathcal{D}_2}\circ g$.
By the induced equivalence 
$\phi\circ T_{\lambda}\simeq \mathfrak{h}_{\mathcal{D}_2}\circ g\circ \DD_{\lambda}$
and $T_{\mu} \circ f \simeq \mathfrak{h}_{\mathcal{D}_2}\circ \DD_{\mu}\circ f$, 
we obtain $\mathfrak{h}_{\mathcal{D}_2}\circ g\circ \DD_{\lambda}\to  \mathfrak{h}_{\mathcal{D}_2}\circ \DD_{\mu}\circ f $.
\QED

\begin{Remark}
\label{exchangeduality2}
Suppose that we are given a composite of morphisms of pairings
\[
\xymatrix{
\mathcal{M}  \ar[r] \ar[d]_{\lambda}   &  \mathcal{N} \ar[d]^{\mu} \ar[r] & \mathcal{O} \ar[d]^{\nu} \\
\mathcal{C}_1\times \mathcal{D}_1 \ar[r]^{f \times g} & \mathcal{C}_2\times \mathcal{D}_2 \ar[r]^{h\times i} & \mathcal{C}_3\times \mathcal{D}_3.
}
\]
In this situation, morphisms of pairings give rise to natural transformations
$\xi_1:g\circ \DD_{\lambda}\to \DD_{\mu} \circ f$ and $\xi_2:i\circ \DD_{\mu}\to \DD_\nu\circ h$.
Moreover, the big square induces $\kappa:i\circ g\circ \DD_{\lambda}\to \DD_\nu\circ h\circ f$.
By the construction, $\kappa$ factors as
\[
i\circ g\circ \DD_{\lambda}\stackrel{i\circ \xi_1}{\longrightarrow} i\circ \DD_{\mu} \circ f \stackrel{\xi_2 \circ f}\longrightarrow \DD_\nu\circ h\circ f.
\]
\end{Remark}

\subsection{}
\label{pairing1}

We give a first example of pairings and duality functors
which is used in subsequent sections.

\begin{Construction}
\label{smartconstruction}
(i)
Let $\mathcal{M}^\otimes\to \assoc^\otimes$ be  a monoidal $\infty$-category.
Suppose that the underlying $\infty$-category $\mathcal{M}$ has geometric realizations/colimits of simplicial diagrams, and the tensor product functor $\mathcal{M}\times\mathcal{M}\to \mathcal{M}$
preserves geometric realizations of simplicial objects.
By using the $\infty$-category of bimodules, we will define a pairing which determines a Koszul duality functor $\Alg_1^+(\MMM)^{op}\to \Alg^+_1(\MMM)$
in good cases (see Example~\ref{smartexample}).

Let $\uni$ be a unit object of $\mathcal{M}$.
We let $\BMod(\mathcal{M})$ denote the $\infty$-category of bimodules objects in $\mathcal{M}^\otimes$
(see Section~\ref{sectionbackground}, \cite[4.3]{HA}).
There are canonical functors
\[
\BMod(\mathcal{M})\to \LMod(\mathcal{M}),\ \  \ \ \ \  \BMod(\mathcal{M})\to \RMod(\mathcal{M})
\]
which forget the right module structures and left module structures, respectively.
These functors induce
\[
\LMod(\mathcal{M})\times_{\mathcal{M}}\{\uni\} \stackrel{p}{\leftarrow} \BMod(\mathcal{M})\times_{\mathcal{M}}\{\uni\}\stackrel{q}{\to} \RMod(\mathcal{M})\times_{\mathcal{M}}\{\uni\}.
\]
Here projections to $\mathcal{M}$ are forgetful functors.
There exists an endomorphism algebra object of $\uni$, that is
the unit algebra $\uni\in \Alg_{1}(\MMM)$
(we abuse notation by writing $\uni$ for the unit algebra).
If we regard $\uni$ as the endomorphism algebra, there is a canonical equivalence 
$\LMod(\mathcal{M})\times_{\mathcal{M}}\{\uni_{\mathcal{M}}\}\simeq \Alg_{1}(\MMM)_{/\uni}$
(see \cite[4.7.1.40]{HA}).
By reversing module actions in the operadic level, we have an equivalence
$\LMod(\MMM)\stackrel{\sim}{\to} \RMod(\MMM)$ which commutes with $(-)^{op}:\Alg_{1}(\MMM)\stackrel{\sim}{\to}\Alg_{1}(\MMM)$ which carries $A$ to the opposite algebra $A^{op}$ (cf. \cite[4.6.3]{HA}).
The equivalence $\LMod(\MMM)\stackrel{\sim}{\to} \RMod(\MMM)$
carries a left $A$-module $M$ to the right $A^{op}$-module $M$.
It gives rise to
$\RMod(\mathcal{M})\times_{\mathcal{M}}\{\uni_{\mathcal{M}}\}\simeq \LMod(\mathcal{M})\times_{\mathcal{M}}\{\uni_{\mathcal{M}}\}\simeq \Alg_{1}(\MMM)_{/\uni}$
that lies over the equivalence $(-)^{op}:\Alg_{1}(\MMM)\to \Alg_{1}(\MMM)$.
Set $\Alg_{1}^+(\MMM)=\Alg_{1}(\MMM)_{/\uni}$.
We obtain
\[
\phi:\LMod(\mathcal{M})\times_{\mathcal{M}}\{\uni\}\simeq \Alg_1^+(\MMM)\ \ \ \ \ \textup{and}\ \ \ \ \ \ \ \psi:\RMod(\mathcal{M})\times_{\mathcal{M}}\{\uni\}\simeq \Alg_1^+(\MMM)
\]
We obtain
\[
\xymatrix{
p_{\MMM}\times q_{\MMM}:K(\MMM):=\BMod(\mathcal{M})\times_{\mathcal{M}}\{\uni_{\mathcal{M}}\} \ar[rr]^(0.6){(\phi\circ p)\times(\psi\circ q)} &  & \Alg^+_{1}(\MMM)\times \Alg^+_{1}(\MMM).
}
\]
This functor is a right fibration.
The canonical projection $\BMod(\MMM)\to \Alg_{1}(\MMM)\times \Alg_{1}(\MMM)$
is a Cartesian fibration such that a morphism $f$ in $\BMod(\MMM)$ is a Cartesian morphism
exactly when the image of $f$ in $\MMM$ is an equivalence.
The induced functor $\BMod(\mathcal{M})\times_{\mathcal{M}}\{\uni\}\to \Alg_{1}(\MMM)\times \Alg_{1}(\MMM)$ is also a Cartesian fibration.
Let $e:s\to t$ be a morphism in 
$(\LMod(\mathcal{M})\times_{\mathcal{M}}\{\uni\})\times (\RMod(\mathcal{M})\times_{\mathcal{M}}\{\uni\})$
and let $\overline{t}\in \BMod(\MMM)\times_{\MMM}\{\uni\}$ be an object lying over $t$.
We denote by $e'$ the image of $e$ in $\Alg_{1}(\MMM)\times \Alg_{1}(\MMM)$,
and let $\overline{e}':\overline{s}\to \overline{t}$ in $\BMod(\mathcal{M})\times_{\mathcal{M}}\{\uni\}$ 
be a Cartesian morphism
(with respect to the projection to $\Alg_{1}(\MMM)\times \Alg_{1}(\MMM)$) lying over 
$e$.
 We easily see that $\overline{e}'$ is
  a $(p\times q)$-Cartesian morphism.
To verify that $p\times q$ is a right fibration, it will suffice to prove 
that each fiber of $p\times q$ is an $\infty$-groupoid.
Indeed,  by \cite[4.3.2.7, 4.8.4.6, 4.8.5.16]{HA}, for $(A,B)\in \Alg_{1}(\MMM)\times \Alg_{1}(\MMM)$,
there are canonical equivalences
\begin{eqnarray*}
{}_A\BMod_B(\MMM)\simeq \LMod_A(\RMod_B(\MMM))&\simeq&  \LMod_A(\MMM)\otimes_{\MMM}\RMod_B(\MMM)\\
&\simeq& \LMod_A(\MMM)\otimes_{\MMM}\LMod_{B^{op}}(\MMM) \\
&\simeq& \LMod_{A\otimes B^{op}}(\MMM).
\end{eqnarray*}
which commutes with the projections to $\MMM$ up to canonical homotopy.
We have
\[
{}_A\BMod_B(\MMM)\times_{\MMM}\{\uni\}\simeq \LMod_{A\otimes B^{op}}(\MMM)\times_{\MMM}\{\uni\}\simeq \Map_{\Alg_{1}(\MMM)}(A\otimes B^{op},\uni).
\]
If we regard  an object of $(\LMod(\mathcal{M})\times_{\mathcal{M}}\{\uni\}) \times (\RMod(\mathcal{M})\times_{\mathcal{M}}\{\uni\})$
as $(\epsilon_A:A\to \uni, \epsilon_{B^{op}}:B^{op} \to \uni)\in  \Alg_{1}^+(\MMM)\times \Alg^+_{1}(\MMM)$, the fiber of $p\times q$ over $(\epsilon_A,\epsilon_{B^{op}})$ is naturally equivalent to $\Map_{\Alg_{1}(\MMM)}(A\otimes B^{op},\uni)\times_{\Map(A,\uni)\times\Map(B^{op},\uni)}\{(\epsilon_A,\epsilon_{B^{op}})\}$,
that is an $\infty$-groupoid.
As a byproduct of the argument, the composite right fibration
$\BMod(\mathcal{M})\times_{\mathcal{M}}\{\uni\} \to \Alg^+_{1}(\MMM)\times \Alg^+_{1}(\MMM)$
corresponds to the functor $(\Alg^+_{1}(\MMM)\times \Alg^+_{1}(\MMM))^{op}\to \SSS$ informally given by
\[
(\epsilon_A:A\to \uni, \epsilon_B:B \to \uni)\mapsto \Map_{\Alg_{1}(\MMM)}(A\otimes B,\uni)\times_{\Map(A,\uni)\times\Map(B,\uni)}\{(\epsilon_A,\epsilon_B)\}.
\]

(ii)
Let $\NNN^\otimes\to \assoc^\otimes\simeq\eone^\otimes$ be another monoidal $\infty$-category
which has geometric realizations/colimits of simplicial objects, and the tensor product
functor preserves geometric realizations of simplicial objects. 
Let $F:\MMM^\otimes\to \NNN^\otimes$ be a monoidal functor
which preserves geometric realizations of simplicial objects.
The symmetric monoidal functor $F$ gives rise to a commutative diagram
\[
\xymatrix{
\BMod(\MMM)\times_{\MMM}\{\uni\} \ar[r] \ar[d]^{p\times q} & \BMod(\NNN)\times_{\NNN}\{\uni_{\NNN}\} \ar[d]^{p'\times q'} \\
\LMod(\MMM)\times_{\MMM}\{\uni\}\times \RMod(\MMM)\times_{\MMM}\{\uni\} \ar[r]^{f_l\times f_r} & \LMod(\NNN)\times_{\NNN}\{\uni_{\NNN}\}\times \RMod(\NNN)\times_{\NNN}\{\uni_{\NNN}\}
}
\]
where the horizontal arrows are induced by $F$, and the vertical arrows are right fibrations, that is, pairings of $\infty$-categories.
Suppose that both pairings are left representable so that there are duality functors
\[
\DD_{p\times q}:\Alg_{1}^+(\MMM)^{op}\simeq (\LMod(\MMM)\times_{\MMM}\{\uni\})^{op}\to \RMod(\MMM)\times_{\MMM}\{\uni\}\simeq \Alg_{1}^+(\MMM),
\]
\[
\DD_{p'\times q'}:\Alg_{1}^+(\NNN)^{op}\simeq (\LMod(\NNN)\times_{\NNN}\{\uni_{\NNN}\})^{op}\to \RMod(\NNN)\times_{\NNN}\{\uni_{\NNN}\}\simeq \Alg_{1}^+(\NNN).
\]
We write $\DD_{1,\MMM^\otimes}$ and $\DD_{1,\NNN^\otimes}$ for $\DD_{p\times q}$
and $\DD_{p'\times q'}$, respectively.
By Lemma~\ref{exchangeduality}, 
it gives rise to a natural transformation
\[
f_r\circ \DD_{1,\MMM^\otimes}\to \DD_{1,\NNN^\otimes}\circ f_l.
\]
\end{Construction}

\begin{Example}
\label{smartexample}
Let $\mathcal{A}^\otimes$ be a symmetric monoidal presentable $\infty$-category
such that the tensor product functor $\mathcal{A}\times \mathcal{A}\to \mathcal{A}$
preserves small colimits in each variable.
Let $\Alg^\otimes_{n}(\mathcal{A})$ denote the symmetric monoidal $\infty$-category
of $\eenu$-algebra obejcts in $\mathcal{A}$.
The underlying $\infty$-category $\Alg_{n}(\mathcal{A})$ admits sifted colimits (e.g., geometric realizations of simplicial objects), and the tensor product $\Alg_{n}(\mathcal{A})\times \Alg_{n}(\mathcal{A})\to \Alg_{n}(\mathcal{A})$ preserves sifted colimits since
the forgetful functor $\Alg_{n}(\mathcal{A})\to \mathcal{A}$ preserves sifted colimits
(cf. \cite[3,2.3.2]{HA}).
Suppose that $\MMM^\otimes$ is $\Alg^\otimes_{n}(\mathcal{A})$ ($n\ge0$).
By convention, $\mathcal{A}^\otimes=\Alg_0^\otimes(\mathcal{A})$.
Then we can apply Construction~\ref{smartconstruction} to $\MMM^\otimes$.
We obtain
\[
P(\mathcal{M}):=p_{\MMM}\times q_{\MMM}:K(\MMM):=\BMod(\mathcal{M})\times_{\mathcal{M}}\{\uni\}\to \Alg^+_{1}(\MMM)\times \Alg^+_{1}(\MMM).
\]
This pairing is left representable (and right representable).
Let  $\MMM^\otimes=\Alg_{n}^\otimes(\mathcal{A})$ and 
consider  $D:\Alg^+_{1}(\Alg_{n}(\mathcal{A}))^{op}\to\Fun(\Alg^+_{1}(\Alg_{n}(\mathcal{A}))^{op},\SSS)$
determined by $P(\Alg_{n}(\mathcal{A}))$ (cf. Section~\ref{sectionpairing}).
By Dunn additivity theorem \cite{HA}, 
$\Alg^+_{1}(\Alg_{n}(\mathcal{A}))\simeq \Alg^+_{n+1}(\mathcal{A})$.
As discussed in the proof of \cite[X, Lemma 3.1.5]{DAG}, the image of  $f:B\to \uni_{\mathcal{A}} \in \Alg^+_{n+1}(\mathcal{A})$ in 
$\Fun(\Alg^+_{n+1}(\mathcal{A})^{op},\SSS)$ under $D$ is represented by 
a centralizer $Z(f)$ of $f:B\to \uni_{\mathcal{A}}$ with the associated argmentation
$Z(f)\to B\otimes Z(f)\to \uni_{\mathcal{A}}$ (see \cite[5.3.1.15]{HA} for the existence of a centralizer).
This shows that $P(\Alg_n(\mathcal{A}))$ is left representable, which defines
\[
\DD_{n+1,\mathcal{A}^\otimes}:\Alg_{n+1}^+(\mathcal{A})^{op} \to \Alg_{n+1}^+(\mathcal{A}).
\]
For example, when $\mathcal{A}^\otimes=\Mod_A^\otimes$,
$\DD_{n+1,\mathcal{A}^\otimes}$
is $\DD_n:\Alg_{n+1}^+(\Mod_A)^{op}\to \Alg_{n+1}^+(\Mod_A)^{op}$
in Section~\ref{koszuldual}.

\end{Example}

\subsection{Koszul duals from deformation functors}
\label{pairing2}

We define a Koszul dual associated to a deformation.

\begin{Construction}
\label{smartconstruction2}

(i)
Consider the situation in Section~\ref{abstractdeform}.
Given a deformation $M'\in \RMod_B(\MMM)\times_{\MMM}\{M\}$ of $M$,
we define the module action of the Koszul dual of $B$ on $M$.
For this purpose, we first define a pairing of $\mathcal{D}ef_{M}(\MMM)$
and $\{M\}\times_{\MMM} \LMod^+(\MMM):=\{M\}\times_{\MMM} \LMod(\MMM)\times_{\Alg_1(\MMM)}\Alg_1^+(\MMM)$.

Recall $K(\MMM)=\BMod(\MMM)\times_{\MMM}\{\uni\}$ from Example~\ref{smartexample}.
There is the canonical functor $K(\MMM)\to \LMod(\MMM)\times_{\MMM}\{\uni\}\simeq \Alg_1^+(\MMM)$.
Define
\[
\mathcal{P}_M(\MMM):=\mathcal{D}ef_{M}(\MMM)\times_{\Alg_1^+(\MMM)}K(\MMM).
\]

Next, we will define $\mathcal{P}_{M}(\MMM)\to \{M\}\times_{\MMM} \LMod^+(\MMM)$.
We first consider the composite
\begin{eqnarray*}
\mathcal{P}_M(\MMM)=\DEFOR_{M}(\MMM)\times_{\Alg_{1}^+(\MMM)}K(\MMM) &\to& \{M\}\times _{\MMM}\RMod(\MMM)\times_{\Alg_{1}(\MMM)}\BMod(\MMM)\times_{\MMM}\{\uni\} \\
&\to&  \{M\}\times _{\MMM}\RMod(\MMM)\times_{\Alg_{1}(\MMM)}\BMod(\MMM) \\
&\to& \{M\}\times_{\MMM}\RMod(\MMM),
\end{eqnarray*}
where the first functor and the second functor are the forgetful functors, and the third functor is 
determined by the relative tensor product functor
$\RMod(\MMM)\times_{\Alg_{1}(\MMM)}\BMod(\MMM)\to \RMod(\MMM)$
(cf. Section~\ref{bimodulesection}).
Let $\RMod(\MMM)\stackrel{\sim}{\to}\LMod(\MMM)$ be the equivalence
obtained by the reversing module actions in the operadic level. 
The composition with this equivalence yields  $\mathcal{P}_{M}(\MMM)\to \{M\}\times_{\MMM} \LMod(\MMM)$.
If an object of $\mathcal{P}_M(\MMM)$ is defined as the pair of $(M'\in \RMod_{B}(\MMM), M\simeq M'\otimes_{B}\uni)$ and a $B$-$C$-bimodule $\uni$,
then its image in $\{M\}\times_{\MMM} \LMod(\MMM)$
is the left $C^{op}$-module corresponding to the right $C$-module $M'\otimes_{B}\uni$ determined by the right $C$-module $\uni$.
Consider the composite $\mathcal{P}_{M}(\MMM)\to \{M\}\times_{\MMM} \LMod(\MMM)\to \Alg_1(\MMM)$.
This composite is naturally promoted to 
\[
\mathcal{P}_M(\MMM)\to K(\MMM)\to \RMod(\MMM)\times_{\MMM}\{\uni\}\simeq \LMod(\MMM)\times_{\MMM}\{\uni\}\simeq \Alg_1^+(\MMM)
\]
where the first functor and the second functor are projections.
Thus, we have 
\[
\omega_M(\MMM):\mathcal{P}_{M}(\MMM)\to \{M\}\times_{\MMM} \LMod(\MMM)\times_{\Alg_1(\MMM)}\Alg_1^+(\MMM)=\{M\}\times_{\MMM}\LMod^+(\MMM).
\]
The projection $\mathcal{P}_M(\MMM)\to \mathcal{D}ef_M(\MMM)$
and $\omega_M(\MMM)$ determine
\[
P_{M}(\MMM):\mathcal{P}_{M}(\MMM)\to \mathcal{D}ef_M(\MMM)\times(\{M\}\times_{\MMM}\LMod^+(\MMM)).
\]
This is a right fibration, up to an equivalence,
whose promotion of a morphism to a Cartesian morphism
is induced by the promotion of the restriction of scalars.
The fiber of this right fibration over the
 pair $(M', B\to \uni, M\simeq M'\otimes_B\uni)$ and a left $C^{op}$-module $M$ is the space equivalent to
the fiber of the induced morphism
\begin{eqnarray*}
\Def_M(\MMM)(B\to \uni)\times_{\Alg_1(\MMM)^+} \Map_{\Alg_1(\MMM)}(B\otimes C^{op},\uni) \ \ \ \ \ \ \ \ \ \  \ \ \ \ \ \ \ \ \ \ \ \ \ \ \ \ \ \ \ \  \\
\ \ \ \ \ \ \ \ \ \ \ \ \ \ \ \ \ \ \ \ \ \ \ \to \Def_M(\MMM)(B\to \uni)\times (\{M\} \times_{\MMM}\LMod^+(\MMM)\times_{\Alg^+_1(\MMM)}\{C^{op}\to \uni\})
\end{eqnarray*}
over the point determined by the pair (keep in mind that
${}_B\BMod_{C}(\MMM)\times_{\MMM}\{\uni\}\simeq \Map_{\Alg_1(\MMM)}(B\otimes C^{op},\uni)$, cf. Construction~\ref{smartconstruction}). 
When it is left representable pairing of $\infty$-categories,
it gives rise to
\begin{eqnarray*}
\DD_{M,\MMM^\otimes}:\mathcal{D}ef_{M}(\MMM)^{op}\longrightarrow \{M\}\times_{\MMM}\LMod^+(\MMM)
\end{eqnarray*}
(cf. Example~\ref{smartexample2}).

(ii)
Let $\NNN^\otimes\to \assoc^\otimes\simeq\eone^\otimes$ be another monoidal $\infty$-category
which has geometric realizations/colimits of simplicial objects, and the tensor product
functor preserves geometric realizations of simplicial objects. 
Let $F:\MMM^\otimes\to \NNN^\otimes$ be a monoidal functor.
We consider the induced morphism of pairings.
Suppose that that $\RMod^+(\MMM)\to \RMod^+(\NNN)$ induced by $F$
carries coCartesian morphisms over $\Alg_1^+(\MMM)$ to 
coCartesian morphisms over $\Alg_1^+(\NNN)$.
Then $F$ induces $\mathcal{D}ef_M(F):\mathcal{D}ef_M(\MMM)\to \mathcal{D}ef_{F(M)}(\NNN)$.
Moreover, $F$ also induces $K(\MMM)\to K(\NNN)$ lying over $\Alg_1^+(\MMM) \to \Alg_1^+(\NNN)$ in the natural way.
Consequently, $F$ determines $\mathcal{P}_M(\MMM)\to \mathcal{P}_{F(M)}(\NNN)$. Note that $F$ induces $F(M_B)\otimes_{F(B)}F(\uni)\stackrel{\sim}{\to}F(M_B\otimes_B\uni)$ for any $M_B\in\RMod_B(\MMM)$.
Thus, we obtain a morphism of pairings
\[
\xymatrix{
\mathcal{P}_M(\MMM)  \ar[r] \ar[d] &   \mathcal{P}_{F(M)}(\NNN) \ar[d] \\
\mathcal{D}ef_M(\MMM)\times(\{M\}\times_{\MMM}\LMod^+(\MMM))\ar[r] & \mathcal{D}ef_{F(M)}(\NNN)\times(\{F(M)\}\times_{\NNN}\LMod^+(\NNN)).
}
\]
where the lower horizontal functors are determined by $\mathcal{D}ef_M(F)$
and $F_r:\{M\}\times_{\MMM}\LMod^+(\MMM)\to \{F(M)\}\times_{\NNN}\LMod^+(\NNN)$ induced by $F$.
When both pairings are left representable,
it follows from Lemma~\ref{exchangeduality}
that there is a natural transformation
$F_r\circ \DD_{M,\MMM^\otimes}\to \DD_{F(M),\NNN^\otimes}\circ \mathcal{D}ef_M(F)$.

\end{Construction}

\begin{Example}
\label{smartexample2}
Consider the situation in Example~\ref{smartexample}.
That is, $\mathcal{A}^\otimes$ is a symmetric monoidal presentable $\infty$-category
such that the tensor product functor $\mathcal{A}\times \mathcal{A}\to \mathcal{A}$
preserves small colimits in each variable.
Suppose that $\MMM^\otimes$ is 
$\mathcal{A}^\otimes$.
Then the right fibration $P_M(\MMM)$ is left representable
so that it gives rise to
\[
\DD_{M,\MMM^\otimes}:\mathcal{D}ef_{M}(\MMM)^{op}\longrightarrow \{M\}\times_{\MMM}\LMod^+(\MMM).
\]
Indeed, $\DD_{M,\MMM^\otimes}$ sends $(M_B\in \RMod_{B}(\MMM), M_B\otimes_B\uni\simeq M)$ to
 $M \in \LMod_{\DD_{1}(B)}(\MMM)$
 whose left $\DD_1(B)$-module structure is defined as follows.
If we think of $\uni$ as an object of $B$-$\DD_{1}(B)^{op}$-bimodule determined by
$B\otimes\DD_{1}(B)\to \uni$, then the``integral kernel'' $\uni$ defines $\RMod_B\to \RMod_{\DD_{1}(B)^{op}}$
given by $M_B\mapsto M_B\otimes_B\uni$. We can regard $M=M_B\otimes_B\uni$
as a right $\DD_1(B)^{op}$-module, that is, a left $\DD_1(B)$-module.

\end{Example}

\begin{Remark}
We record a simple version of the above construction for the future reference.
Consider the situation in Example~\ref{smartexample}.
Set $\mathcal{P}(\mathcal{A})= \RMod^+(\mathcal{A}) \times_{\Alg_1^+(\mathcal{A})}K(\mathcal{A})$ and let 
$\mathcal{P}(\mathcal{A})\to \RMod^+(\mathcal{A})  \times  \LMod^+(\mathcal{A})$
be a pairing defined similarly to the above.
It gives rise to the duality functor
\[
\xymatrix{
\RMod^+(\mathcal{A})^{op} \ar[r]^(0.5){\DD} & \LMod^+(\mathcal{A})
}
\]
lying over $\DD_{1,\mathcal{A}}:\Alg^+_{1}(\mathcal{A})^{op} \to \Alg^+_{1}(\mathcal{A})$.
\end{Remark}

\section{Deformations of categories and cyclic deformations}
Using the formalism in Section~\ref{DAC},
we define the deformation functors that describe the deformation problem
of a stable $\infty$-category
and that of a module endowed with $S^1$-action.
Also, based on Section~\ref{DAC},
to a deformation of a stable $\infty$-category $\CCC$,
we associate the Koszul dual of the deformation,
that is a certain module structure on $\CCC$.
Similarly, we define the Koszul dual of a deformation of module
endowed with $S^1$-action. These are carried out in Section~\ref{deformationcategorysection}
and Section~\ref{cyclicdeformationmodulesection}.
In Section~\ref{maintranspose}, we study how Hochschild chain functor
$\HH_\bullet(-/A)$ relates two deformation functors and their Koszul
duals. The main result is Proposition~\ref{firstkoszuldual}.

\subsection{Deformations of categories}
\label{deformationcategorysection}
We start with the definition of deformations of stable $\infty$-categories.
Let $\MMM^\otimes=\ST_A^\otimes$.
This symmetric monoidal $\infty$-category is compactly generated,
and the tensor product $\otimes :\ST_A\times\ST_A\to \ST_A$
preserves small colimits separately in each variable (cf. Section~\ref{sectionmodule}).

Let $\CCC$ be an object of $\ST_A$, that is, an $A$-linear (small) stable $\infty$-category.
We  will define the deformation functor of $\CCC$ over $\Alg_1^+(\ST_A)$.
We apply Definition~\ref{ADF} to $\MMM^\otimes=\ST_A^\otimes$. Then
we obtain the left fibration $\mathcal{D}ef_{\CCC}(\ST_A)\to \Alg_1^+(\ST_A)$.
It is classified by the functor $\Def_{\CCC}(\ST_A):\Alg_1^+(\ST_A)
\to \wSSS$. 

Next, we consider the symmetric monoidal fucntor
$\Alg_1(\Mod_A)\to \ST_A$ given by $B\mapsto \Perf_B$ (cf. Section~\ref{sectionmodule}).
It determines $\Alg_2^+(\Mod_A)\stackrel{Dn}{\simeq} \Alg_1^+(\Alg_1(\Mod_A))\to \Alg_1^+(\ST_A)$ where the equivalence $Dn$ follows from Dunn additivity theorem.
We define the composite
\[
\Def_{\CCC}^{\etwo}:\Alg_2^+(\Mod_A)\to \Alg_1^+(\ST_A)\stackrel{\Def_{\CCC}(\ST_A)}{\longrightarrow}  \wSSS
\]
(cf. Definition~\ref{ADF2}).
It sends $[B\to A]\in \Alg_2^+(\Mod_A)$ to the space of deformations
$\{\CCC\}\times_{\ST_A}\RMod_{\Perf_B}(\ST_A)^\simeq$.
	We refer to $\Def_{\CCC}^{\etwo}$ as the $\etwo$-deformation functor of $\CCC$.
An object of $\Def_{\CCC}^{\etwo}(B)$
is described as $(\CCC'\in \RMod_{\Perf_B}(\ST_A),\ \CCC\simeq \CCC'\otimes_{\Perf_B}\Perf_A)$.

Next, we define the Koszul dual arising from a deformation of $\CCC$.
To this end, we apply Construction~\ref{smartconstruction2} and Example~\ref{smartexample2}
to $\MMM^\otimes=\ST_A^\otimes$. 
We have the pairing 
$P_{\CCC}(\ST_A):\mathcal{P}_{\CCC}(\ST_A)\to \mathcal{D}ef_{\CCC}(\ST_A)\times (\{\CCC\}\times_{\ST_A}\LMod^+(\ST_A))$.
Here $\LMod^+(\ST_A)=\LMod(\ST_A)\times_{\Alg_1(\ST_A)}\Alg_1^+(\ST_A)$
and we use the notation in Construction~\ref{smartconstruction2}.
The projection $\mathcal{P}_{\CCC}(\ST_A)\to \mathcal{D}ef_{\CCC}(\ST_A)$, the forgetful functor 
$\BMod(\Alg_1(\Mod_A))\to \LMod(\Alg_1(\Mod_A))$, and $p_{\ST_A}:K(\ST_A)\to \Alg_1^+(\ST_A)$ 
 induce
 \[
 \mathcal{P}_{\CCC}(\ST_A)\times_{K(\ST_A)}K(\Alg_1(\Mod_A)) \to \mathcal{D}ef_{\CCC}(\ST_A)\times_{\Alg_1^+(\ST_A)}\LMod(\Alg_1(\Mod_A))\times_{\Alg_1(\Mod_A)}\{A\}
 \]
  (see Construction~\ref{smartconstruction}
for $p_{\ST_A}$).
Consider
$Q:\mathcal{P}_{\CCC}(\ST_A) \to \Alg_1^+(\ST_A)$
defined as  the composite
$\mathcal{P}_{\CCC}(\ST_A)\stackrel{\textup{forget}}{\to} K(\ST_A)\stackrel{q_{\ST_A}}{\to} \Alg_1^+(\ST_A)$
(see Construction~\ref{smartconstruction}
for $q_{\ST_A}$, keep in mind that it is not same with $\mathcal{P}_{\CCC}(\ST_A)\stackrel{\textup{forget}}{\to} K(\ST_A)\stackrel{p_{\ST_A}}{\to} \Alg_1^+(\ST_A)$).
The morphism $\mathcal{P}_{\CCC}(\ST_A) \to \{\CCC\}\times_{\ST_A}\LMod^+(\ST_A)$
naturally commutes with $Q$ and the forgetful functor. 
 Let $\mathcal{P}_{\CCC}(\ST_A)\times_{Q,\Alg_1^+(\ST_A)}\Alg_1^+(\Alg_1(\Mod_A))$ be the fiber product induced by $Q$. 
Then we have
\begin{eqnarray*}
\mathcal{P}_{\CCC}(\ST_A)\times_{K(\ST_A)}K(\Alg_1(\Mod_A)) &\stackrel{\textup{id}\times_{q_{\ST_A}}q_{\Alg_1(\Mod_A)}}{\longrightarrow}&
\mathcal{P}_{\CCC}(\ST_A)\times_{Q,\Alg_1^+(\ST_A)}\Alg_1^+(\Alg_1(\Mod_A)) \\
&\stackrel{P_{\CCC}(\ST_A)\times_{\textup{id}}Dn}{\longrightarrow}& \{\CCC\}\times_{\ST_A}\LMod^+(\ST_A)\times_{\Alg_1^+(\ST_A)}\Alg_2^+(\Mod_A)
\end{eqnarray*}
Set $\mathcal{D}ef_{\CCC}:=\mathcal{D}ef_{\CCC}(\ST_A)\times_{\Alg_1^+(\ST_A)}\Alg^+_2(\Mod_A)$ and 
\[
\LMod^{\etwo}(A)_{\CCC}:=\{\CCC\}\times_{\ST_A}\LMod^+(\ST_A)\times_{\Alg_1^+(\ST_A)}\Alg_2^+(\Mod_A).
\]
Write $\mathcal{P}_{\CCC}:=\mathcal{P}_{\CCC}(\ST_A)\times_{K(\ST_A)}K(\Alg_1(\Mod_A))$.
We obtain the pairing
\begin{eqnarray*}
P_{\CCC}:\mathcal{P}_{\CCC} \to \mathcal{D}ef_{\CCC} \times 
\LMod^{\etwo}(A)_{\CCC}
\end{eqnarray*}
and the associated duality functor
\[
\DD_{\CCC}:(\mathcal{D}ef_{\CCC})^{op} \to \LMod^{\etwo}(A)_{\CCC}.
\]
Given a deformation of $\CCC$, we call the image under $\DD_{\CCC}$  the Koszul dual of the deformation.
Taking into account Example~\ref{smartexample2}
and $\Alg_2(\Mod_A)\to \Alg_1(\ST_A)$,
we see that the duality functor
$\DD_{\CCC}$ carries
$(\CCC'\in \RMod_{\Perf_B}(\ST_A),\ \CCC\simeq \CCC'\otimes_{\Perf_B}\Perf_A)$
to the left $\Perf_{\DD_{2}(B)}^\otimes$-module $\CCC$ (together with the augmentation
$\DD_2(B)\to A$) defined as follows.
The universal pairing $B\otimes_A\DD_{2}(B)\to A$
induces a morphism $\Perf_{B\otimes_A\DD_2(B)}^\otimes\simeq \Perf_{B}^\otimes\otimes_{\Perf_A}\Perf_{\DD_2(B)}^\otimes\to \Perf_A^\otimes$  in $\Alg_1(\ST_A)$
where $\Perf_{B}^\otimes$, $\Perf_{\DD_2(B)}^\otimes$ and $\Perf_A^\otimes$
are regarded as objects of $\Alg_1(\ST_A)$.
Thus, it determines a structure of $\Perf_{B}^\otimes$-$(\Perf_{\DD_2(B)}^\otimes)^{op}$-bimodule on $\Perf_A$. Here 
$(\Perf_{\DD_2(B)}^\otimes)^{op}$ is the opposite algebra of
 $\Perf_{\DD_2(B)}^\otimes$ in $\Alg_1(\ST_A)$.
Therefore, $\CCC\simeq \CCC'\otimes_{\Perf_{B}}\Perf_A$
inherits a right $(\Perf_{\DD_2(B)}^\otimes)^{op}$-module structure
(i.e., a left $\Perf_{\DD_2(B)}^\otimes$-module)
from that of $\Perf_A$.
From this construction and the fully faithful embedding
$\Alg_2(\Mod_A)\hookrightarrow \Alg_1(\ST_A)$,
the pairing $P_{\CCC}$ is left representable.

\subsection{Cyclic deformations of modules with $S^1$-actions}
\label{cyclicdeformationmodulesection}
Let $\Mod_A^{S^1}=\Fun(BS^1,\Mod_A)$
and let $H$ be an object of $\Mod_A^{S^1}$.
Given $[B\to A] \in \Alg_1^+(\Mod_A^{S^1})$, 
a deformation of $H$ to $B$ is defined as an object of $\RMod_{B}(\Mod_A^{S^1})\times_{\Mod_A^{S^1}}\{H\}$ where $\RMod_{B}(\Mod_A^{S^1})\to \Mod_A^{S^1}$
is the reduction functor (cf. Section~\ref{ADF}).
This deformation problem plays a central role in this paper.

Suppose that $\MMM^\otimes$ is $\Mod_A^{S^1}=\Fun(BS^1,\Mod_A)$
endowed with the symmetric monoidal structure induced by that of $\Mod_A^\otimes$. 
The $A$-module $A$ endowed with the trivial $S^1$-action is a unit object in $\Mod^{S^1}_A$.

Applying Definition~\ref{ADF}
to $\MMM=\Mod_A^{S^1}$,
we have the left fibration
\[
\mathcal{D}ef_H^{\circlearrowleft}:=\mathcal{D}ef_{H}(\Mod_A^{S^1})\to \Alg_1^+(\Mod_A^{S^1}),
\]
which is classified by
$\Def_{H}^{\circlearrowleft}:=\Def_{H}(\Mod_A^{S^1}):\Alg_1^+(\Mod_A^{S^1})\to \SSS$. This deformation functor carries $[B\to A]$
to $\{H\}\times_{\Mod_A^{S^1}}\RMod_B(\Mod_A^{S^1})^\simeq$.
We refer to $\Def_{H}^{\circlearrowleft}(B\to A)$ (simply denoted by $\Def_{H}^{\circlearrowleft}(B)$)
as the space of cyclic deformations of $H$ to $B$.

We apply Construction~\ref{smartconstruction2} 
to $\MMM=\Mod_A^{S^1}$.
Then it gives rise to the pairing
\[
\mathcal{P}_{H}(\Mod_A^{S^1})\to \mathcal{D}ef_H^{\circlearrowleft}\times (\{H\}\times_{\Mod_A^{S^1}}\LMod^+(\Mod_A^{S^1})).
\]
We write $\LMod^+(\Mod_A^{S^1})_H=\{H\}\times_{\Mod_A^{S^1}}\LMod^+(\Mod_A^{S^1})$. Here $\LMod^+(\Mod_A^{S^1})=\LMod(\Mod_A^{S^1})\times_{\Alg_1(\Mod_A^{S^1})}\Alg_1^+(\Mod_A^{S^1})$.
By Example~\ref{smartexample2}, we have the duality functor associated to the pairing
\[
\DD_H^{\circlearrowleft}:=\DD_{H,\Mod_A^{S^1}}: (\mathcal{D}ef_H^{\circlearrowleft})^{op}\to \LMod^+(\Mod_A^{S^1})_H.
\]
Given a deformation of $H$, we call the image under $\DD_{H}^\circlearrowleft$  the Koszul dual of the deformation.
Given a defortmation $(H_B\in \RMod_{B}(\Mod_A^{S^1}), H_B\otimes_BA\simeq H)\in \mathcal{D}ef_H^{\circlearrowleft}(B)$,
the Koszul dual of the deformation
 is defined as follows.
Consider $A$ to be an object of $B$-$\DD_{1}(B)^{op}$-bimodule determined by
$B\otimes\DD_{1}(B)\to A$, then the``integral kernel'' $A$ defines $\RMod_B(\Mod_A^{S^1})\to \RMod_{\DD_{1}(B)^{op}}(\Mod_A^{S^1})$
given by $H_B\mapsto H_B\otimes_BA$. We can regard $H\simeq H_B\otimes_BA$
as a right $\DD_1(B)^{op}$-module, that is, a left $\DD_1(B)$-module in $\Mod_A^{S^1}$.

\subsection{From deformations of categories to cyclic deformations}
\label{maintranspose}
Let $h=\HH_\bullet(-/A):\ST_A^\otimes\to (\Mod_A^{S^1})^{\otimes}$ be the symmetric monoidal Hochschild chain functor (see Section~\ref{hochschildhomologysection}).
We observe that the functor $h$ carries a deformation of an $A$-linear
stable $\infty$-category $\CCC$ to a cyclic deformation
of $\HH_\bullet(\CCC/A)$ (Construction~\ref{transpose}). This is the main motivation for the notion of
cyclic deformations.
Proposition~\ref{firstkoszuldual} provides a comparison between
the image of the Koszul dual of a deformation of $\CCC$ (i.e., the image under $\DD_{\CCC}$,
see Section~\ref{deformationcategorysection}) under $h$ (more precisely,
$\LMod'(h)$, see Construction~\ref{transpose})
and the Koszul dual of the induced cyclic deformation of $\HH_\bullet(\CCC/A)$
(i.e., the image under $\DD_{\HH_\bullet(\CCC/A)}^{\circlearrowleft}$,
see Section~\ref{cyclicdeformationmodulesection}).

\vspace{2mm}

We start with an exchange of the duality functors
arising from $h$.

\begin{Definition}
\label{algebratranspose}
We define an exchange of the duality functors
$\DD_{2,\Mod_A}$ and $\DD_{1,\Mod_A^{S^1}}$, which is induced by $\Alg_1(\Mod_A)\to \ST_A\stackrel{h}{\to} \Mod_A^{S^1}$
(for ease of notation, we write $h$ also for the composite).
By Construction~\ref{smartconstruction}, Example~\ref{smartconstruction} and Dunn additivity theorem, $h$ induces the diagram
\[
\xymatrix{
K(\Alg_{1}(\Mod_A))\ar[r] \ar[d] & K(\Mod_A^{S^1})   \ar[d] \\
\Alg_{2}^+(\Mod_A)\times \Alg_{2}^+(\Mod_A) \ar[r] & \Alg_{1}^+(\Mod_A^{S^1})\times \Alg_{1}^+(\Mod_A^{S^1})
}
\]
which commutes up to canonical homotopy.
The horizontal bottom arrow is determined by the product of
$\Alg_1^+(h):\Alg_{2}^+(\Mod_A)\to \Alg_{1}^+(\Mod_A^{S^1})$
induced by $h$.
Both vertical arrows are left representable pairings so that 
it induces duality functors $\DD_2:=\DD_{2,\Mod_A}:\Alg_{2}^+(\Mod_A)^{op}\to \Alg_{2}^+(\Mod_A)$
and $\DD_{1}^{S^1}:=\DD_{1,\Mod_A^{S^1}}:\Alg^+_{1}(\Mod^{S^1}_A)^{op} \to \Alg^+_{1}(\Mod^{S^1}_A)$.
By Lemma~\ref{exchangeduality},
the diagram gives rise to a natural transformation
\[
\Alg_1^+(h)\circ \DD_2\to \DD_1^{S^1}\circ \Alg_1^+(h)
\]
between functors $\Alg_2^+(\Mod_A)^{op}\to \Alg_1^+(\Mod_A^{S^1})$.
\end{Definition}

\begin{Remark}
\label{algebratranspose2}
Let $B$ be an object of $\Alg_2^+(\Mod_A)$.
Consider $\Alg_1^+(h)\circ \DD_2(B)\to \DD_1^{S^1}\circ \Alg_1^+(h)(B)$.
By the proof of Lemma~\ref{exchangeduality}, this map can be described as follows.
If $B\otimes_A\DD_2(B)\to A$ is a universal pairing,
it gives rise to a morphism
$\HH_\bullet(B/A)\otimes_A\HH_\bullet(\DD_2(B)/A)\to \HH_\bullet(A/A)=A$
in $\Alg_1(\Mod_A^{S^1})$.
The resulting pairing determines $\Alg_1^+(h)\circ \DD_2(B)=\HH_\bullet(\DD_2(B)/A)\to \DD_1^{S^1}(\HH_\bullet(B/A))= \DD_1^{S^1}\circ \Alg_1^+(h)(B)$.
\end{Remark}

\begin{Lemma}
\label{deformationlemma}
For any morphism $R\to R'$ in $\Alg_{2}(\Mod_A)$ and any object $\CCC$ of $\RMod_{\Perf^\otimes_{R}}(\ST_A)$, the canonical map
$\HH_{\bullet}(\CCC/A)\otimes_{\HH_\bullet(\Perf_R/A)}\HH_\bullet(\Perf_{R'}/A)\to \HH_\bullet(\CCC\otimes_{\Perf_R}\Perf_{R'}/A)$
is an equivalence.
\end{Lemma}

\Proof
Note that $\CCC$ can be regarded as a filtered colimit (poset) of full subcategories $\colim_{\lambda}\CCC_\lambda$
such that each $\CCC_{\lambda}$ has a single compact gernerator so that
it is equivalent to $\Perf_{B_\lambda}$ for some $B_{\lambda}\in \Alg_{1}(\Mod_A)$.
Moreover, the functor $\HH_\bullet(-/A)$ preserves filtered colimits
(cf. e.g. \cite[Prop. 10.2]{BGT1}). Thus, we may and will assume that $\CCC=\Perf_B$ such that $B\in \Alg_{1}(\Mod_A)$.
The relative tensor product
$\Perf_B\otimes_{\Perf_R}\Perf_{R'}$ can be obtained by means of bar construction (cf. e.g. \cite[4.4.2]{HA}):
$\Perf_B \otimes_{\Perf_R}\Perf_{R'}$ is a colimit (geometric realization) of a simplicial diagram
$s:\NNNN(\Delta^{op})\to \ST_A$ whose $n$-th term is equivalent to $\Perf_B\otimes_A \Perf_R^{\otimes n}\otimes_A\Perf_{R'}$.
This simplicial diagram can be obtained from a simplicial diagram in $\Alg_{1}(\Mod_A)$
by applying $\Perf_{(-)}$.
Indeed, taking into account the map $B=B\otimes_AA\to B\otimes_RR'$, $s$ is extended to 
a functor $s':\NNNN(\Delta^{op})\to (\ST_A)_{\Perf_B/}$.
By the composition with the functor $\Perf_A \to \Perf_B$ determined by the base change along $A\to B$,
$s$ is promoted to a functor $s'':\NNNN(\Delta^{op})\to(\ST_A)_{\Perf_A/}$.
We note that each $\Perf_A\to \Perf_B\otimes_A \Perf_R^{\otimes n}\otimes_A\Perf_{R'}$
is an essentially unique morphism in $\ST_A$ which carries $A\in \Perf_A$ to $B\otimes_A R^{\otimes n}\otimes_A R'\in
\Perf_{B\otimes_A R^{\otimes n}\otimes_AR'}$.
According to \cite[4.8.5.11]{HA}
and a categorical equivalence $\ST\stackrel{\sim}{\to} \textup{Cgt}_{\textup{St}}^{\textup{L,cpt}}$
(cf. Section~\ref{sectionmodule}), we have
$(\ST_A)_{\Perf_A/}\to  (\PR_A)_{\Mod_A/}\hookleftarrow \Alg_{1}(\Mod_A)$
where the first arrow is induced by $\ST_A\stackrel{\Ind}{\to} \textup{Cgt}_{\textup{St}}^{\textup{L,cpt}}\to \PR_A$,
 and the second left arrow is a fully faithful functor which carries $C$ to $\LMod_A\to \LMod_C$ 
 that sends $A$ to $C$ (cf. \cite[4.8.5]{HA}).
Therefore,
$s'':\NNNN(\Delta^{op})\to(\ST_A)_{\Perf_A/}$ determines a simplicial diagram $\NNNN(\Delta^{op})\to \Alg_{1}(\Mod_A)$ which induces $s''$ (by applying $\Perf$). 
(The $n$-term is equivalent to $B\otimes_A R^{\otimes n}\otimes_A R'\in \Alg_{1}(\Mod_A)$.)
Since $\HH_\bullet(-/A)$ is a symmetric monoidal functor and it is given by bar construction so that it preserves geometric realizations in $\Alg_{1}(\Mod_A)$, we obtain equivalences
\begin{eqnarray*}
\HH_{\bullet}(B/A)\otimes_{\HH_\bullet(R/A)}\HH_\bullet(R'/A) &\simeq& \colim_{[n]\in\Delta^{op}}\HH_\bullet(B/A)\otimes_A\HH_\bullet(R/A)^{\otimes n}\otimes_A\HH_\bullet(R'/A) \\
&\simeq& \colim_{[n]\in \Delta^{op}}\HH_\bullet(B\otimes_A R^{\otimes n}\otimes_AR'/A) \\
&\simeq& \HH_\bullet(B\otimes_RR'/A)\simeq \HH_\bullet(\Perf_B \otimes_{\Perf_R}\Perf_{R'}/A).
\end{eqnarray*}
Thus, the assertion follows.
\QED

\begin{Notation}
\label{simplenotation1}
For ease of notation, we often write $\HHH$ for $\HH_\bullet(\CCC/A)$.
Accordingly, we often use the symbol such as $\End(\HHH)$ instead of $\End(\HH_\bullet(\CCC/A))$.
For ease of notation, we often write $\HHHH$ for $\HH^\bullet(\CCC/A)$.

\end{Notation}

\begin{Construction}
\label{transpose}
We will construct a morphism of pairings defined in Section~\ref{deformationcategorysection}
and Section~\ref{cyclicdeformationmodulesection}.
We suppose that $H=\HH_\bullet(\CCC/A)=\HHH$
in Section~\ref{cyclicdeformationmodulesection}.
It yields an exchange of duality functors $\DD_{\CCC}$
and $\DD_{\HHH}^\circlearrowleft$.
According to Lemma~\ref{deformationlemma},
$h=\HH_\bullet(-/A):\ST_A\to \Mod_A^{S^1}$ induces
$\mathcal{D}ef_{\CCC}(\ST_A)\to \mathcal{D}ef_{\HHH}(\Mod_A^{S^1})$.
Thus, invoking Construction~\ref{smartconstruction} (ii),
we obtain a morphism of pairings
\[
\xymatrix{
\mathcal{P}_{\CCC}(\ST_A) \ar[r] \ar[d] &  \mathcal{P}_{\HHH}(\Mod_A^{S^1}) \ar[d] \\
 \mathcal{D}ef_{\CCC}(\ST_A)\times (\{\CCC\}\times_{\ST_A}\LMod^+(\ST_A)) \ar[r] & \mathcal{D}ef^\circlearrowleft_{\HHH}\times \LMod^+(\Mod_A^{S^1})_{\HHH}
}
\]
(see Section~\ref{deformationcategorysection} and Section~\ref{cyclicdeformationmodulesection} for these pairings).
On the other hand, taking into account
definitions of $\mathcal{P}_{\CCC}$,
$\mathcal{D}ef_{\CCC}$ and $\LMod^{\etwo}(A)_{\CCC}$
described as fiber products (see Section~\ref{deformationcategorysection}), we see that projections induce a morphism
of pairings in the natural way:
\[
\xymatrix{
\mathcal{P}_{\CCC} \ar[r] \ar[d] &  \mathcal{P}_{\CCC}(\ST_A)  \ar[d] \\
 \mathcal{D}ef_{\CCC}\times \LMod^{\etwo}(A)_{\CCC} \ar[r] & \mathcal{D}ef_{\CCC}(\ST_A)\times (\{\CCC\}\times_{\ST_A}\LMod^+(\ST_A)).
}
\]
The composition of two morphisms of pairings gives
\[
\xymatrix{
\mathcal{P}_{\CCC} \ar[r] \ar[d] &  \mathcal{P}_{\HHH}(\Mod_A^{S^1})   \ar[d] \\
 \mathcal{D}ef_{\CCC}\times \LMod^{\etwo}(A)_{\CCC} \ar[r] & \mathcal{D}ef^\circlearrowleft_{\HHH} \times \LMod^+(\Mod_A^{S^1})_{\HHH}.
}
\]
We write $\mathcal{D}ef(h)$ and $\LMod'(h)$
for the induced maps
$\mathcal{D}ef_{\CCC}\to \mathcal{D}ef^\circlearrowleft_{\HHH}$
and 
$\LMod^{\etwo}(A)_{\CCC} \to \LMod^+(\Mod_A^{S^1})_{\HHH}$, respectively.
By Lemma~\ref{exchangeduality},
the composite induces a natural transformation
\[
\LMod'(h) \circ \DD_{\CCC} \to \DD^\circlearrowleft_{\HHH}\circ \mathcal{D}ef(h)
\]
which fills the square

\[
\xymatrix{
(\DEFOR_{\CCC})^{op}\ar[r]^{\mathcal{D}ef(h)} \ar[d]^{\DD_\CCC}  & (\DEFOR_{\HHH}^{\circlearrowleft})^{op}  \ar[d]^{\DD_{\HHH}^\circlearrowleft} \\
 \LMod^{\etwo}(A)_{\CCC} \ar[r]_(0.4){\LMod'(h)}  & \LMod^+(\Mod_A^{S^1})_{\HHH}.
}
\]
\end{Construction}

\begin{Remark}
\label{detailedreduction}
We clarify the relationship
between  $\LMod'(h)\circ \DD_{\CCC}\to \DD_{\HHH}^\circlearrowleft\circ \DEFOR(h)$
and $\Alg_1^+(h)\circ \DD_2\to \DD_{1}^{S^1}\circ \Alg_1^+(h)$ (see Definition~\ref{algebratranspose}).
The consequence is stated in Lemma~\ref{detailedLemma}.
Consider the square of pairings which is described as the commutative diagram
\[
\xymatrix@R=3mm @C=-8mm{
  &  \mathcal{P}_{\CCC} \ar[ld] \ar[rr] \ar[dd]|\hole &   & \mathcal{P}_{\HHH}(\Mod_A^{S^1})  \ar[ld] \ar[dd] \\
 K(\Alg_1(\Mod_A)) \ar[rr] \ar[dd] &  & K(\Mod_A^{S^1})  \ar[dd] &  \\
  &  \DEFOR_{\CCC}\times \LMod^{\etwo}(A)_{\CCC} \ar[ld]_{p_1\times p_2} \ar[rr]^(0.5){\DEFOR(h)\times \LMod'(h)}|(0.47)\hole&   & \DEFORM^\circlearrowleft_{\HHH} \times \LMod^+(\Mod_A^{S^1})_{\HHH} \ar[ld]^(0.4){q_1\times q_2} \\
\Alg_{2}^+(\Mod_A)\times \Alg_{2}^+(\Mod_A) \ar[rr]^{\Alg_1^+(h)\times \Alg_1^+(h)} &  & \Alg_{1}^+(\Mod_A^{S^1})\times \Alg_{1}^+(\Mod_A^{S^1})  &  .
}
\]
in $\wCat$. Here $p_1$, $p_2$, $q_1$ and $q_2$ are canonical projections.
It gives rise to the diagram 
\[
\xymatrix@R=3mm @C=0mm{
&   (\DEFOR_{\CCC})^{op}\ar[rr]^{\mathcal{D}ef(h)} \ar[ld]_{\DD_\CCC} \ar[dd]|\hole&   & (\DEFOR_{\HHH}^{\circlearrowleft})^{op} \ar[dd] \ar[dl]_(0.6){\DD_{\HHH}^\circlearrowleft} \\
 \LMod^{\etwo}(A)_{\CCC} \ar[rr]_(0.55){\LMod'(h)} \ar[dd] & & \LMod^+(\Mod_A^{S^1})_{\HHH} \ar[dd] & \\
&  (\Alg^+_2(\Mod_A))^{op}\ar[rr]^(0.4){\Alg_1^+(h)}|\hole \ar[dl]_{\DD_2}  &  & (\Alg^+_1(\Mod_A^{S^1}))^{op}  \ar[ld]^{\DD_1^{S^1}} \\
 \Alg^+_2(\Mod_A) \ar[rr]^{\Alg_1^+(h)}  & & \Alg^+_1(\Mod_A^{S^1}).&
}
\]
The bottom square is filled by the natural transformation
$\Alg_1^+(h)\circ \DD_2\to \DD_1^{S^1}\circ \Alg_1^+(h)$ 
in Definition~\ref{algebratranspose}.
From the description of $\DD_{\CCC}$ and 
$\DD_{\HHH}^\circlearrowleft$ in Section~\ref{deformationcategorysection} and Section~\ref{cyclicdeformationmodulesection},
we observe that $q_2\circ\DD^\circlearrowleft_{\HHH}\to \DD^{S^1}_{1}\circ q_1$
and $p_2\circ \DD_{\CCC}\to \DD_{2}\circ p_1$ are invertible.
Applying the observation in Remark~\ref{exchangeduality2}
to the above square of pairings, we obtain the diagram
\[
\xymatrix{
q_2\circ \LMod'(h)\circ \DD_{\CCC} \ar[r] \ar[d]_{\simeq} & q_2\circ \DD_{\HHH}^\circlearrowleft\circ \DEFOR(h) \ar[r]^(0.55){\simeq}  & \DD_{1}^{S^1}\circ q_1\circ \DEFOR(h) \ar[d]^{\simeq} \\
\Alg_1^+(h)\circ p_2\circ \DD_{\CCC} \ar[r]^{\simeq} &   \Alg_1^+(h) \circ \DD_{2}\circ p_1 \ar[r]   &   \DD_{1}^{S^1}\circ \Alg_1^+(h)\circ p_1
}
\]
which commutes up to homotopy.
Note that 
$q_2\circ \LMod'(h)\circ \DD_{\CCC} \to q_2\circ \DD^\circlearrowleft_{\HHH}\circ \DEFOR(h)$ is invertible
if and only if
$\LMod'(h)\circ \DD_{\CCC} \to \DD^\circlearrowleft_{\HHH}\circ \DEFOR(h)$ is so.
Thus, taking into account the diagram we conclude that 
$\LMod'(h)\circ \DD_{\CCC} \to \DD^\circlearrowleft_{\HHH}\circ \DEFOR(h)$
is invertible if $\Alg_1^+(h)\circ \DD_{2} \to \DD_{1}^{S^1}\circ \Alg_1^+(h)$ is invertible. 
\end{Remark}

We record the observation:

\begin{Lemma}
\label{detailedLemma}
If $\Alg_1^+(h)\circ \DD_{2} \to \DD_{1}^{S^1}\circ \Alg_1^+(h)$ is invertible, then
$\LMod'(h)\circ \DD_{\CCC} \to \DD^\circlearrowleft_{\HHH}\circ \DEFOR(h)$
is invertible.
\end{Lemma}

\begin{Proposition}
\label{firstkoszuldual}
Consider the composite $\EXT_A\hookrightarrow \CAlg_A^+\stackrel{\textup{forget}}{\to} \Alg_2^+(\Mod_A)$,
the induced base change $(\DEFOR_{\CCC})^{op}\times_{\Alg_{2}^+(\Mod_A)}\EXT_A\to \EXT_A$, and the projection
$v:\DEFOR_{\CCC}\times_{\Alg_{2}^+(\Mod_A)}\EXT_A\to \DEFOR_{\CCC}$.
Then the natural transformation
$\LMod'(h) \circ \DD_{\CCC}\circ v \to \DD^\circlearrowleft_{\HHH}\circ \mathcal{D}ef(h)\circ v$
is a natural equivalence.
Namely, for any $R\in \EXT_A$ and any deformation $\CCC_R\in \ST_R$ of $\CCC$,
 $\LMod'(h) \circ \DD_{\CCC}(\CCC_R)\to \DD_{\HHH}\circ \mathcal{D}ef(h)(\CCC_R)$
 is an equivalence in $\LMod^+(\Mod_A^{S^1})_{\HHH}$.
\end{Proposition}

\Proof
Unfolding the definitions, 
$\LMod'(h)\circ \DD_{\CCC}(\CCC_R)$ is a $\HH_\bullet(\DD_2(R)/A)$-module $\HH_{\bullet}(\CCC/A)$,
and $\DD_{\HHH}\circ \mathcal{D}ef(h)(\CCC_R)$ is a $\DD_1(\HH_\bullet(R/A))$-module
$\HH_\bullet(\CCC/A)$. 
Moreover,
the morphism $\LMod'(h)\circ \DD_{\CCC}(\CCC_R)\to \DD_{\HHH}\circ \mathcal{D}ef(\CCC_R)$ is
the restriction map along $\eta:\HH_\bullet(\DD_2(R)/A)\to \DD_1(\HH_\bullet(R/A))$
which is determined by the pairing $\HH_\bullet(R/A)\otimes\HH_\bullet(\DD_2(R)/A)\to \HH_\bullet(A/A)=A$
induced by the universal pairing $R\otimes_A\DD_2(R)\to A$ 
and $\HH_\bullet(-/A)$ (cf. Definition~\ref{algebratranspose}, Remark~\ref{algebratranspose}, Remark~\ref{detailedreduction}).
By Lemma~\ref{detailedLemma}, it suffices to prove that $\eta$ induces an equivalence of underlying complexes.
We write $\textup{Bar}$ for the bar construction functor $\Alg_{1}(\Mod_A)\to \Alg_{1}(\Mod_A^{op})^{op}$.
Let $(-)^\vee$ denote the $A$-linear dual functor which carries $M$ to $M^\vee$.
Using $\DD_n\simeq (-)^\vee\circ \textup{Bar}^n$,
the map $\eta$ can be identified with the canonical map $\HH_\bullet(\DD_1(\textup{Bar}R)/A)\to (\textup{Bar}\HH_\bullet(R/A))^\vee\simeq \HH_\bullet(\textup{Bar}(R)/A)^{\vee}$ where the latter equivalence
follows from the compatibility of $\HH_\bullet(-/A)$ with respect to sifted colimits
(see the proof of Lemma~\ref{deformationlemma}).
Suppose that $A$ is the field $k$ for the moment.
By applying the standard reduced bar resolution to $R\in \EXT_k$
we see that the complex $\textup{Bar}(R)\simeq k\otimes_Rk$ is connected 
and finite dimensional in each degree.
Under this condition, by duality theorem \cite[4.1.1, 4.1.3]{AF2},
the canonical map $\HH_\bullet(\DD_1(\textup{Bar}R)/k)\to \HH_\bullet(\textup{Bar}(R)/k)^{\vee}$ is an equivalence.
Next, we consider the case when $A$ is an arbitrary connective dg algebra over $k$.
For $R\in \EXT_k$, the standard reduced bar resolution shows also that
$\textup{Bar}^2(R)$ is connected and finite dimensional in each degree.
Moreover, according to
the formula of the standard Hochschild complexes,
$\HH_\bullet(\textup{Bar}(R)/k)$ is connected and finite dimensional in each degree.
By the base change along $k\to A$ and these finiteness properties of $\textup{Bar}^2(R)$ and $\HH_\bullet(\textup{Bar}(R)/k)$,
we deduce the general case from the case of $A=k$.
\QED

\section{Various algebras and Hochschild homology}

In this section, we study a comparison of Koszul duals of Hochschild homology (Proposition~\ref{keycircle}). Combined with 
Proposition~\ref{firstkoszuldual} and Proposition~\ref{UEAProp},
we prove Proposition~\ref{commutativediagram} and Proposition~\ref{eoneetwouniversal}.
The square of equivalences in Proposition~\ref{commutativediagram}
plays an important role in the next section.

\subsection{}
\label{universaleone}
We will construct morphisms $f_C:U_1(\DD_\infty(C))\to \DD_{1}(C)$, which are natural in $C\in \CAlg_A^+$,
and prove that $f_C$ is an equivalence when $C=R\otimes_AS^1$ 
such that $R$ in $\EXT_A$:

\begin{Proposition}
\label{keycircle}
Suppose that $C=R\otimes_AS^1\simeq R\otimes_{R\otimes_AR}R$ such that
$R \in \EXT_A$. Then $f_C:U_1(\DD_\infty(C))\to \DD_{1}(C)$ is an equivalence.
\end{Proposition}

\begin{Definition}
We will define $f_C:U_1(\DD_{\infty}(C))\rightarrow \DD_1(C)$.
Consider the adjoint pair
\[
l: \Fun((\EXT_A)^{op},\Alg_{1}^+(\Mod_A)) \rightleftarrows  \Fun((\CAlg_A^+)^{op},\Alg_{1}^+(\Mod_A)):r
\]
where $r$ is induced by the composition with the inclusion $(\EXT_A)^{op} \to (\CAlg_A^+)^{op}$,
and the left adjoint $l$ carries a functor to its left Kan extension along $(\EXT_A)^{op} \to (\CAlg_A^+)^{op}$
(note that $\Alg_{1}^+(\Mod_A)$ has small colimits).
Then the counit map of ther adjunction applied to
$\DD_{1}:(\CAlg_A^+)^{op} \to \Alg_{1}^+(\Mod_A)$ (we usually omit the forgetful functor
$(\CAlg_A^+)^{op}\to (\Alg_{1}(\Mod_A))^{op}$) gives us a natural transformation
$l\circ r(\DD_{1})\to \DD_{1}$.
For each $C \in \CAlg_A^+$, we write
$f_C:l\circ r(\DD_{1})(C) \to \DD_{1}(C)$
for the morphism obtained by the evaluation.

Recall that $r(\DD_{1})$ is naturally equivalent to
the functor $U_1\circ \DD_\infty:(\EXT_A)^{op}\to \Alg_{1}^+(\Mod_A)$
(cf. Proposition~\ref{UEAProp} (1)).
Thus, for $C\in \CAlg_A^+$, there exist equivalences
\[
l\circ r (\DD_{1})(C)\simeq \varinjlim_{R\in (\EXT_A)^{op}_{/C}}U_1\circ \DD_\infty(R)\simeq \varinjlim_{L\in (Lie_A^f)_{/\DD_{\infty}(C)}}U_1(L) \simeq U_1(\DD_{\infty}(C))
\]
where the second equivalence follows from the equivalence $(\EXT_A)^{op}\simeq Lie_A^f$
and the universality of $\DD_{\infty}(C)$ (determined by the adjoint pair
$(Ch^\bullet,\DD_\infty)$, see Section~\ref{formalstack}), and the third one is an equivalence
since
 $\DD_{\infty}(C)$ is a (sifted) colimit of $(Lie_A^f)_{/\DD_{\infty}(C)}\to Lie_A$.
Thus, we have
\[
f_C:U_1(\DD_{\infty}(C))\longrightarrow \DD_1(C).
\]
\end{Definition}
 
The rest of Section~\ref{universaleone} is devoted to the proof of Proposition~\ref{keycircle}.
We define a closely related map $g_C$.

\begin{Definition}
We will define $g_C$.
Given $C\in \CAlg_A^+$ we consider the morphism $p_C:C\to Ch^\bullet(\DD_\infty(C))$ in $\Alg^+_{1}(\Mod_A)$ determined by the counit map of the adjoint pair $(Ch^\bullet,\DD_\infty)$.
Since there is a natural equivalence $Ch^\bullet\simeq\DD_1\circ  U_1$
 between functors $Lie_A\to \Alg_1^+(\Mod_A)$ (we omit the forgetful functor $\CAlg^+_A\to \Alg_1^+(\Mod_A)$), it follows that $p_C$
can be identified with $C\to \DD_1(U_1(\DD_\infty(C))$,
which is determined  by the pairing
\[
C \otimes_AU_1(\DD_\infty(C))\to Ch^\bullet(\DD_{\infty}(C))\otimes_AU_1(\DD_\infty(C))\simeq \DD_{1}(U_1(\DD_\infty(C)))\otimes U_1(\DD_{\infty}(C))\to A
\]
induced by
$C \to Ch^\bullet(\DD_\infty(C))\simeq \DD_{1}(U_1(C))$
and the universal pairing $\DD_{1}(U_1(\DD_\infty(C)))\otimes U_1(\DD_{\infty}(\DD_\infty(C)))\to A$.
This pairing also determines a morphism $g_C:U_1(\DD_\infty(C)) \to \DD_{1}(C)$.
\end{Definition}

\begin{Lemma}
\label{ftog}
For any $C\in \CAlg_A^+$, $f_C$ is equivalent to $g_C$ in $\Map_{\Alg^+_{1}(\Mod_A)}(U_1(\DD_\infty(C)),\DD_{1}(C))$.
\end{Lemma}

\Proof
For $C,C'\in \CAlg_{A}^+$, there is an equivalence
\[
\Map_{\Alg^+_{1}(\Mod_A)}(U_1(\DD_\infty(C)),\DD_{1}(C'))\simeq \Map_{\Alg^+_{1}(\Mod_A)}(C',\DD_{1}(U_1(\DD_\infty(C)))
\]
given by the adjoint pair $(\DD_{1},\DD_{1})$.
The right hand side is equivalent to
$\Map_{\Alg^+_{1}(\Mod_A)}(C',Ch^\bullet(\DD_\infty(C)))$.
When $C=C'$, by the definitions $p_C:C\to Ch^\bullet(\DD_\infty(C))\simeq \DD_{1}(U_1(\DD_\infty(C))$
corresponds to $g_C$ through equivalences of mapping spaces.
Thus, we will prove that $f_C$ naturally corresponds to $p_C$ through the equivalences.
We note that $f_C$ is given by the morphism $U_1(\DD_{\infty}(C))\simeq \varinjlim_{R\in (\EXT_A)^{op}_{/C}}\DD_{1}(R)\to \DD_{1}(C)$ determined by the canonical morphism $C\to \varprojlim_{R\in (\EXT_A)_{C/}}R$.
Here $\varinjlim_{R\in (\EXT_A)^{op}_{/C}}\DD_{1}(R)\to \DD_{1}(C)$ maps to (the component of)
$C\to \varprojlim_{R\in (\EXT_A)_{C/}}R$ in $ \Map_{\Alg^+_{1}(\Mod_A)}(C,\varprojlim_{R\in (\EXT_A)_{C/}} R)$ through equivalences
\begin{eqnarray*}
\Map_{\Alg_{1}^+(\Mod_A)}( \varinjlim_{R\in (\EXT_A)^{op}_{/C}}\DD_{1}(R),\DD_{1}(C))  &\simeq& \varprojlim_{R\in (\EXT_A)^{op}_{/C}}\Map_{\Alg^+_{1}(\Mod_A)}(\DD_{1}(R),\DD_{1}(C)) \\ 
&\simeq& \Map_{\Alg^+_{1}(\Mod_A)}(C,\varprojlim_{R\in (\EXT_A)_{C/}} \DD_{1}\DD_{1}(R))) \\
&\simeq& \Map_{\Alg^+_{1}(\Mod_A)}(C,\varprojlim_{R\in (\EXT_A)_{C/}} R). \\
\end{eqnarray*}
On the other hand, taking into account
\begin{eqnarray*}
C\to Ch^\bullet(\DD_\infty(C)) \simeq Ch^\bullet\bigl(\varinjlim_{L\in (Lie_A^f)_{/\DD_{\infty}(C)}}L\bigr)  \simeq Ch^\bullet\bigl(\varinjlim_{R\in (\EXT_A)^{op}_{/C}}\DD_\infty(R)\bigr) &\simeq& \varprojlim_{R\in (\EXT_A)_{C/}} Ch^\bullet(\DD_\infty(R))   \\
&\simeq& \varprojlim_{R\in (\EXT_A)_{C/}}R,
\end{eqnarray*}
$p_C$ can be identified with the canonical map
$C\to \varprojlim_{R\in (\EXT_A)_{C/}}R$.
It follows that
$p_C: C\to Ch^\bullet(\DD_\infty(C))$ 
corresponds to
$f_C$
\QED

We will prove Proposition~\ref{keycircle}.
For this purpose, we start with Lemma~\ref{liebasechange}, Corollary~\ref{liebasechange2}
and Remark~\ref{reductiontofield}.

\begin{Lemma}
\label{liebasechange}
Let $L\in Lie_k$ and let $L^{S^1}$ be the cotensor of $L$ by $S^1$.
Denote by $(L\otimes_kA)^{S^1}$ the cotensor in $Lie_A$.
Then there exists a canonical equivalence
$L^{S^1}\otimes_kA\simeq (L\otimes_kA)^{S^1}$.
\end{Lemma}

\begin{Corollary}
\label{liebasechange2}
Let $R_k$ be an object of $\EXT_k$, and $\DD_{\infty,k}(R_k)^{S^1}$ is the cotensor by $S^1$
where $\DD_{\infty,k}:(\CAlg_k^+)^{op} \to Lie_k$ is the Koszul duality functor $\DD_\infty$ in the case of $A=k$.
Let $U_{1,k}:Lie_k\to \Alg_{1}^+(\Mod_k)$ be the universal enveloping algebra functor over $k$.
Then there exists an equivalence
\[
U_1(\DD_{\infty}(R_k\otimes_kA)^{S^1}) \simeq U_{1,k}(\DD_{\infty,k}(R_k)^{S^1})\otimes_kA
\]
in $\Alg_{1}^+(\Mod_A)$.
\end{Corollary}

\Proof
If we put $R_k=k\oplus M$, then $\DD_{\infty,k}(R_k)$
is a free Lie algebra generated by $M^\vee[-1]$
in $Lie_k$, and $\DD_{\infty}(R_k\otimes_kA)$ is a free Lie algebra generated by $(M^\vee\otimes_kA)[-1]$
in $Lie_A$ (see e.g. \cite[X]{DAG}) so that there is a canonical equivalence $\DD_{\infty,k}(R_k)\otimes_kA\stackrel{\sim}{\to} \DD_{\infty}(R_k\otimes_kA)$.
By Lemma~\ref{liebasechange}, $(\DD_{\infty,k}(R_k)\otimes_kA)^{S^1} \simeq  \DD_{\infty,k}(R_k)^{S^1}\otimes_kA$.
Consequently, we have
\[
U_1(\DD_{\infty}(R_k\otimes_kA)^{S^1})\simeq U_1((\DD_{\infty,k}(R_k)\otimes_kA)^{S^1})\simeq U_1(\DD_{\infty,k}(R_k)^{S^1}\otimes_kA)\simeq U_{1,k}(\DD_{\infty,k}(R_k)^{S^1})\otimes_kA.
\]
This completes the proof.
\QED

{\it Proof of Lemma~\ref{liebasechange}.}
It is enough to prove that the canonical map
$L^{S^1}\otimes_kA\to (L\otimes_kA)^{S^1}$
is an equivalence in $\Mod_A$.
Since $\Mod_A$ is stable, it follows that
$(-)^{S^1}:\Mod_k\to \Mod_k$ given by $C\mapsto C^{S^1}=C\times_{C\times C}C$ preserves small colimits (in fact, $(-)^{S^1}$ is equivalent to the
functor given by $C\mapsto C\oplus C[-1]$).
In particular, $L^{S^1}\otimes_kA\to (L\otimes_kA)^{S^1}$
is an equivalence.
\QED

We assume that $A$ is the base field $k$
in Lemma~\ref{koszulcircle}, Lemma~\ref{koszulcircle2}, and Lemma~\ref{fieldfinal}.
See Lemma~\ref{reductiontofield} for the reduction to the case of $A=k$.

\begin{Lemma}
\label{koszulcircle}
For $R\in \EXT_A=\EXT_k$,
we let $u:U_1(\DD_\infty(R\otimes_AS^1)) \to \DD_{1}\DD_{1}(U_1(\DD_\infty(R\otimes_AS^1)))$
be the unit map associated to the adjoint pair $(\DD_1,\DD_1)$.
This morphism is an equivalence in $\Alg^+_{1}(\Mod_A)$.
\end{Lemma}

\Proof
Thanks to \cite[X, 3.1.15]{DAG}, if $U_1(\DD_\infty(R\otimes_kS^1))$ is coconnective 
(i.e., $H_i(\Ker(U_1(\DD_\infty(R\otimes_kS^1))\to A)=0$ for $i\ge 0$) and 
$H_i(U_1(\DD_\infty(R\otimes_kS^1)))$ is finite dimensional  in each degree,
then $u$ is an equivalence. Thus, it is enough to show that
$U_1(\DD_\infty(R\otimes_kS^1))$ satisfies these conditions.
There is a distinguished triangle
\[
\DD_\infty(R)[-1] \to \DD_\infty(R\otimes_kS^1)\simeq \DD_\infty(R)^{S^1}\to \DD_\infty(R)
\]
in the homotopy (triangulated) category of $\Mod_A$. The right morphism is induced by
a point $\ast\to S^1$. Thus, it admits a section $\DD_\infty(R)\to \DD_\infty(R\otimes_kS^1)$
determined by the map $S^1\to \ast$.
It follows that there exists an equivalence $\DD_\infty(R)[-1]\oplus\DD_\infty(R)\simeq  \DD_\infty(R\otimes_kS^1)$ in $\Mod_k$.
The (dg) Lie algebra $\DD_\infty(R)\in Lie_k$ is a
free Lie algebra generated by some perfect complex $N$ in $\Mod_k$
such that $H_i(N)=0$ for $i\ge0$ (cf. \cite[X, Section 2]{DAG}).
Using this description and Poincar\'e-Birkhoff-Witt theorem we see that $U_1(\DD_\infty(R\otimes_kS^1))$
is coconnective, and 
$H_i(U_1(\DD_\infty(R\otimes_kS^1)))$ is finite dimensional in each degree.
\QED

\begin{Lemma}
\label{koszulcircle2}
Suppose that $R$ belongs to $\EXT_A$.
Let $v:R\otimes_AS^1\to Ch^\bullet(\DD_\infty(R\otimes_AS^1))$
be the canonical morphism associated to the adjoint pair $(Ch^\bullet,\DD_\infty)$.
Then $v$ is an equivalence.
In particular, the induced morphism $\DD_1(Ch^\bullet(\DD_\infty(R\otimes_AS^1)))\to \DD_1(R\otimes_AS^1)$
is an equivalence in $\Alg^+_{1}(\Mod_A)$.
\end{Lemma}

\Proof
In this proof, we write $k$ for $A$ since we assume $A=k$.
We denote by $\Mod_k^{\textup{con}}$ the full subcategory of $\Mod_k$
spanned by connective objects, i.e., objects $P$ such that $H_i(P)=0$ for $0>i$.
We denote by $\Mod_k^{(n)}$ the full subcategory of  $\Mod_k^{\textup{con}}$
spanned by
those objects $P$ such that $H_i(P)=0$ for $i>n$.
By definition, $\CCAlg_k=\CAlg(\Mod_k^{\textup{con}})$. 
In this proof, we put $\CAlg_k^{\textup{con}}=\CAlg(\Mod_k^{\textup{con}})$.
Let $\CAlg_k^{(n)}$ be the full subcategory of $\CAlg_k$ spanned by
those objects whose underlying complexes belong to $\Mod_k^{(n)}$.
The fully faithful embedding $\Mod_k^{(n)}\to \Mod_k^{\conn}$
has a left adjoint, that is the truncation functor $\tau_{\le n}:\Mod_k^{\conn}\to \Mod_k^{(n)}$.
The fully faithful functor
$\CAlg_k^{(n)}\to \CAlg_k^{\conn}$ admits a left adjoint $\CAlg_k^{\conn}\to \CAlg_k^{(n)}$
which extends the truncation functor $\tau_{\le n}:\Mod_k^{\conn}\to \Mod_k^{(n)}$
in the natural way.
Thus, there exists an adjoint pair
\[
\tau_{\le n}:\CAlg_k^{\conn}\rightleftarrows \CAlg_k^{(n)}
\]
where we sightly abuse notation by writing $\tau_{\le n}$ for the left adjoint (cf. \cite[Volume 1, Chap.2, 1.2]{Gai2}).
Since the statdard $t$-structure on $\Mod_k$ (determined by the homological amplitudes)
is left complete (see \cite[1.2.1.17]{HA} for this notion), it follows that
there exists a canonical equivalence $R\otimes_kS^1\simeq \varprojlim_{n}\tau_{\le n}(R\otimes_kS^1)$
in $\CAlg_k$. 
We write $A_n$ for $\tau_{\le n}(R\otimes_kS^1)$.
 By \cite[4.4.2.9]{HTT}, $\CAlg_A^+\to \CAlg_A$ preserves sequential limits
 so that
 $R\otimes_kS^1\simeq \varprojlim_{n}A_n$
 is promoted to an equivalence in $\CAlg_k^+$ from that in $\CAlg_k$.
 Next, we will show that each $A_n$ belongs to $\EXT_A$. 
Observe that each homology group of $R\otimes_kS^1$ is finite dimensional.
Indeed, the underlying chain complex of $R\otimes_kS^1$ is equivalent to the standard Hochschild chain complex.
Since $R$ is a trivial square zero extension $k\oplus M$ (such that
$M$ is a connective perfect complex over the field $k$), we deduce from this presentation that 
$H_i(R\otimes_kS^1)$ is a finite dimensional $k$-vector space for any $i\in \ZZ$.
Moreover, $H_0(R)$ is a usual finite dimendional algebra. Since the reduced ring
$H_0(R\otimes_{R\otimes_kR}R)_{\textup{red}}$ is isomorphic to
the underived tensor product $H_0(R)_{\textup{red}}\otimes_{H_0(R)_{\textup{red}}\otimes_kH_0(R)_{\textup{red}}}H_0(R)_{\textup{red}}\simeq k$,
it follows that $H_0(R)$ is a usual artin algebra with residue field $k$.
Thus, each $A_n$ belongs to $\EXT_k$.
Consequently, there is an equivalence $A_n\stackrel{\sim}{\to} Ch^\bullet(\DD_\infty(A_n))$ induced by the adjoint pair
$(Ch^\bullet,\DD_\infty)$.
Next, we will show that 
$\varinjlim_{n}\DD_\infty(A_n) \simeq \DD_\infty(R\otimes_kS^1)$.
By the adjoint pair $(Ch^\bullet,\DD_\infty)$, for $B\in \EXT_k$ there are canonical equivalences
\[
 \Map_{Lie_k}(\DD_\infty(B),\DD_\infty(R\otimes_kS^1))\simeq \Map_{\CAlg_k^+}(R\otimes_kS^1,Ch^\bullet(\DD_\infty(B)))\simeq \Map_{\CAlg_k^+}(R\otimes_kS^1,B).
\]
Thus, the formal stack $\mathcal{F}_{\DD_\infty(R\otimes_kS^1)}$ in $\FST_k$ determined by $\DD_\infty(R\otimes_kS^1)$
can be identified with 
the functor $F:\EXT_k\to \SSS$ given by $B \mapsto \Map_{\CAlg_k^+}(R\otimes_kS^1,B)$.
Similarly, $\mathcal{F}_{\DD_\infty(A_n)}$ can be identified with
the functor $F_n:\EXT_k\to \SSS$ given by $B \mapsto \Map_{\CAlg_k^+}(A_n,B)$.
Suppose that $B$ is the trivial square zero extension $k\oplus k[m]$.
Then
\[
\Map_{Lie_k}(\DD_\infty(B),\DD_\infty(R\otimes_kS^1))\simeq \Map_{Lie_k}(\Free_{Lie}(k[-m-1]),\DD_\infty(R\otimes_kS^1)) \simeq \Omega_{Lie}^\infty(\DD_\infty(R\otimes_kS^1)[m+1])
\]
where $\Free_{Lie}:\Mod_k\to Lie_k$ is the free Lie algebra functor, that is a left adjoint to
the forgetful functor. 
Here $\Omega_{Lie}^\infty$ is the composite
$Lie_k\stackrel{\textup{forget}}{\to} \Mod_k\stackrel{\textup{forget}}{\to} \SP\stackrel{\Omega^\infty}{\to} \SSS$.
The left equivalence is induced by $\Free_{Lie}(k[-m-1])\simeq \DD_\infty(k\oplus k[m])$
(see \cite[X, Section 2]{DAG}).
Similarly, $\Map_{Lie_k}(\DD_\infty(B),\DD_\infty(A_n))\simeq \Omega_{Lie}^\infty(\DD_\infty(A_n)[m+1])$.
To see $\varinjlim_{n}\DD_\infty(A_n) \simeq \DD_\infty(R\otimes_kS^1)$,
it is enough to prove that for any $m\ge0$,
the induced morphism $l_{n,m}: \Omega_{Lie}^\infty(\DD_\infty(A_n)[m+1])\to \Omega_{Lie}^\infty(\DD_\infty(R\otimes_kS^1)[m+1])$
is an equivalence if $n\ge m$.
By equivalences, $l_{n,m}$ can be identified with
\[
l'_{n,m}:\Map_{\CAlg_k^+}(A_n,k\oplus k[m]) \to \Map_{\CAlg_k^+}(R\otimes_kS^1,k\oplus k[m]) 
\]
induced by the composition with $R\otimes_kS^1\to A_n$.
Taking into account the adjoint pair $\tau_{\le n}:\CAlg^{\conn}\rightleftarrows \CAlg^{(n)}$,
 we conclude that $l'_{m,n}$ is an equivalence if $n\ge m$.
 Hence $\DD_\infty(R\otimes_kS^1)\simeq \varinjlim_{n}\DD_\infty(A_n)$.
Now the left adjoint functor $Ch^\bullet:Lie_k\to (\CAlg_k^+)^{op}$ gives rise to equivalences
\[
Ch^\bullet(\DD_{\infty}(R\otimes_kS^1))\simeq Ch^\bullet(\varinjlim_{n}(\DD_\infty(A_n)))\simeq \varprojlim_{n}Ch^\bullet(\DD_\infty(A_n))
\]
in $\CAlg_k^+$.
By the canonical equivalences $A_n\stackrel{\sim}{\to}Ch^\bullet(\DD_\infty(A_n))$, we obtain $ \varprojlim_{n}Ch^\bullet(\DD_\infty(A_n))\simeq \varprojlim_{n}A_n\simeq R\otimes_kS^1$.
This completes the proof.
\QED

\begin{Lemma}
\label{fieldfinal}
The morphism
$g_C:U_1(\DD_\infty(R\otimes_AS^1))\to \DD_1(R\otimes_AS^1)$
is an equivalence.
\end{Lemma}

\Proof
Note that the morphism $g_C$
factors as
\[
U_1(\DD_\infty(R\otimes_AS^1)) \to \DD_{1}\DD_{1}(U_1(\DD_\infty(R\otimes_AS^1)))\simeq \DD_1(Ch^\bullet(\DD_\infty(R\otimes_AS^1)))\to \DD_1(R\otimes_AS^1)
\]
where the left morphism is that of Lemma~\ref{koszulcircle}, the right morphism
is that of Lemma~\ref{koszulcircle2}, and the middle equivalence comes from $\DD_1\circ U_1\simeq Ch^\bullet$.
Our claim follows from Lemma~\ref{koszulcircle} and Lemma~\ref{koszulcircle2}.
\QED

\begin{Lemma}
\label{reductiontofield}
Propostion~\ref{keycircle}
follows from the case when $A$ is $k$.
\end{Lemma}

\Proof
To observe this,
we set $R=R_k\otimes_kA$ (recall that any object $R\in \EXT_A$ is of this form). 
Suppose that $\DD_{1,k}(R_k\otimes_kS^1)$ is coconnective and is finite dimensional in each degree.
Taking into account the description of $\DD_{1,k}(-)$ as the $k$-linear dual $(\textup{Bar}(-))^\vee$,
we see that
there exists a canonical equivalence
$\DD_{1,k}(R_k\otimes_kS^1)\otimes_kA\simeq  \DD_{1}((R_k\otimes_kA)\otimes_AS^1)$
(note also that the noetherian dg algebra $A$ has a bounded amplitude in $\Mod_k$).
Corollary~\ref{liebasechange2} 
allows us to identify $U_1(\DD_{\infty}(R\otimes_AS^1))\to \DD_{1}(R\otimes_AS^1)$
in Proposition~\ref{keycircle} with the composite
\[
U_1(\DD_{\infty}((R_k\otimes_kA)\otimes_AS^1))\simeq U_1(\DD_{\infty}(R_k\otimes_kA)^{S^1}) \simeq U_{1,k}(\DD_{\infty,k}(R_k)^{S^1})\otimes_kA\to \DD_{1,k}(R_k\otimes_kS^1)\otimes_kA
\]
where the final morphism is induced by 
$U_{1,k}(\DD_{\infty,k}(R_k)^{S^1})\to \DD_{1,k}(R_k\otimes_kS^1)$
(that is a morphism in Proposition~\ref{keycircle} in the case of $A=k$).
Lemma~\ref{fieldfinal} shows that
that $U_{1,k}(\DD_{\infty,k}(R_k)^{S^1})\simeq U_{1,k}(\DD_{\infty,k}(R_k\otimes_kS^1))\simeq \DD_{1,k}(R_k\otimes_kS^1)$,
and in the proof of Lemma~\ref{koszulcircle} we prove that
$U_{1,k}(\DD_{\infty,k}(R_k)^{S^1})$ is coconnevive
and is finite dimensional in each degree. 
It follows that $\DD_{1,k}(R_k\otimes_kS^1)$ is coconnective and is finite dimensional in each degree.
Therefore, to verify Propostion~\ref{keycircle},
it will suffice to prove the case of $A=k$. 
\QED

{\it Proof of Proposition~\ref{keycircle}.} (In this proof, $A$ is arbitrary.)
According to Lemma~\ref{reductiontofield}, we may and will assume $A=k$.
By Lemma~\ref{ftog}, we are reduced to proving that $g_C$ is an equivalence.
Lemma~\ref{fieldfinal} shows that $g_C$ is an equivalence.
\QED

\subsection{}

\begin{Proposition}
\label{commutativediagram}
We continue to abuse notation by writing $h$ for both functors
 $\CAlg_A^+\to \CAlg_A^+(\Mod_A^{S^1})$
 and  $\Alg^+_{2}(\Mod_A)\to \Alg_{1}^+(\Mod^{S^1}_A)$
 induced by $\HH_\bullet(-/A):\Alg_{1}(\Mod_A)\to \Mod_A^{S^1}$.
Here we denote simply by $\DD_1^{S^1}\circ h$ or $\DD_{1}\circ h$ the composite
\[
(\EXT_A)^{op}\hookrightarrow (\CAlg_A^+)^{op} \stackrel{h}{\to} (\CAlg_A^+(\Mod_A^{S^1}))^{op} \stackrel{\textup{forget}}{\to} (\Alg_{1}^+(\Mod_A^{S^1}))^{op} \stackrel{\DD_{1}^{S^1}}{\to} \Alg_{1}^+(\Mod_A^{S^1}).
\]
We denote simply by  $h\circ \DD_{2}$
\[
(\EXT_A)^{op}\hookrightarrow (\CAlg_A^+) \stackrel{\textup{forget}}{\to} (\Alg_{2}^+(\Mod_A))^{op} \stackrel{\DD_{2}}{\to} \Alg_{2}^+(\Mod_A)\stackrel{h}{\to} \Alg_{1}^+(\Mod_A^{S^1}).
\]
Then there exist an equivalence
$\delta:h\circ U_2\circ \DD_\infty\stackrel{\sim}{\to}U_1\circ \DD_\infty\circ h$
in $\Fun((\EXT)^{op},\Alg_{1}^+(\Mod_A^{S^1}))$ and the square diagram
of equivalences
\[
\xymatrix{
U_1\circ \DD_{\infty}\circ h \ar[r]_{\simeq}^{\nu}  & \DD^{S^1}_{1}\circ h \\
h\circ U_2\circ \DD_{\infty} \ar[r]_{\simeq}^{\mu} \ar[u]^{\delta}_{\simeq} & h\circ \DD_{2} \ar[u]_{\lambda}^{\simeq} \\
}
\]
that commutes up to canonical homotopy.
\end{Proposition}


\Proof
By Proposition~\ref{firstkoszuldual} and its proof (see Remark~\ref{detailedreduction}),
there is a canonical equivalence
$\lambda:h\circ \DD_{2}\stackrel{\sim}{\to} \DD_1^{S^1}\circ h$
in $\Fun((\EXT_A)^{op},\Alg_{1}^+(\Mod_A^{S^1}))$.
Note that $\DD_{2}:(\EXT_A)^{op}\to\Alg_{2}^+(\Mod_A)$
is naturally equivalent to $U_2\circ \DD_{\infty}:(\EXT_A)^{op}\to \Alg_{2}^+(\Mod_A)$
(see Proposition~\ref{UEAProp}). It follows immediately that there exists a canonical equivalence
$\mu:h\circ U_2\circ \DD_{\infty}\stackrel{\sim}{\to} h\circ \DD_{2}$
in $\Fun((\EXT)^{op},\Alg_{1}^+(\Mod_A^{S^1}))$.
Let $U_1\circ \DD_{\infty}\to \DD_{1}$ be the morphism in $\Fun((\CAlg^+_A)^{op},\Alg_{1}^+(\Mod_A))$
defined in Section~\ref{universaleone}.
It induces a morphism in $\Fun(\Fun(BS^1,\CAlg^+_A)^{op},\Fun(BS^1,\Alg_{1}^+(\Mod_A)))$
for which we write $U_1\circ \DD_{\infty}\to \DD_{1}$ by abuse of notation.
Here we identify with $\CAlg^+(\Mod_A^{S^1}) \stackrel{\textup{forget}}{\to} \Alg_{1}^+(\Mod_A^{S^1})^{op}\stackrel{\DD_1^{S^1}}{\to} \Alg_{1}^+(\Mod_A^{S^1})$ (see
Example~\ref{smartexample} for the notation $\DD_1^{S^1}$) with $\Fun(BS^1,(\CAlg_A^+)^{op})\to \Fun(BS^1,\Alg_{1}^+(\Mod_A))$
induced by $(\CAlg^+_A)^{op} \stackrel{\textup{forget}}{\to} \Alg_{1}^+(\Mod_A)^{op} \stackrel{\DD_1}{\to} \Alg_{1}^+(\Mod_A)$
in the natural way, see Remark~\ref{noharm}.
The composition with
$\EXT_A\hookrightarrow \CAlg_A^+\stackrel{h}{\to} \CAlg^+(\Mod^{S^1}_A)\simeq \Fun(BS^1,\CAlg^+_A)$
(that is equivalent to $\otimes S^1$)
induces a natural transformation
$\nu:U_1\circ \DD_{\infty}\circ h \to \DD_{1}\circ h$
between functors
$(\EXT_A)^{op}\to\Alg_{1}^+(\Mod_A^{S^1})$.
According to Proposition~\ref{keycircle}, $\nu$ is an equivalence in $\Fun((\EXT_A)^{op},\Alg_{1}^+(\Mod_A^{S^1}))$.
We take an inverse $\nu^{-1}$ and let
$\delta$ be the composite $\nu^{-1}\circ \lambda\circ \mu$.
\QED

Using the equivalence $\delta$ in Proposition~\ref{commutativediagram}, we prove the following:
 
\begin{Proposition}
\label{eoneetwouniversal}
Let $L$ be a dg Lie algebra over $A$, that is, an object of $Lie_A$.
Let $L^{S^1}$ denote the cotensor of $L$ by $S^1$.
Then there exists a natural equivalence
$h\circ U_2=\HH_\bullet(-/A)\circ U_2\stackrel{\sim}{\to}U_1\circ (-)^{S^1}$ 
between functors $Lie_A\to \Alg_{1}^+(\Mod_A^{S^1})$.
In particular,
\[
\HH_\bullet(U_2(L)/A)
\stackrel{\sim}{\longrightarrow} U_1(L^{S^1}).
\]
for any $L\in Lie_A$.
\end{Proposition}

\Proof 
Consider the projection
$(Lie_A^f)_{/L}\simeq (\EXT_A)^{op}\times_{Lie_A^f}(Lie_A^f)_{/L}\stackrel{\textup{pr}_1}{\to} (\EXT_A)^{op}$
whose composite is given by the formula $[M\to L] \mapsto Ch^\bullet(M)$.
The composition of this functor and $\delta:h\circ U_2\circ \DD_{\infty}\stackrel{\sim}{\to}U_1\circ \DD_{\infty}\circ h$
induces
\[
\varinjlim_{M\to L\in (Lie_A^f)_{/L}}\HH_\bullet(U_2(M)/A) \stackrel{\sim}{\longrightarrow}    \varinjlim_{M\to L\in (Lie_A^f)_{/L}}U_1(M^{S^1})
\]
where $\DD_{\infty}(Ch^\bullet(M)\otimes S^1)$ is naturally identified with $\DD_\infty(Ch^\bullet(M))^{S^1}\simeq  M^{S^1}$.
Since $L$ is a  sifted colimit of  $(Lie_A^f)_{/L}\to Lie_A$, and both $U_2$ and $\HH_\bullet(-/A)$
preserve sifted colimits, it follows that
there exists a canonical equivalence $\varinjlim_{M\to L\in (Lie_A^f)_{/L}}\HH_\bullet(U_2(M)/A)\simeq \HH_\bullet(U_2(L)/A)$.
Note that the canonical morphism  $\varinjlim_{M\to L\in (Lie_A^f)_{/L}}M^{S^1}\to L^{S^1}$ 
is an equivalence (see Lemma~\ref{separatelemma}).
The functor $U_1$ also preserves colimits so that we have 
$\varinjlim_{M\to L\in (Lie_A^f)_{/L}}U_1(M^{S^1})\simeq U_1(L^{S^1})$.
We conclude that $\HH_\bullet(-/A)\circ U_2$ and $U_1\circ (-)^{S^1}$
are left Kan extenstions of their restrictions $(\HH_\bullet(-/A)\circ U_2)|_{Lie_A^f}$ and $(U_1\circ (-)^{S^1})|_{Lie_A^f}$, respectively.
Thus, $\delta$ induces a natural equivalence
$\HH_\bullet(-/A)\circ U_2\stackrel{\sim}{\to}U_1\circ (-)^{S^1}$.
\QED

\begin{Lemma}
\label{separatelemma}
The canonical morphism  $\varinjlim_{M\to L\in (Lie_A^f)_{/L}}M^{S^1}\to L^{S^1}$ 
is an equivalence.
\end{Lemma}

\Proof
To see this, it is enough to show that the functor $(-)^{S^1}:\Mod_A\to \Mod_A$ given by
$C\mapsto C^{S^1}$ commutes with sifted colimits (because $Lie_A\to \Mod_A$ is conservative
and preserves sifted colimits).
Since $C^{S^1}=C\times_{C\times C}C$ is defined as a finite limit
and $\Mod_A$ is a stable $\infty$-category,
$(-)^{S^1}$ preserves small colimits.
\QED

\begin{Remark}
\label{noharm}
There is no harm to use $\DD_1$ instead of $\DD_1^{S^1}$ because 
it is easy to see that for $B\in \Alg_{1}^+(\Mod_A^{S^1})$
the underlying object $\DD_{1}^{S^1}(B)$ in $\Alg^+_{1}(\Mod_A)$ (which forgets the $S^1$-action)
is naturally equivalent to $\DD_{1}(B)$ when we regard $B$ as an object of $\Alg_{1}^+(\Mod_A)$.
For example, it can be shown by using the description
of $\DD_{1}^{S^1}$ in terms of bar construction or by applying \cite[X, 4.4.18]{DAG} to our situation. 
Moreover, the composite functor $\CAlg^+(\Mod_A^{S^1})^{op} \stackrel{\textup{forget}}{\to} \Alg_{1}^+(\Mod_A^{S^1})^{op}\stackrel{\DD_1^{S^1}}{\to}  \Alg_{1}^+(\Mod_A^{S^1})$
can naturally be identified with the functor
\[
\CAlg^+(\Mod_A^{S^1})^{op}\simeq \Fun(BS^1,\CAlg^+(\Mod_A)^{op})\stackrel{\Fun(BS^1,\DD_1)}{\longrightarrow} \Fun(BS^1,\Alg_{1}^+(\Mod_A))
\]
induced by 
$\CAlg^+(\Mod_A)^{op} \stackrel{\textup{forget}}{\to} \Alg_{1}^+(\Mod_A)^{op}\stackrel{\DD_1}{\to}  \Alg_{1}^+(\Mod_A)$. To see this, by \cite[X, 4.4.20]{DAG}, $\DD_1^{S^1}$ is the 
composite of the bar construction $\textup{Bar}_{\Mod_A^{S^1}}:\Alg_{1}^+(\Mod_A^{S^1})=\Alg_{1}((\Mod_A^{S^1})_{A//A})\to \Alg_{1}(((\Mod_A^{S^1})_{A//A})^{op})^{op}=:\textup{CoAlg}_{1}^+(\Mod^{S^1}_A)$ (see \cite[X, 4.3.1--4.3.3]{DAG}) and the (weak) dual functor $(-)^{\vee}: \textup{CoAlg}^+_{1}(\Mod_A^{S^1})\to \Alg_{1}^+(\Mod_A^{S^1})^{op}$, where $(-)^\vee$ carries each object $M$ to its dual object
in the symmetric monoidal $\infty$-category $\Mod_A^{S^1}$.
($\textup{CoAlg}_{1}^+(\Mod_A^{S^1})$ indicates the $\infty$-category of coaugmented coalgebra
objects in $\Mod_A^{S^1}$.)
Then it suffices to prove that $\CAlg^+(\Mod_A^{S^1}) \stackrel{\textup{forget}}{\to} \Alg_{1}^+(\Mod_A^{S^1}) \to \textup{CoAlg}_{1}^+(\Mod_A^{S^1})$ is equivalent to
\[
\Fun(BS^1,\CAlg^+(\Mod_A)) \to \Fun(BS^1, \Alg_{1}^+(\Mod_A)) \to\Fun(BS^1, \textup{CoAlg}_{1}^+(\Mod_A))
\]
induced by 
$\CAlg^+(\Mod_A)\stackrel{\textup{forget}}{\to}\Alg_{1}^+(\Mod_A)\stackrel{\textup{Bar}_{\Mod_A}}{\to} \Alg_{1}(((\Mod_A)_{A//A})^{op})^{op}=\textup{CoAlg}_{1}^+(\Mod_A)$
where
$\textup{Bar}_{\Mod_A}$ is the bar construction functor.
Let
\[
\textup{Bar}_{\CAlg^{S^1}_A}: \Alg_{1}(\CAlg(\Mod_A^{S^1})_{A//A}) \to \Alg_{1}((\CAlg(\Mod_A^{S^1})_{A//A})^{op})^{op}
\]
be the bar construction for augmented $\eone$-algebra
objects in $\CAlg^{S^1}_A=\CAlg(\Mod_A^{S^1})$.
By Dunn additivity theorem,
there exists a canonical (second) equivalence
$\CAlg^+(\Mod_A^{S^1})\simeq \CAlg(\Mod_A^{S^1})_{A//A}\simeq \Alg_{1}(\CAlg((\Mod^{S^1}_A)_{A//A}))$.
By \cite[X, 4.4.19]{DAG}, there exists a natural equivalence between 
$\textup{Bar}_{\Mod^{S^1}_A}\circ \textup{forget}$ and the composite functor
\[
\CAlg^+(\Mod_A^{S^1})\simeq \Alg_{1}(\CAlg(\Mod_A^{S^1})_{A//A}) \to \Alg_{1}((\CAlg(\Mod_A^{S^1})_{A//A})^{op})^{op} \to  \textup{CoAlg}^+_{1}(\Mod_A^{S^1})
\]
where the middle functor is $\textup{Bar}_{\CAlg^{S^1}_A}$, and the right functor is the forgetful functor
induced by $(\CAlg_A^{S^1})_{A//A}\to (\Mod_A^{S^1})_{A//A}$.
In order to verify our assertion, it is enough to prove that
$\textup{Bar}_{\CAlg^{S^1}_A}$ is equivalent to the functor obtained from the bar construction
$\textup{Bar}_{\CAlg_A}:\Alg_{1}(\CAlg(\Mod_A)_{A//A}) \to \Alg_{1}((\CAlg(\Mod_A)_{A//A})^{op})^{op}$
(defined in a similar way)
by taking $\Fun(BS^1,-)$.
To this end, consider the adjoint pair of bar construction and cobar construction functors
\[
\textup{Bar}_{\CAlg^{S^1}_A}: \CAlg^+(\Mod_A^{S^1}) \rightleftarrows \Alg_{1}((\CAlg^+(\Mod_A^{S^1})^{op})^{op}:\textup{CoBar}_{\CAlg^{S^1}_A}.
\]
Since the symmetric monoidal $\infty$-category $\CAlg^+(\Mod_A^{S^1})$ is 
coCartesian, thus $\Alg_{1}((\CAlg^+(\Mod_A^{S^1})^{op})^{op}$
can be identified with a full subcategory of $\Fun(\Delta,\CAlg^+(\Mod_A^{S^1}))$.
Then if we regard an object of $\Alg_{1}(\CAlg^+(\Mod_A^{S^1})^{op})^{op}$
as a cosimplicial object $c:\Delta\to \CAlg^+(\Mod_A^{S^1})$,
$\textup{CoBar}_{\CAlg^{S^1}_A}$ is given by totalizations/limits,
that is, it carries $c$ to the limit in $\CAlg^+(\Mod_A^{S^1})$.
Therefore, the left adjoint $\textup{Bar}_{\CAlg^{S^1}_A}$ is equivalent to the functor induced by
\[
\CAlg^+(\Mod_A^{S^1}) \to \Fun(\Delta^1, \CAlg(\Mod_A^{S^1}))\to \Fun(\Delta_+,\CAlg(\Mod_A^{S^1}))\to \Fun(\Delta,\CAlg(\Mod_A^{S^1}))
\]
where the left functor is the natural fully faithful functor,
the middle functor is a left adjoint to the restriction functor along $\Delta^1=\{[-1]\to [0]\}\hookrightarrow \Delta_+$,
and the right functor is induced by the restriction along $\Delta\subset \Delta_+$.
 More precisely, $\textup{Bar}_{\CAlg^{S^1}_A}$ can be identified with
$(\CAlg^+(\Mod_A^{S^1}))_{/A}\to\Fun(\Delta,\CAlg(\Mod_A^{S^1})_{/A})$
induced by the composite.
Using this description of $\textup{Bar}_{\CAlg^{S^1}_A}$,
we conclude that $\textup{Bar}_{\CAlg^{S^1}_A}$ is equivalent to the functor obtained from $\textup{Bar}_{\CAlg_A}$ by taking $\Fun(BS^1,-)$.
\end{Remark}

\begin{Remark}
There is another method to prove that
$\HH_\bullet(U_2(L)/A)$ is equivalent to 
$U_1(L^{S^1})$.
The nonabelian Poincar\'e duality for the factorization homology
from Lie algebras and 
the description of $U_n(L)$ as the Chevalley-Eilenberg chain complex $Ch_\bullet(\Omega_0^n(L))$ 
($\Omega_0^n$ indicates the $n$-fold 
based loop) \cite[Propostion 5.13, Remark 5.14]{AF}
show that
\[
\int_{S^1\times\RRR}U_2(L)\simeq Ch_\bullet(\Omega_0(L^{S^1}))\simeq U_1(L^{S^1}).
\]
However,
it is not clear to the author that this equivalence makes the square diagram
in Proposition~\ref{commutativediagram} commute (that is the reason that
I do not directly apply this equivalence).
The commutativity of the square in Proposition~\ref{commutativediagram}
will be used to
prove that the Lie algebra actions constructed from the Hochschild pair
$(\HH^\bullet(\CCC/A),\HH_\bullet(\CCC/A))$ determines
the lower arrow in Theorem~\ref{intromain2}
More specifically, it is necessary for the consistency of identifications
which appear in arguments in Proposition~\ref{algebraint2} and Remark~\ref{modularintrem}.

\end{Remark}

\section{Moduli and Lie algebra actions}

In this section, we prove Theorem~\ref{intromain2}.
We first review an algebraic structure on the Hochschild pair
$(\HH^\bullet(\CCC/A),\HH_\bullet(\CCC/A))$ (cf. Section~\ref{KSalgebrasection}).
The Hochschild cohomology (cochain complex)
$\HH^\bullet(\CCC/A)$ is an $\etwo$-algebra.
It gives rise to the dg Lie algebra $\GG_{\CCC}$
whose underlying complex is equivalent (quasi-isomorphic) to 
the shiftedf Hochschild cohomology $\HH^\bullet(\CCC/A)[1]$.
To $(\HH^\bullet(\CCC/A),\HH_\bullet(\CCC/A))$
we associate the Lie algebra actions of $\GG_{\CCC}$ and $\GG_{\CCC}^{S^1}$
on $\HH_\bullet(\CCC/A)$. It turns out that these actions
determine the morphisms in lower horizontal arrows in Theorem~\ref{intromain2}.
In Section~\ref{leftsquaresection}, we will construct the left square in the diagram in Theorem~\ref{intromain2}. To do this, we will use Proposition~\ref{firstkoszuldual} and the content
of Section~\ref{maintranspose}.
We also interpret the action of $\GG_{\CCC}^{S^1}$
on $\HH_\bullet(\CCC/A)$ as the lower horizontal arrow (see Propotion~\ref{algebraint2}).
In Section~\ref{rightsquaresection}, we will construct and discuss the right square in the diagram in Theorem~\ref{intromain2}.
The Koszul dualities in Proposition~\ref{detailedreduction}
play a pivotal role.
In Section~\ref{finalsubsection}, based on Section~\ref{leftsquaresection} and~\ref{rightsquaresection}, we obtain the main result:
Theorem~\ref{supermain} (and Proposition~\ref{algebraint3}).
The upper horizontal arrows in Theorem~\ref{supermain}
have moduli-theoretic interpretations (see Remark~\ref{modulimeaning1}
and Remark~\ref{modulimeaning2}).

\subsection{}
\label{KSalgebrasection}
We briefly review the algebra of $(\HH^\bullet(\CCC/A),\HH_\bullet(\CCC/A))$
encoded by Kontsevich-Soibelman operad $\KS$.
For details, we refer the reader to \cite{I}
where we construct an algebra $(\HH^\bullet(\CCC/A),\HH_\bullet(\CCC/A))$
over $\KS$. 
We do not recall the operad $\KS$. Instead,
we will give an equivalent formulation which is sutaible for our purpose. 
According to \cite[Theorem 1.2]{I},
a $\KS$-algebra $(\HH^\bullet(\CCC/A),\HH_\bullet(\CCC/A))$ in $\Mod_A$ is equivalent to
giving the following triple:
an $\etwo$-algebra $\HH^\bullet(\CCC/A)\in \Alg_{2}(\Mod_A)$,
an $A$-module with $S^1$-action $\HH_\bullet(\CCC/A)\in \Mod_A^{S^1}$,
and a left $\HH_\bullet(\HH^\bullet(\CCC/A)/A)$-module $\HH_\bullet(\CCC/A)$
in $\Mod_A^{S^1}$ that is an object of $\LMod_{\HH_\bullet(\HH^\bullet(\CCC/A)/A)}(\Mod_A^{S^1})$
lying over $\HH_\bullet(\CCC/A)\in \Mod_A^{S^1}$.
In other words, the data of triple amounts to giving an object of 
\[
\Alg_{2}(\Mod_A)\times_{\Alg_1(\Mod_A^{S^1})}\LMod(\Mod_A^{S^1})
\]
where $\Alg_{2}(\Mod_A) \to \Alg_1(\Mod_A^{S^1})$ is induced by $\HH_\bullet(-/A)$.
We sketch the construction of a left $\HH_\bullet(\HH^\bullet(\CCC/A)/A)$-module $\HH_\bullet(\CCC/A)$.
Recall the adjoint pair $(\theta_A,E_A)$ between $\Alg_2(\Mod_A)$ and $\Alg_1(\PR_A)$
from Section~\ref{hochschildhomologysection}.
Consider the associated counit map
\[
\LMod_{\HH^\bullet(\CCC/A)}=\theta_A\circ E_A(\mathcal{E}nd_A(\Ind(\CCC)))\to \mathcal{E}nd_A(\Ind(\CCC))
\]
in $\Alg_1(\PR_A)$
and the tautological left action of $\mathcal{E}nd_A(\Ind(\CCC))$ on $\Ind(\CCC)$.
We obtain a left $\LMod_{\HH^\bullet(\CCC/A)}$-module $\Ind(\CCC)$
in $\PR_A$.
Passing to full subcategories of compact objects we have
a left $\Perf_{\HH^\bullet(\CCC/A)}$-module $\CCC$ in $\ST_A$.
The map $\LMod(\ST_A)\to \LMod(\Mod^{S^1}_A)$
induced by $\HH_\bullet(-/A)$
gives rise to a left $\HH_\bullet(\HH^\bullet(\CCC/A)/A)$-module $\HH_\bullet(\CCC/A)$.

\begin{Construction}
\label{liealgebraaction}
We construct actions of $\GG_{\CCC}$ and $\GG_{\CCC}^{S^1}$ on $\HH_\bullet(\CCC/A)$ out of 
the left $\HH_\bullet(\HH^\bullet(\CCC/A)/A)$-module $\HH_\bullet(\CCC/A)$,
that is defined as a morphism $U_1(\GG_{\CCC})\to \End(\HHH)^{S^1}$
in $\Alg_1(\Mod_A)$ (see Section~\ref{UEA}).
The counit map of the adjoint pair $(U_2,res_{\etwo/Lie})$
induces $U_2(\GG_{\CCC})=U_2(res_{\etwo/Lie}(\HHHH))\to \HHHH$.
By applying $\HH_\bullet(-/A)$, it gives rise to
$\HH_\bullet(U_2(\GG_{\CCC})/A)\to \HH_\bullet(\HH^\bullet(\CCC/A)/A)=\HH_\bullet(\HHHH/A)$.
Consider the sequence in $\Alg_1(\Mod_A^{S^1})$:
\[
U_1(\GG_{\CCC})\to U_1(\GG_{\CCC}^{S^1})\simeq \HH_\bullet(U_2(\GG_{\CCC}))\to \HH_\bullet(\HHHH/A) \stackrel{A_{\CCC}}{\longrightarrow} \End(\HH_\bullet(\CCC/A))=\End(\HHH)
\]
where the first morphism is determined by the cotensor by $S^1\to \ast$,
the second morphism is the equivalence which comes from Proposition~\ref{eoneetwouniversal},
and
$A_{\CCC}$ corresponds to the left $\HH_\bullet(\HH^\bullet(\CCC/A))$-module $\HH_\bullet(\CCC/A)$.
Since the $S^1$-action on $U_1(\GG_{\CCC})$ is trivial, the composite of the sequence
amounts to
a morphism
\[
A_{\CCC}^{\eone}:U_1(\GG_{\CCC})\longrightarrow \End(\HH_\bullet(\CCC/A))^{S^1}
\]
in $\Alg_1(\Mod_A)$. By the adjoint pair $(U_1,res_{\eone/Lie})$, this morphism
gives us
\[
A_{\CCC}^L:\GG_{\CCC}\longrightarrow \End^L(\HH_\bullet(\CCC/A))^{S^1}
\]
in $Lie_A$.
Here $\End(\HH_\bullet(\CCC/A))$ is the endomorphism
algebra object which is defined as an object of $\Alg_1(\Mod_A^{S^1})$,
and $\End^L(\HH_\bullet(\CCC/A))\in Lie_A^{S^1}$ is the dg Lie algebra obtained from $\End(\HH_\bullet(\CCC/A))\in \Alg_{1}(\Mod_A^{S^1})$.
Let $\End^L(\HH_\bullet(\CCC/A))^{S^1}\in Lie_A$ be the homotopy fixed points of the $S^1$-action.
Similarly, $U_1(\GG_{\CCC}^{S^1})\to  \End(\HH_\bullet(\CCC/A))$ gives rise to
\[
\widehat{A}_{\CCC}^{L}:\GG^{S^1}_{\CCC}\longrightarrow \End^L(\HH_\bullet(\CCC/A))
\]
in $Lie_A^{S^1}=\Fun(BS^1,Lie_A)$.
\end{Construction}

\subsection{}
\label{leftsquaresection}
The main object of this Section~\ref{leftsquaresection} is the left square in Theorem~\ref{intromain2} (see Proposition~\ref{firstsquare}, Remark~\ref{modulimeaning2}, Proposition~\ref{algebraint1}). We begin by
considering right fibrations
$\LMod^{\etwo}(A)_{\CCC}\to \Alg_2^+(\Mod_A)$
and
$\LMod^+(\Mod_A^{S^1})_{\HHH}\to \Alg_1^+(\Mod_A^{S^1})$
defined in Section~\ref{deformationcategorysection}
and Section~\ref{cyclicdeformationmodulesection}.
We observe that in both cases, the functor $\Alg_2^+(\Mod_A)^{op}\to \SSS$
and $\Alg_1^+(\Mod_A^{S^1})^{op} \to \SSS$ corresponding right fibrations
are representable.

\begin{Lemma}
\label{Hochschildrepresentable}
The right fibration
$\LMod^{\etwo}(A)_{\CCC}\to \Alg^{+}_{2}(\Mod_A)$
is equivalent to
\[
\Alg^+_{2}(\Mod_A)\times_{\Alg_{1}(\PR_A)}  \Alg_{1}(\PR_A)_{/\mathcal{E}nd_A(\Ind(\CCC))}\to \Alg^{+}_{2}(\Mod_A)
\]
where $\theta_A:\Alg^+_{2}(\Mod_A) \to \Alg_{2}(\Mod_A) \to \Alg_{1}(\PR_A)$ is obtained by applying
$\Alg_{1}(-)$ to $\Alg_{1}(\Mod_A)\to (\PR_A)_{\Mod_A/}$ informally given by $B\mapsto [\Mod_A\to \LMod_B]$, cf. Section~\ref{sectionmodule}
(note that $\Alg_{2}(\Mod_A)\simeq \Alg_{1}(\Alg_{1}(\Mod_A))$),
and $\mathcal{E}nd_A(\Ind(\CCC))$
is the endomorphism functor category (that is, the endomorphism algebra object of $\Ind(\CCC)$ in $\PR_A$).
Moreover, it is equivalent to the right fibration
\[
\Alg_{2}^+(\Mod_A)_{/A\oplus \HH^\bullet(\CCC/A)}\to \Alg^{+}_{2}(\Mod_A)
\]
where $A\oplus \HH^\bullet(\CCC/A)\to A$ is the augmented Hochschild cohomology complex
associated to the Hochschild cohomology complex $\HH^\bullet(\CCC/A)\in \Alg_{2}(\Mod_A)$.
Namely, if we let $\xi:\Alg^+_{2}(\Mod_A)^{op}\to \SSS$ be a functor that corresponds to $u$,
then $\xi$ is represented by the augmented Hochschild cohomology $A\oplus \HH^\bullet(\CCC/A)\to A$.
\end{Lemma}

\Proof
We first note that $\Ind$-construction provides
a symmetric monoidal equivalence $\ST\stackrel{\sim}{\to} \textup{Cgt}_{\textup{St}}^{\textup{L,cpt}}$
(cf. Section~\ref{sectionmodule}).
Moreover, if $\DDD_B$ is a left module over 
the monoidal $\infty$-category $\LMod^\otimes_B$ (with $B\in \Alg^+_{2}(\Mod_A)$),
by passing to the full subcategories of compact objects,
the module action functor $\Mod_B\times \DDD_B\to \DDD_B$
induces $\Perf_B\times \DDD_B^\omega\to \DDD_B^\omega$.
Therefore, there is an equivalence
\begin{eqnarray*}
\LMod_{\Perf_B}(\ST)\times_{\LMod_{\Perf_A}(\ST)}\{\CCC\} &\simeq&   \LMod_{\LMod_B}(\PR_A)\times_{\PR_A}\{ \Ind(\CCC) \}.
\end{eqnarray*}
By the definition of endomorphsim algebra objects,
the functor $\Alg_{1}(\PR_A)^{op} \to \wSSS$ corresponding to
$\LMod_{\Mod_B}(\PR_A)\times_{\PR_A}\{ \Ind(\CCC) \}\to \Alg_{1}(\PR_A)$
is represented by the endomorphism functor category $\mathcal{E}nd_A(\Ind(\CCC))\in \Alg_{1}(\PR_A)$.
This proves the first assertion.

Next, we will prove the second assertion.
By \cite[4.8.5.11]{HA}, the functor $\theta_A:\Alg_{2}(\Mod_A)\to \Alg_{1}(\PR_A)$ is fully faithful.
By the definition of $\HHHH:=\HH^\bullet(\CCC/A)$ (see Section~\ref{hochschildhomologysection} or \cite[Definition 5.4]{I}), $\HHHH$
is the image of $\mathcal{E}nd_A(\Ind(\CCC))$ under the right adjoint
$\Alg_1(\PR_A)\to \Alg_{2}(\Mod_A)$ of $\theta_A$.
The functor
$\Alg_{2}(\Mod_A)_{/\HHHH}\to \Alg_{1}(\PR_A)_{/\LMod_{\HHHH}}$
induced by $\theta$ is fully faithful.
Moreover, the counit map $\LMod_{\HHHH}\to \mathcal{E}nd_A(\Ind(\CCC))$
induces an equivalence in $\wCat$:
\[
\Alg_{2}(\Mod_A)\times_{\Alg_{1}(\PR_A)} \Alg_{1}(\PR_A)_{/\LMod_{\HHHH}}\to \Alg_{2}(\Mod_A)\times_{\Alg_{1}(\PR_A)} \Alg_{1}(\PR_A)_{/ \mathcal{E}nd_A(\Ind(\CCC))}.
\]
Consequently, we see that $u$ is equivalent to the right fibration
$\Alg^+_{2}(\Mod_A)_{/A\oplus \HHHH}\to \Alg^+_{2}(\Mod_A)$.
This implies the second assertion.
\QED

A similar (but easier) argument also shows:

\begin{Lemma}
\label{endorepresentable}
The right fibration
$\LMod^+(\Mod_A^{S^1})_{\HHH}\to \Alg_1^+(\Mod_A^{S^1})$
is equivalent to
\[
\Alg_1^+(\Mod_A^{S^1})_{/A\oplus \End(\HHH)}\to \Alg_1^+(\Mod_A^{S^1})
\]
where $A\oplus \End(\HHH)$ indicates the augmented 
endomorphism algebra object $A\oplus \End(\HHH)\to A$.
\end{Lemma}

According to Lemma~\ref{Hochschildrepresentable}
and Lemma~\ref{endorepresentable},
$\LMod^{\etwo}(A)_{\CCC}\to \Alg_2^+(\Mod_A)$
and
$\LMod^+(\Mod_A^{S^1})_{\HHH}\to \Alg_1^+(\Mod_A^{S^1})$
are classified by representable functors
$\hhhh_{A\oplus \HHHH}:\Alg_2^+(\Mod_A)^{op}\to \SSS$
and $\hhhh_{A\oplus\End(\HHH)}:\Alg_1^+(\Mod_A^{S^1})^{op}\to\SSS$, respectively.
Consider 
$\LMod'(h):\LMod^{\etwo}(A)_{\CCC} \to \LMod^+(\Mod_A^{S^1})_{\HHH}$
lying over $\Alg_1^+(h):\Alg_2^+(\Mod_A)\to \Alg_1^+(\Mod_A^{S^1})$
(see Remark~\ref{detailedreduction}). By abuse of notation,
we often write $h$ for $\Alg_1^+(h)$.
It gives rise to a natural transformation
\[
\mathfrak{T}^{\etwo,0}_{\CCC}:\hhhh_{A\oplus \HHHH}\longrightarrow \hhhh_{A\oplus\End(\HHH)}\circ h
\]

\begin{Proposition}
\label{algebraint1}
The natural transformation $\mathfrak{T}^{\etwo,0}_{\CCC}$ is determined by the left $\HH_\bullet(\HHHH/A)$-module $\HHH$ in Section~\ref{KSalgebrasection},
that is, a morphism
$A_{\CCC}:\HH_\bullet(\HHHH/A)\to \End(\HHH)$
in $\Alg_1(\Mod_A^{S^1})$ in the natural way (see the proof).
\end{Proposition}

\Proof
By the construction of the left $\HH_\bullet(\HHHH/A)$-module $\HHH$, the left module action of monoidal $\infty$-category
$\Perf_{\HHHH}^\otimes$
on $\CCC$ is induced by the the canonical
monoidal functor $\LMod_{\HHHH}^\otimes\to \End^\otimes(\CCC)$ and the tautological (left) action of $\End(\Ind(\CCC))$ on $\Ind(\CCC)$
(cf. Section~\ref{KSalgebrasection}).
By the restriction, we have a left module action of the $A$-linear monoidal $\infty$-category
$\Perf_{\HHH}^\otimes$ on $\CCC$, which is universal:
given $B\in \Alg_2(\Mod_A)$,
any left action of $\Perf_{B}^\otimes$ on $\CCC$ is uniquely 
given by the monoidal functor $\Perf_{B}^\otimes\to \Perf^\otimes_{\HHHH}$ induced by a unique morphism
$B\to \HHHH$ (namely, $\Map_{\Alg_2(\Mod_A)}(B,\HHHH)\simeq \LMod_{\Perf_B^\otimes}(\ST_A)\times_{\ST_A}\{\CCC\}$, see Lemma~\ref{Hochschildrepresentable}).
Passing to Hochschid chain functor, we have $\HH_\bullet(\HHHH/A)\to \End(\HHH)$ in $\Alg_1(\Mod_A^{S^1})$
which is nothing but a part of data of the algebraic structure
on the Hochschild pair (see Section~\ref{KSalgebrasection}).
In this way, the natural transformation is given by
\begin{eqnarray*}
\hhhh_{A\oplus \HHHH}(B)=\Map_{\Alg_2(\Mod_A)}(B,\HHHH) &\to& \Map_{\Alg_1(\Mod_A^{S^1})}(h(B),h(\HHHH)) \\
&\to&  \Map_{\Alg_1(\Mod^{S^1}_A)}(h(B),\End(\HHH))
\end{eqnarray*}
where the second arrow is given by the composition with
$h(\HHHH)=\HH_\bullet(\HHHH/A)\to \End(\HHH)$.
Finally, taking into account augmentations, the claim follows.
\QED

\begin{Definition}
We define
a natural tranformation $\mathfrak{T}^{\etwo}_{\CCC}:\hhhh_{A\oplus \HHHH}\circ \DD_2\to \hhhh_{A\oplus\End(\HHH)}\circ h\circ \DD_2$
by composing $\mathfrak{T}_\CCC^{\etwo,0}$
with $\DD_2:\Alg_2^+(\Mod_A)^{op} \to\Alg_2^+(\Mod_A)$.
\end{Definition}

\begin{Construction}
\label{fromcube}
Consider the second cube in Remark~\ref{detailedreduction}
which defines a square diagram of right fibrations.
The base changes of right fibrations induce a diagram of right fibrations
over $\Alg_2^+(\Mod_A)^{op}$. It corresponds to the diagram in $\Fun(\Alg_2^+(\Mod_A),\SSS)$:
\[
\xymatrix{
   &  \Def_{\CCC}^{\etwo}\ar[d] \ar[dr]  \ar[r]^{M_{\CCC}^{\etwo,\circlearrowleft}} \ar[ld]_{J_{\CCC}^{\etwo}} &  \Def^{\etwo,\circlearrowleft}(\HHH)  \ar[d]^{J^{\etwo}_{\HHH}} \\
\hhhh_{A\oplus \HHHH}\circ\DD_2  \ar[r]_(0.4){\mathfrak{T}_{\CCC}^{\etwo}} & \mathfrak{h}_{A\oplus\End(\HHH)} \circ h\circ \DD_{2}    &     \mathfrak{h}_{A\oplus\End(\HHH)}  \circ  \DD^{S^1}_{1}\circ h \ar[l]
}
\]
such that each triangle commutes up to canonical homotopy,
where we write $h$ for $\Alg_1^+(h)$. The functor
$\Def^{\etwo,\circlearrowleft}(\HHH)$ corresponds to
$\mathcal{D}ef_{H}^{\circlearrowleft}\times_{(\Alg_1^+(\Mod_A)^{S^1})^{op}}\Alg_2^+(\Mod_A)^{op}\to \Alg_2^+(\Mod_A)^{op}$.
If we take the restriction along the forgetful functor
$\EXT_A\to \Alg_2^+(\Mod_A)$,
it follows from Proposition~\ref{firstkoszuldual} that
$\mathfrak{h}_{A\oplus\End(\HHH)}  \circ  \DD_{1}^{S^1}\circ h\to 
\mathfrak{h}_{A\oplus\End(\HHH)} \circ h\circ \DD_{2}$ 
becomes an equivalence after the restriction.
Let us consider the diagram in $\Fun(\EXT_A,\SSS)$ induced by the restriction
along $\EXT_A\to \Alg_2^+(\Mod_A)$:
\[
\xymatrix{
   &  \Def_{\CCC} \ar[d] \ar[dr]  \ar[r]^{M_{\CCC}^{\circlearrowleft}} \ar[ld]_{J_{\CCC}} &  \Def^{\circlearrowleft}(\HHH)  \ar[d]^{J_{\HHH}^\circlearrowleft} \\
\hhhh_{A\oplus \HHHH}\circ\DD_2|_{\EXT_A}  \ar[r]_(0.4){\mathfrak{T}_{\CCC}} & \mathfrak{h}_{A\oplus\End(\HHH)} \circ h\circ \DD_{2}|_{\EXT_A}    &     \mathfrak{h}_{A\oplus\End(\HHH)}  \circ  \DD^{S^1}_{1}\circ h|_{\EXT_A} \ar[l]^{\simeq}
}
\]
where $(-)|_{\EXT_A}$ indicates the restriction along 
$\EXT_A\to \Alg_2^+(\Mod_A)$, and $\Def_{\CCC}=\Def_{\CCC}^{\etwo}|_{\EXT_A}$ and $\Def^{\circlearrowleft}(\HHH)=\Def^{\etwo,\circlearrowleft}(\HHH)|_{\EXT_A}$.
\end{Construction}

\begin{Remark}
\label{modulimeaning1}
Let us consider the moduli-theoretic meaning of $M_{\CCC}^\circlearrowleft$.
By definition,
$\Def^{\circlearrowleft}(\HHH):\EXT_A\to \SSS$
is given by
\[
\EXT_A\ni [R\to A]\mapsto \RMod_{\HH_\bullet(R/A)}(\Mod_A^{S^1})^\simeq \times_{(\Mod^{S^1}_A)^{\simeq}}\{\HHH\}\in \SSS.
\]
By Lemma~\ref{commutativemodulitransform2}, $\EXT_A\to \Alg_1(\Mod_A^{S^1})$ given by $h|_{\EXT_A}$
is equivalent to the functor determined by the tensor with $S^1$, that is,
the functor given by $R\mapsto R\otimes_AS^1$.
Namely, the space
$\Def^{\circlearrowleft}(\HHH)(R)$
parametrizes cyclic deformations of $\HH_\bullet(\CCC/A)$
to $R\otimes_AS^1$ (cf. Section~\ref{cyclicdeformationmodulesection}).
We refer to $\Def^{\circlearrowleft}(\HHH)$
to the deformation functor of cyclic deformations of $\HHH$
over $\EXT_A$.
Let $\overline{\CCC}=(\CCC'\in \RMod_{\Perf_R^\otimes}(\ST_A),\ R\to A,\ \CCC\simeq \CCC'\otimes_{\Perf_R}\Perf_A)$ be a deformation of $\CCC$ to
$R$, that is, an object of $\Def_{\CCC}(R)$. Then $M_{\CCC}^\circlearrowleft$
sends $\overline{\CCC}$
to the cyclic deformation of $\HH_\bullet(\CCC/A)$:
\[
(\HH_\bullet(\CCC'/A)\in \RMod_{R\otimes_AS^1}(\Mod_A^{S^1}),\ R\otimes_AS^1=\HH_\bullet(R/A)\to A,\ \HHH\simeq\HH_\bullet(\CCC'/A)\otimes_{R\otimes_AS^1}A).
\]
where the final equivalence follows from Lemma~\ref{deformationlemma}.
\end{Remark}

\vspace{2mm}

\begin{Notation}
\label{cyclicliefunctor}
We write $\mathcal{F}_{A\oplus \End(\HH_\bullet(\CCC/A))}^\circlearrowleft:\EXT_A\to \SSS$
for $\mathfrak{h}_{A\oplus\End(\HH_\bullet(\CCC/A))}  \circ  \DD_{1}^{S^1}\circ h|_{\EXT_A}$.
For ease of notation we usually write $\mathcal{F}_{A\oplus \End(\HHH)}^\circlearrowleft$ for $\mathcal{F}_{A\oplus \End(\HH_\bullet(\CCC/A))}^\circlearrowleft$.
\end{Notation}

\begin{Notation}
\label{liecatdef}
We write $\mathcal{F}_{A\oplus \HHHH}$
for $\hhhh_{A\oplus \HHHH}\circ\DD_2|_{\EXT_A}$.
\end{Notation}

\begin{Notation}
\label{GCdef}
We write $\GG_{\CCC}$ for the image $res_{\etwo/Lie}(\HHHH)=res_{\etwo/Lie}(\HH^\bullet(\CCC/A))$
of $\HH^\bullet(\CCC/A)$
under the right adjoint $\Alg_{2}(\Mod_A)\to Lie_A$.
Let $\mathcal{F}_{\GG_{\CCC}}$ denote the formal stack associated to
$\GG_{\CCC}$, that is given by
the functor given by $R\mapsto \Map_{Lie_A}(\DD_\infty(R),\GG_{\CCC})$
 (cf. Section~\ref{formalstack}).
\end{Notation}

\begin{Lemma}
\label{e2lie}
There is a canonical equivalence
$\mathcal{F}_{\GG_{\CCC}}\simeq \mathcal{F}_{A\oplus \HHHH}$.
\end{Lemma}

\Proof
This follows from the definition of $\GG_{\CCC}$ (see Section~\ref{UEA}).
\QED

\begin{Proposition}
\label{firstsquare}
We use the notation in Construction~\ref{fromcube}, Notation~\ref{cyclicliefunctor}, \ref{liecatdef}, \ref{GCdef}, and Lemma~\ref{e2lie}.
Then we have the diagram in $\Fun(\EXT_A,\SSS)$:
\[
\xymatrix{
\Def_{\CCC} \ar[r]^(0.4){M_{\CCC}^{\circlearrowleft}} \ar[d]^{J_{\CCC}} & \Def^{\circlearrowleft}(\HHH)  \ar[d]^{J^\circlearrowleft_{\HHH}}  \\
\mathcal{F}_{\GG_{\CCC}}  \ar[r]^(0.35){\mathfrak{T}_{\CCC}} & \mathcal{F}_{A\oplus \End(\HHH)}^\circlearrowleft.
}
\]
which commutes up to canonical homotopy.
By abuse of notation, the lower horizontal arrow denotes the composite
$\mathcal{F}_{\GG_{\CCC}}\simeq
\hhhh_{A\oplus \HH^\bullet(\CCC)}\circ\DD_2|_{\EXT_A}  \stackrel{\mathfrak{T}_{\CCC}}{\to} \mathfrak{h}_{A\oplus\End(\HHH)} \circ h\circ \DD_{2}|_{\EXT_A}  \simeq     \mathfrak{h}_{A\oplus\End(\HHH)}  \circ  \DD_{1}\circ h|_{\EXT_A}=\mathcal{F}_{A\oplus \End(\HHH)}^\circlearrowleft$.
\end{Proposition}

\Proof
Combine Construction~\ref{fromcube} and Lemma~\ref{e2lie}.
\QED

The lower horizontal arrow 
$\mathcal{F}_{\GG_{\CCC}}  \to \mathcal{F}_{A\oplus \End(\HHH)}^\circlearrowleft$
is induced by $\mathfrak{T}_{\CCC}$ (determined by the algebra of Hochschild pair) and the duality $\mathfrak{h}_{A\oplus\End(\HHH)} \circ h\circ \DD_{2}|_{\EXT_A}  \simeq \mathfrak{h}_{A\oplus\End(\HHH)}  \circ  \DD_{1}\circ h|_{\EXT_A}$ (cf. Proposition~\ref{algebraint1}). 
To pursue the relationship with the algebraic structure of the Hochschild pair,
we prove Proposition~\ref{algebraint2} below.

\begin{Lemma}
\label{modularint}
Let $E^L$ be an object of $Lie_A^{S^1}$.
Let $M$ be a dg Lie algebra, that is, $M\in Lie_A$.
We define $\mathcal{F}^{\circlearrowleft}_{E^L}:(Lie_A^{f})^{op}\to \SSS$
by the formula $L\mapsto \Map_{Lie_A^{S^1}}(L^{S^1},E^L)$.
Then there exists a canonical equivalence of spaces
\[
\Map_{Lie_A^{S^1}}(M^{S^1}, E^L)\simeq \Map_{\Fun(\EXT_A,\SSS)}(\mathcal{F}_{M}, \mathcal{F}^\circlearrowleft_{E^L}).
\]
Given $M^{S^1}\to E^L$, the corresponding morphism $\mathcal{F}_{M}\to \mathcal{F}^\circlearrowleft_{E^L}$ in $\Fun((Lie_A^f)^{op},\SSS)$ is given by
\[
\mathcal{F}_M(L)=\Map_{Lie_A}(L,M)\to \Map_{Lie_A^{S^1}}(L^{S^1},M^{S^1}) \to \Map_{Lie_A^{S^1}}(L^{S^1}, E^L)=\mathcal{F}^\circlearrowleft_{E^L}(L)
\]
where the first arrow is induced by the cotensor with $S^1$,
and the second arrow is induced by the composition with $M^{S^1}\to E^L$.
\end{Lemma}

\begin{Remark}
\label{modularintrem}
Suppose that $E^L=\End^L(\HHH)$.
Then there exists an equivalence $\mathcal{F}^{\circlearrowleft}_{E^L}\simeq \mathcal{F}_{A\oplus \End(\HHH)}^\circlearrowleft$ as fucntors $\EXT_A\simeq (Lie_A^f)^{op}\to \SSS$ since
$U_1(\DD_\infty(R)^{S^1})\simeq \DD_1(h(R))$ for $R\in \EXT_A$ (see Proposition~\ref{commutativediagram}).
Namely, $\mathcal{F}_{A\oplus \End(\HHH)}^\circlearrowleft$
is given by
\[
L\mapsto \Map_{Lie_A^{S^1}}(L^{S^1}, \End^L(\HHH))
\]
as a functor $(Lie_A^{f})^{op}\to \SSS$.
\end{Remark}

\Proof
We first observe that we may assume that
$M$ belongs to $Lie_A^f$.
Let $(-)^{S^1}:Lie_A\to Lie^{S^1}_A$ be the functor given by cotensor with $S^1$.
By (the proof of) Lemma~\ref{separatelemma}, $(-)^{S^1}$ preserves sifted colimits.
For any $M\in Lie_A$,
$M$ is a (sifted) colimit of $(Lie_A^f)_{/M} \to Lie_A\simeq  \mathcal{P}^{st}_{\Sigma}(Lie_A^f)$.
The inclusion $Lie_A\simeq \mathcal{P}^{st}_{\Sigma}(Lie_A^f)\hookrightarrow \Fun(\EXT_A,\SSS)$ carries $M$ to the functor $\mathcal{F}_M:\EXT_A\simeq (Lie_A^f)^{op}\to \SSS$ corepresented by $M$ (cf. Section~\ref{formalstack}). 
An object of $Lie_A^f$ is compact and projective in $Lie_A$, that is, 
the functor
it corepresents preserves sifted colimits.
It follows that $Lie_A\hookrightarrow  \Fun(\EXT_A,\SSS)$ preserves sifted colimits.
Consequently, we may and will assume that $M$ belongs to $Lie_A^f$.

The functor $(-)^{S^1}$ gives rise to 
the adjoint pair
\[
\tau:\Fun((Lie_A)^{op},\wSSS) \rightleftarrows \Fun((Lie_A^{S^1})^{op},\wSSS)
\]
where the
left arrow $\leftarrow$ is the right adjoint functor given by the composition with
$(-)^{S^1}$.
The right arrow $\tau$ is the left adjoint given by left Kan extensions along $(-)^{S^1}:(Lie_A)^{op}\to (Lie_A^{S^1})^{op}$.
Let $\overline{\mathcal{F}}_{E^L}^\circlearrowleft:(Lie_A)^{op} \to \SSS$
denote the image of $\hhhh_{E^L}$ under the right adjoint
(where $\hhhh_{E^L}$ is the functor represented by $E^L$).
Let $\hhhh_{M}:(Lie_A)^{op}\to \SSS$ denote
the functor represented by $M$.
Then by the adjoint pair,
$\Map(\hhhh_{M}, \overline{\mathcal{F}}_{E^L}^\circlearrowleft)\simeq 
\Map(\tau(\hhhh_{M}), \hhhh_{E^L})$.
Note that $(-)^{S^1}:Lie_A\to Lie_A^{S^1}$ commutes with the left adjoint
$\tau$ through the Yoneda embeddings.
Thus, $\tau(\hhhh_{M})$ is equivalent to the functor $\hhhh_{M^{S^1}}$
(reperesented by $M^{S^1}$).
It follows that $\Map(\tau(\hhhh_{M}), \hhhh_{E^L})\simeq \Map_{Lie_A^{S^1}}(M^{S^1},E^L)$.

To prove our assertion, it suffices to
obtain $\Map_{\Fun((Lie_A)^{op},\SSS)}(\hhhh_{M}, \overline{\mathcal{F}}_{E^L}^\circlearrowleft)\simeq \Map_{\Fun(\EXT_A,\SSS)}(\mathcal{F}_{M}, \mathcal{F}^\circlearrowleft_{E^L})$.
The composition with $\EXT_A\stackrel{\DD_\infty}{\simeq }(Lie_A^f)^{op}\hookrightarrow (Lie_A)^{op}$
induces $\Fun((Lie_A)^{op},\SSS)\to \Fun(\EXT_A,\SSS)$.
The image of $\hhhh_M$ under this functor is $\mathcal{F}_M$.
The image of $\overline{\mathcal{F}}_{E^L}^\circlearrowleft$
is $\mathcal{F}_{E^L}^\circlearrowleft$.
Consider the adjoint pair $\gamma:\Fun(\EXT_A,\SSS)\rightleftarrows \Fun(Lie_A^{op},\SSS)$ such that the right adjoint is induced
by the composition with $\EXT_A\simeq (Lie_A^f)^{op}\hookrightarrow Lie_A^{op}$.
It is enough to show that $\gamma(\mathcal{F}_{M})\simeq \hhhh_{M}$ for $M\in Lie_A^f$.
The left adjoint $\gamma$ is the left Kan extension of $Lie_A^f\subset Lie_A\to \Fun(Lie_A^{op},\SSS)$ along $Lie_A^f\simeq (\EXT_A)^{op}\stackrel{Y}{\to} \Fun(\EXT_A,\SSS)$
where $Y$ is the Yoneda embedding.
Thus, $\gamma(\mathcal{F}_{M})\simeq \hhhh_{M}$ when $M\in Lie_A^f$.
The final assertion is straightforward to check from the construction of the equivalence.
\QED

\begin{Proposition}
\label{algebraint2}
The morphism
$\mathfrak{T}_{\CCC}:\mathcal{F}_{\GGG_{\CCC}} \to \mathcal{F}^\circlearrowleft_{A\oplus \End(\HHH)}$ in Proposition~\ref{firstsquare}
is determined by the morphism 
$\widehat{A}_{\CCC}^L:\GG_{\CCC}^{S^1} \to \End^L(\HHH)$
in Construction~\ref{liealgebraaction}
through the equivalence in Lemma~\ref{modularint}
(see also Remark~\ref{modularintrem}).
\end{Proposition}

\Proof
By the definition of $\mathfrak{T}_{\CCC}$ and Proposition~\ref{algebraint1},
$\mathfrak{T}_{\CCC}:\mathcal{F}_{\GGG_{\CCC}} \to \mathcal{F}^\circlearrowleft_{A\oplus \End(\HHH)}$
is given by
\begin{eqnarray*}
\alpha:\mathcal{F}_{\GG_{\CCC}}(-)\simeq \Map_{\Alg_2(\Mod_A)}(\DD_2(-),\HHHH) &\stackrel{h}{\longrightarrow}& \Map_{\Alg_1(\Mod_A^{S^1})}(h\circ \DD_2(-),h(\HHHH)) \\
&\to& \Map_{\Alg_1(\Mod_A^{S^1})}(h\circ \DD_2(-),\End(\HHH)) \\
&\simeq& \mathcal{F}^\circlearrowleft_{A\oplus \End(\HHH)}(-)
\end{eqnarray*}
where the second arrow is induced by $A_{\CCC}:h(\HHHH)=\HH_\bullet(\HHHH/A)\to \End(\HHH)$
(see Remark~\ref{algebraint1}).
This sequence indicates the sequence in $\Fun(\EXT_A,\SSS)$,
and $(-)$ means  the ``argument''.
By the adjoint pair $U_2:Lie_A\rightleftarrows \Alg_2(\Mod_A)$ and
the natural equivalence $U_2\circ\DD_{\infty}\simeq \DD_2$
between functors $(\EXT_A)^{op}\to \Alg_2(\Mod_A)$ (see Proposition~\ref{UEAProp} (1)),
the first equivalence is defined as the composite
$\Map_{Lie}(\DD_\infty(-),\GG_{\CCC})\stackrel{U_2}{\to}\Map_{\Alg_2(\Mod_A)}(U_2(\DD_\infty(-)),U_2(\GG_{\CCC}))\to \Map_{\Alg_2(\Mod_A)}(\DD_2(-), \HHHH)$ where the second arrow is given by the composition with the counit map $U_2(\GG_{\CCC})\to \HHHH$ and $U_2\circ\DD_{\infty}\simeq \DD_2$.
Thus,
$\alpha$ is equivalent to 
\begin{eqnarray*}
\beta:\Map_{\Alg_2(\Mod_A)}(\DD_2(-),\HHHH) &\simeq&\Map_{Lie}(\DD_\infty(-),\GG_{\CCC}) \\
&\stackrel{h\circ U_2}{\longrightarrow}& \Map_{\Alg_1(\Mod_A^{S^1})}(h\circ U_2\circ \DD_\infty(-),h(U_2(\GG_{\CCC}))) \\
&\to& \Map_{\Alg_1(\Mod_A^{S^1})}(h\circ U_2\circ \DD_\infty(-),\End(\HHH)) \\
&\simeq& \Map_{\Alg_1(\Mod_A^{S^1})}(h\circ\DD_2(-),\End(\HHH))
\end{eqnarray*}
where the third arrow is induced by $h(U_2(\GG_{\CCC}))\to h(\HHHH)\to \End(\HHH)$,
and the fourth arrow is induced by $U_2\circ\DD_{\infty}\simeq \DD_2$.
Next we use the natural equivalence
$h\circ U_2=\HH_\bullet(-/A)\circ U_2\stackrel{\sim}{\to}U_1\circ (-)^{S^1}$ 
between functors $Lie_A\to \Alg_{1}^+(\Mod_A^{S^1})$
in Proposition~\ref{eoneetwouniversal}.
Then using this natural equivalence, we see that
$\beta$ is equivalent to
\begin{eqnarray*}
\gamma:\Map_{\Alg_2(\Mod_A)}(\DD_2(-),\HHHH) &\simeq&\Map_{Lie}(\DD_\infty(-),\GG_{\CCC}) \\
&\stackrel{U_1\circ (-)^{S^1}}{\longrightarrow}& \Map_{\Alg_1(\Mod_A^{S^1})}(U_1(\DD_{\infty}(-)^{S^1}),U_1(\GG_{\CCC}^{S^1})) \\
&\to& \Map_{\Alg_1(\Mod_A^{S^1})}(U_1(\DD_{\infty}(-)^{S^1}) ,\End(\HHH)) \\
&\simeq& \Map_{\Alg_1(\Mod_A^{S^1})}(h\circ\DD_2(-),\End(\HHH))
\end{eqnarray*}
where the third arrow is induced by $U_1(\GG_{\CCC}^{S^1})\simeq h(U_2(\GG_{\CCC}))\to h(\HHHH)\to \End(\HHH)$,
and the fourth arrow is induced by the natural equivalence $U_1\circ (-)^{S^1}\simeq h\circ U_2\circ\DD_{\infty}\simeq h\circ \DD_2$
between functors $\EXT_A\to \Alg_1(\Mod_A^{S^1})$ (see Proposition~\ref{commutativediagram}).
Note that $\gamma$ is equivalent to
$\Map_{Lie_A}(\DD_\infty(-),\GG_{\CCC})\to \Map_{Lie_A^{S^1}}(\DD_\infty(-)^{S^1},\GG_{\CCC}^{S^1})\to \Map_{Lie_A^{S^1}}(\DD_\infty(-)^{S^1},\End^L(\HHH))$
where the second functor is induced by $\widehat{A}_{\CCC}^L$.
Now our assertion follows from Lemma~\ref{modularint}.
\QED

\subsection{}
\label{rightsquaresection}

We will construct the right square in Theorem~\ref{intromain2}
(see Proposition~\ref{secondsquare2}, Remark~\ref{modulimeaning2}).

The functors $h$ and $i$ in Lemma~\ref{commutativemodulitransform2}
are extended to functors $h^+,i^+:\CAlg_A^+\to \CAlg^+(\Mod_A^{S^1})$,
and $\sigma:h\to i$ is extended to $\sigma^+:h^+\to i^+$ in the natural way. 
We denote by
$h|_{\EXT_A}\to i|_{\EXT_A}$
a natural transformation between functors $\EXT_A\to \Alg_1^+(\Mod_A^{S^1})$, which is 
obtained from $h^+\to i^+$
by the compositions with $\EXT_A\to \CAlg_A^+$ and $\CAlg^+(\Mod_A^{S^1})\stackrel{\textup{forget}}{\to} \Alg_1^+(\Mod_A^{S^1})$.

\begin{Construction}
\label{secondsquare}
We construct the square diagram in $\Fun(\EXT_A,\SSS)$.
Recall the diagram such that the vertical arrows are right fibrations 
\[
\xymatrix{
(\mathcal{D}ef_{\HHH}^\circlearrowleft)^{op} \ar[r]^(0.4){\DD_{\HHH}^\circlearrowleft} \ar[d] & \LMod^+(\Mod_A^{S^1})_{\HHH} \ar[d]  \\
\Alg_1^+(\Mod_A^{S^1})^{op} \ar[r]_{\DD_1^{S^1}} & \Alg_1^+(\Mod_A^{S^1})
}
\]
from the second cube in Remark~\ref{detailedreduction}.
Taking into account the corresponding functors $\Alg_1^+(\Mod_A^{S^1})\to \SSS$ and Lemma~\ref{endorepresentable}, we have a natural transformation
 $\Def^{\circlearrowleft}_{\HHH}  \stackrel{J^{\circlearrowleft}_{\HHH}}{\longrightarrow} \mathfrak{h}_{A\oplus\End(\HHH)}  \circ  \DD_{1}^{S^1}$
(compare Construction~\ref{fromcube}).
The composition with $h|_{\EXT_A}\to i|_{\EXT_A}$
yields the diagram in $\Fun(\EXT_A,\SSS)$
\[
\xymatrix{
\Def^{\circlearrowleft}_{\HHH} \circ h|_{\EXT_A} \ar[r] \ar[d] & \Def^{\circlearrowleft}(\HHH) \circ i|_{\EXT_A} \ar[d]  \\
\mathfrak{h}_{A\oplus\End(\HHH)}  \circ  \DD_{1}^{S^1}\circ h|_{\EXT_A} \ar[r] & \mathfrak{h}_{A\oplus\End(\HHH)}  \circ  \DD_{1}^{S^1}\circ i|_{\EXT_A}
}
\]
which commutes up to canonical homotopy.
Note that $\Def^{\circlearrowleft}(\HHH)=\Def^{\etwo,\circlearrowleft}(\HHH)|_{\EXT_A}=\Def^{\circlearrowleft}_{\HHH} \circ h|_{\EXT_A}$.
Note also that $\mathcal{F}^{\circlearrowleft}_{A\oplus \End(\HHH)}=\mathfrak{h}_{A\oplus\End(\HHH)}  \circ  \DD_{1}^{S^1}\circ h|_{\EXT_A}$.
The left vertical arrow is $J^\circlearrowleft_{\HHH}$
in Proposition~\ref{firstsquare}.
\end{Construction}

\begin{Remark}
The functor
$\Def^{\circlearrowleft}(\HHH) \circ i|_{\EXT_A}:\EXT_A\to \SSS$
is informally given by
\[
\EXT_A\ni [R\to A] \mapsto \LMod_{i(R)}(\Mod_A^{S^1})^{\simeq}\times_{(\Mod^{S^1}_A)^{\simeq}}\{\HHH\}\in \SSS
\]
where $\LMod_{i(R)}(\Mod_A^{S^1}) \to \Mod_A$ is the base change induced by $R\to A$. The right-hand side is the space of $S^1$-equivariant deformations
of $\HH_\bullet(\CCC/A)$ to $R$.
Namely, it describes the $S^1$-equivariant deformation problem.
\end{Remark}

\begin{Notation}
We write $\Def^{S^1}(\HHH)$ for $\Def^{\circlearrowleft}(\HH_\bullet(\CCC/A)) \circ i|_{\EXT_A}$.
\end{Notation}

\begin{Lemma}
\label{s1endlie}
Let $\End^L(\HHH)$ denote the dg Lie algebra
endowed with $S^1$-action
(i.e., an object of $Lie_A^{S^1}\simeq \Fun(BS^1,Lie_A)$)
associated to $\End(\HHH)\in \Alg_1(\Mod_A^{S^1})$
(that is, $\End^L(\HHH)=res_{\eone/Lie}(\End(\HHH))$).
Let $\End^L(\HHH)^{S^1}$ be the $S^1$-invariants
(homotopy fixed points).
Then there exists a canonical equivalence
$\mathfrak{h}_{A\oplus\End(\HHH)}  \circ  \DD_{1}^{S^1}\circ i|_{\EXT_A}\simeq \mathcal{F}_{\End^L(\HHH)^{S^1}}$
where
$\mathcal{F}_{\End^L(\HHH)^{S^1}}$ is the formal stack associated to $\End^L(\HHH)^{S^1}$ (cf. Section~\ref{formalstack}).
\end{Lemma}

\Proof
Observe that $\mathfrak{h}_{A\oplus\End(\HHH)}  \circ  \DD_{1}^{S^1}\circ i|_{\EXT_A}$
is given by
\[
[R\to A] \mapsto \Map_{\Alg_{1}(\Mod_A^{S^1})}(\DD_{1}(R),\End(\HHH)).
\]
Here we think of $\DD_{1}(R)$ as an unaugmented $\eone$-algebra in $\Mod_A$ endowed with
the trivial $S^1$-action. 
Since $\DD_1(R)\simeq U_1(\DD_\infty(R))$ for $R \in \EXT_A$ (see Proposition~\ref{UEAProp} (1)),
it follows that $\mathfrak{h}_{A\oplus\End(\HHH)}  \circ  \DD_{1}^{S^1}\circ i|_{\EXT_A}:\EXT_A\to \SSS$
is given by
\begin{eqnarray*}
[R\to A]   &\mapsto&  \Map_{\Alg_{1}(\Mod_A^{S^1})}(\DD_{1}(R),\End(\HHH)) \\
&\simeq& \Map_{\Alg_{1}(\Mod_A)}(\DD_{1}(R),\End(\HHH)^{S^1}) \\
   &\simeq& \Map_{\Alg_{1}(\Mod_A)}(U_1(\DD_{\infty}(R)),\End(\HHH)^{S^1})  \\
   &\simeq& \Map_{Lie_A}(\DD_{\infty}(R),\End^L(\HHH)^{S^1}).
\end{eqnarray*}
It follows that $\mathfrak{h}_{A\oplus\End(\HHH)}  \circ  \DD_{1}^{S^1}\circ i|_{\EXT_A}\simeq \mathcal{F}_{\End^L(\HHH)^{S^1}}$.
\QED

\begin{Proposition}
\label{secondsquare2}
There exists a diagram in $\Fun(\EXT_A,\SSS)$
\[
\xymatrix{
\Def^\circlearrowleft(\HH_\bullet(\CCC/A)) \ar[r]^{N_{\CCC}} \ar[d]^{J^\circlearrowleft_{\HH_\bullet(\CCC/A)}} & \Def^{S^1}(\HH_\bullet(\CCC/A)) \ar[d]^{J_{\HH_\bullet(\CCC/A)}^{S^1}} \\
\mathcal{F}^\circlearrowleft_{A\oplus \End(\HH_\bullet(\CCC/A))} \ar[r] & \mathcal{F}_{\End^L(\HH_\bullet(\CCC/A))^{S^1}}
}
\]
which commutes up to canonical homotopy.
\end{Proposition}

\Proof
Consider the second square diagram in Construction~\ref{secondsquare}.
By Lemma~\ref{s1endlie} and definitions, we can
interpret the square diagram
in Construction~\ref{secondsquare}
as the desired one.
\QED

\begin{Remark}
\label{modulimeaning2}
We consider the moduli-theoretic meaning of $N_{\CCC}$.
Recall that the space
$\Def^{\circlearrowleft}(\HHH)(R)$
parametrizes cyclic deformations of $\HH_\bullet(\CCC/A)$
to $R\otimes_AS^1$, that is,
$\LMod_{\HH_\bullet(R/A)}(\Mod_A^{S^1})^\simeq \times_{(\Mod^{S^1}_A)^{\simeq}}\{\HHH\}$. Let
\[
(\mathcal{H}'\in \RMod_{R\otimes_AS^1}(\Mod_A^{S^1}),\ R\otimes_AS^1=\HH_\bullet(R/A)\to A,\ \HHH\simeq \mathcal{H}'\otimes_{R\otimes_AS^1}A).
\]
be an object of $\Def^{\circlearrowleft}(\HHH)(R)$.
The image under $N_{\CCC}$
is
\[
(\mathcal{H}'\otimes_{R\otimes_AS^1}R\in \RMod_{R}(\Mod_A^{S^1}),\ R\to A,\ \HHH\simeq (\mathcal{H}'\otimes_{R\otimes_AS^1}R)\otimes_{R}A).
\]
Namely, $N_{\CCC}$ is given by the base change along $R\otimes_AS^1\to R$ for each $R\in \EXT_A$. 
\end{Remark}

\begin{Remark}
\label{liediagonal}
By the definition (see Construction~\ref{secondsquare}),
$\mathcal{F}^\circlearrowleft_{A\oplus \End(\HHH)} \to \mathcal{F}_{\End^L(\HHH)^{S^1}}$
is given by
\begin{eqnarray*}
\mathcal{F}^\circlearrowleft_{A\oplus \End(\HHH)}(-) &\simeq& \Map_{Lie_A^{S^1}}(\DD_{\infty}(-)^{S^1},\End^L(\HHH)) \\ &\to& \Map_{Lie_A^{S^1}}(\DD_{\infty}(-),\End^L(\HHH))  \\ &\simeq& \Map_{Lie_A}(\DD_{\infty}(-),\End^L(\HHH)^{S^1})= \mathcal{F}_{\End^L(\HHH)^{S^1}}(-)
\end{eqnarray*}
where the second arrow is induced by the diagonal
map $\DD_{\infty}(-)\to \DD_{\infty}(-)^{S^1}$ (see also Lemma~\ref{modularint} for the first equivalence). 
\end{Remark}

\subsection{}
\label{finalsubsection}

\begin{Theorem}
\label{supermain}
There exists a diagram in $\Fun(\EXT_A,\SSS)$
\[
\xymatrix{
\Def_{\CCC} \ar[r]^(0.3){M_{\CCC}^{\circlearrowleft}} \ar[d]^{J_{\CCC}} &  \Def^\circlearrowleft(\HH_\bullet(\CCC/A)) \ar[r]^{N_{\CCC}} \ar[d]^{J^\circlearrowleft_{\HH_\bullet(\CCC/A)}} & \Def^{S^1} (\HH_\bullet(\CCC/A)) \ar[d]^{J_{\HH_\bullet(\CCC/A)}^{S^1}} \\
\mathcal{F}_{\GGG_{\CCC}} \ar[r]^(0.3){\mathfrak{T}_{\CCC}}  & \mathcal{F}^\circlearrowleft_{A\oplus \End(\HH_\bullet(\CCC/A))} \ar[r] & \mathcal{F}_{\End^L(\HH_\bullet(\CCC/A))^{S^1}}
}
\]
which commutes up to canonical homotopy.
\end{Theorem}

\Proof
Combine the diagram in Proposition~\ref{firstsquare} and the diagram in Propoition~\ref{secondsquare2}.
\QED

\begin{Proposition}
\label{algebraint3}
The composition $\mathcal{F}_{\GGG_{\CCC}} \stackrel{\mathfrak{T}_{\CCC}}{\to} \mathcal{F}^\circlearrowleft_{A\oplus \End(\HH_\bullet(\CCC/A))} \to \mathcal{F}_{\End^L(\HH_\bullet(\CCC/A))^{S^1}}$
corresponds to $A_{\CCC}^L:\GG_{\CCC}\to \End^L(\HH_\bullet(\CCC/A))^{S^1}$
(cf. Construction~\ref{liealgebraaction})
through $\FST_A\simeq Lie_A$. 
\end{Proposition}

\Proof
We write $\HHH$ for $\HH_\bullet(\CCC/A)$.
Consider the square diagram
\[
\xymatrix{
\Map_{Lie^{S^1}_A}(\GG_{\CCC}^{S^1},\End^L(\HHH)) \ar[r]^\simeq \ar[d] & \Map(\mathcal{F}_{\GG_{\CCC}},\mathcal{F}_{A\oplus \End(\HHH)}^\circlearrowleft)  \ar[d] \\
\Map_{Lie^{S^1}_A}(\GG_{\CCC},\End^L(\HHH)) \ar[r]^\simeq &  \Map(\mathcal{F}_{\GG_{\CCC}},\mathcal{F}_{\End^L(\HHH)^{S^1}}).
}
\]
The upper horizontal equivalence comes from
in Lemma~\ref{modularint}.
The lower horizontal equivalence comes from $Lie_A\simeq \FST_A$
and $\Map_{Lie^{S^1}_A}(\GG_{\CCC},\End^L(\HHH))\simeq \Map_{Lie_A}(\GG_{\CCC},\End^L(\HHH)^{S^1})$.
The left vertical arrow is induced by the composition with the diagonal map $\GG_{\CCC}\to \GG_{\CCC}^{S^1}$.
The right vertical arrow is is induced by the composition with $\mathcal{F}_{A\oplus \End(\HHH)}^\circlearrowleft \to \mathcal{F}_{\End^L(\HHH)^{S^1}}$.
Unfolding the construction of the upper horizontal equivalence
in Lemma~\ref{modularint} (see also Remark~\ref{modularintrem}),
we observe that this square commutes up to homotopy. 
According to Proposition~\ref{algebraint2}, $\mathfrak{T}_\CCC$
corresponds to $\widehat{A}^L_{\CCC}$ in Construction~\ref{liealgebraaction}
(via upper horizontal equivalence)
so that $\mathcal{F}_{\GGG_{\CCC}} \stackrel{\mathfrak{T}_{\CCC}}{\to} \mathcal{F}^\circlearrowleft_{A\oplus \End(\HHH)} \to \mathcal{F}_{\End^L(\HHH)^{S^1}}$
corresponds to $A_{\CCC}^L:\GG_{\CCC}\to \End^L(\HHH)^{S^1}$, as desired.
\QED

\begin{Remark}
\label{refinedremark}
Let us observe the diagram in Theorem~\ref{supermain}.
Roughly, the upper row has a moduli-theoretic description
(Remark~\ref{modulimeaning2} and Remark~\ref{modulimeaning2}),
while the lower row is defined in terms of (Lie) algebras
(Proposition~\ref{algebraint1} and Proposition~\ref{algebraint2}). 
Suppose that we are given a deformation $\CCC_{R}\in \Def(\CCC_R)$ of $\CCC$ to $R\in \EXT_A$
(informally, we think of $\CCC_R$ as an object of $\ST_R$).
Taking into account the modular interpretations of $\mathcal{F}_{A\oplus \End(\HH_\bullet(\CCC/A))}^{\circlearrowleft}$ and $\mathcal{F}_{\End^L(\HH_\bullet(\CCC/A))^{S^1}}$ (Lemma~\ref{modularint}), 
the images of $\CCC_R$ under the above diagram can informally be depicted as follows:
\[
\xymatrix{
\CCC_R\ar@{.>}[r] \ar@{.>}[d] & \overset{\text{cyclic deformation}}{\HH_\bullet(\CCC_R/A)} \ar@{.>}[d]  \ar@{.>}[r]  &  \overset{\text{equivariant deformation}}{\HH_\bullet(\CCC_R/R)} \ar@{.>}[d]   \\
\{\DD_\infty(R) \stackrel{p}{\to} \GG_{\CCC}\} \ar@{.>}[r] &  \{ \DD_\infty(R)^{S^1}\stackrel{\widehat{A}_{\CCC}^L \circ p^{S^1}}{\longrightarrow} \End^L(\HH_\bullet(\CCC/A))\} \ar@{.>}[r]   &  \{\DD_\infty(R)\stackrel{\widehat{A}_{\CCC}^L\circ p^{S^1}\circ \iota}{\longrightarrow} \End^L(\HH_\bullet(\CCC/A)) \}
}
\]
Here 
$p:\DD_\infty(R) \to \GG_{\CCC}$ indicates the image of the deformation under $J_{\CCC}$, that is an object of the space
$\FF_{\GG_{\CCC}}(R)$.
The middle and right items on the lower row is defined by compositions of the sequence
of morphisms 
$\DD_\infty(R)\stackrel{\iota}{\to} \DD_\infty(R)^{S^1} \stackrel{p^{S^1}}{\to} \GG^{S^1}\stackrel{\widehat{A}^L_{\CCC}}{\to} \End^L(\HH_\bullet(\CCC/A))$ in $Lie_A^{S^1}$
where $\iota$ is the diagonal map.
We also use $\HH_\bullet(\CCC_R/R)\otimes_{\HH_\bullet(R/A)}R\simeq \HH_\bullet(\CCC_R/R)$. 
\end{Remark}

{\it Acknowledgements.}
The author would like to thank Takuo Matsuoka for valuable comments and many discussions since the summer of 2017.
He would like to thank everyone who provided constructive feedback
at his talks about the main content of this paper.
This work is supported by JSPS KAKENHI grant.

\end{document}